\documentclass[11pt]{amsart}

\usepackage{hyperref}
\hypersetup{colorlinks=true,linkcolor=blue,citecolor=blue,urlcolor=blue}
\usepackage{amssymb, amsmath, amsthm, amsfonts}
\usepackage{enumerate}
\usepackage{xcolor}
\usepackage[mathscr]{eucal}
\usepackage{tikz}
\usepackage{graphicx}
\usepackage{enumitem}

\makeatletter
\@ifundefined{subjclassname@2020}{%
  \@namedef{subjclassname@2020}{\textup{2020} Mathematics Subject Classification}%
}{}
\makeatother

\makeatletter
\renewcommand{\@captionfont}{\normalfont\footnotesize}
\makeatother
\newtheorem{thm}{Theorem}[section]

\newtheorem{prop}[thm]{Proposition}
\newtheorem{cor}[thm]{Corollary}
\newtheorem{lem}[thm]{Lemma}
\theoremstyle{remark}
\newtheorem{rem}[thm]{Remark}
\numberwithin{equation}{section}

\newcommand\supp{\operatorname{supp}}

\newcommand{\rd}{\mathrm d}
\newcommand{\R}{\mathbb R}
\newcommand{\Schatten}{\mathfrak S}

\newcommand{\e}{\mathrm{e}}
\newcommand{\Fmu}{\langle\mu\rangle_\alpha}
\newcommand{\Tr}{\operatorname{Tr}}
\newcommand{\rank}{\operatorname{rank}}
\newcommand{\ind}{\mathbf 1}

\date{\today}

\begin{document}

\keywords{Fractional Schr\"odinger equation, wave equation, Klein--Gordon equation, Strichartz estimates, orthonormal systems, fractal time}
\subjclass[2020]{35Q41 (primary); 35B45, 42B37, 28A80}

\title[Orthonormal Strichartz estimates over fractal time]{Orthonormal Strichartz estimates over fractal time}

\author[G. Feng]{Guoxia Feng} \address[G. Feng]{School of Mathematical Sciences, Xiamen University, Xiamen 361005, China} \email{fengguoxia@stu.xmu.edu.cn}

\author[J.~B. Lee]{Jin Bong Lee} \address[J.~B. Lee]{Research Institute of Mathematics, Seoul National University, Seoul 08826, Republic of Korea} \email{jinblee@snu.ac.kr}

\author[S. Lee]{Sanghyuk Lee} \address[S. Lee]{Department of Mathematical Sciences and RIM, Seoul National University, Seoul 08826, Republic of  Korea} \email{shklee@snu.ac.kr}

\begin{abstract}
We study Strichartz estimates for orthonormal families of initial data
when time is integrated against an $\alpha$-dimensional Frostman measure
$\mu$, $0<\alpha\le1$, with coefficients measured in $\ell^\beta$. For the
fractional Schr\"odinger, wave and Klein--Gordon equations, we obtain
orthogonality gains in the admissible region and necessary conditions on
$\beta$. For the fractional Schr\"odinger equation the exponents are sharp
up to an endpoint; the optimality of the threshold $\beta=q/2$ requires
$\mu$ to be Ahlfors regular. For the Schr\"odinger equation in dimensions
$d\ge2$, we prove a Lorentz endpoint estimate and, if $d>2\alpha$, the
critical exponent $\beta=q/2$ at every pair $2\le q,r<\infty$ with
\[\frac d2(\frac12-\frac1r)\ge\frac\alpha q, \quad r>\frac{dq}{d-\alpha};\]
on $r=\infty$ we obtain a restricted weak-type estimate. Since the
constants depend on $\mu$ only through its Frostman constant, Lebesgue
measure yields global-in-time critical estimates in this range when
$d\ge3$. The proofs combine Frostman fractional integration, a spectral
near--far argument at the endpoint, and asymmetric Hilbert--Schmidt bounds
with orthogonal frequency-block summation. We also obtain further
estimates in the low-dimensional cases $d=1$ (fractional Schr\"odinger)
and $d=2$ (wave and Klein--Gordon) with $\alpha>1/2$, and formulate the
consequences for densities of freely evolving states.
\end{abstract}

\maketitle

\section{Introduction and main results}

Let $\phi$ be a real-valued dispersion relation and denote by $U_\phi$
the solution operator of the linear dispersive equation
\begin{equation*}
 \left\{
 \begin{aligned}
 &i\partial_t u+\phi(D) u=0,\qquad x\in\mathbb{R}^d,\ t\in\mathbb{R},\\
 &u(x,0)=f(x).
 \end{aligned}
 \right.
\end{equation*}
We use the unitary Fourier transform
$\widehat f(\xi)=(2\pi)^{-d/2}\int e^{-ix\cdot\xi}f(x)\,\rd x$.
The solution is
\begin{equation*}
U_\phi(t)f(x)=e^{it\phi(D)}f(x)=(2\pi)^{-d/2}\int_{\mathbb{R}^d}e^{i(x\cdot \xi+t\phi(\xi))}\hat{f}(\xi)\,\mathrm d\xi.
\end{equation*}
Fixed dimensional Fourier constants are absorbed into kernel estimates below; the orthonormality identities use the unitary convention. The classical Strichartz estimates take the form
\begin{equation}\label{class}
    \|U_\phi f\|_{L_t^q(\mathbb{R}; L_x^r(\mathbb{R}^d))}\lesssim\|f\|_{\mathcal{H}},
\end{equation}
where the initial data $f$ belongs to an appropriate Hilbert space $\mathcal{H}$, typically the homogeneous Sobolev space $\dot{H}^s(\mathbb{R}^d)$ or the inhomogeneous Sobolev space ${H}^s(\mathbb{R}^d)$. The usual admissibility condition describes the exponent range, together with the appropriate regularity and endpoint exclusions for each flow. More precisely, for a parameter $\sigma>0$, a pair $(q,r) \in [2,\infty]\times[2,\infty]$ is called  $\sigma$-admissible if
\begin{equation*}
    \sigma\Big(\frac{1}{2}-\frac{1}{r}\Big)\geq\frac{1}{q},
\end{equation*}
and sharp $\sigma$-admissible if equality holds.

Typical examples are given by $\phi(\xi) = |\xi|^a, |\xi|, \sqrt{1+|\xi|^2}$ with $a\in \mathbb R_{>0}\setminus\{1\}$, which correspond to the fractional Schr\"odinger, the wave, and the Klein--Gordon equations, respectively.
The associated Strichartz estimates \eqref{class} with $\mathcal H = H^s, \dot{H}^s$ have been extensively studied in terms of $\sigma$-admissible pairs $(q, r)$ and the optimal regularity $s$ of the initial data.
We refer the reader to \cite{RS1977,GV1995,KT1998,CKS2016,GW2014,HS2015}.

Frank, Lewin, Lieb and Seiringer \cite{FLLS2014} and Frank and Sabin \cite{FraSa2017} initiated the study of generalizations of \eqref{class} to orthonormal families of initial data, that is, estimates of the form
 \begin{align}\label{orthonormal-Strichartz}
	\Bigg\| \sum_{j\in J} n_j |U_\phi f_j|^2 \Bigg\|_{L_t^{\frac{q}{2}}(\mathbb{R}; L_x^{\frac{r}{2}}(\mathbb{R}^d))} \lesssim \|\{n_j\}_{j\in J}\|_{\ell^\beta},
\end{align}
where $\{f_j\}_j$ is an orthonormal family in $\dot{H}^s(\mathbb R^d)$ for the Schr\"odinger and the wave equations and in $H^s(\mathbb R^d)$ for the Klein--Gordon equation.
Combining \eqref{class} with the triangle inequality gives \eqref{orthonormal-Strichartz} with $\beta=1$, so the main problem is to obtain \eqref{orthonormal-Strichartz} with $\beta$ as large as possible.
For Schr\"odinger evolution, Frank--Lewin--Lieb--Seiringer \cite{FLLS2014} proved the mixed-norm result for $2\le r\le2(d+2)/d$, and Frank--Sabin \cite{FraSa2017} extended it to $2\le r<2(d+1)/(d-1)$, read as $r<\infty$ when $d=1$. In this range, \eqref{orthonormal-Strichartz} holds with the optimal exponent $\beta = \frac{2r}{r+2}$ for sharp admissible pairs $(q,r)$. These estimates were extended, in particular to the wave and Klein--Gordon equations and to frequency localized and endpoint settings, by Bez--Hong--Lee--Nakamura--Sawano \cite{BHLNS2019}, Bez--Lee--Nakamura \cite{BLN2021}, and Bez--Kinoshita--Shiraki \cite{BKS2024}; when $\sigma>1/2$, the corresponding range for $\sigma$-admissible pairs is $2\le r<\frac{2(2\sigma+1)}{2\sigma-1}$.
For $\sigma>1/2$ at $r=\frac{2(2\sigma+1)}{2\sigma-1}$, in the abstract energy--dispersion setting, Hadama and Toyoshima \cite{HT2026} proved that \eqref{orthonormal-Strichartz} holds with the $\ell^{\beta,1}$-norm on the right-hand side. For other related results, we refer to \cite{MS2025, MSwain2022, FS2024, KPMM2024, GMS2024, FMSW2026}.

For Lebesgue time, Hadama--Toyoshima \cite{HT2026} also obtain the
critical strong estimate on the sharp segment beyond the Lorentz
endpoint. Gong \cite[Theorem 1.3]{Gong2026} proves the critical estimate
at $q=4$ with finite spatial exponent. In \cite{FSW2025}, a related critical interior result is obtained
via weak mixed-norm interpolation. We do not use the general interpolation identity of \cite{FSW2025}:
Section~\ref{sec-quadratic-limitations} gives a counterexample to the needed
function-space embedding.

In this paper, we establish orthonormal Strichartz estimates over $\mathbb{R}^d\times E$, where the time variable is restricted to a set $E$ of fractal dimension and the Lebesgue measure $\mathrm dt$ is replaced by a measure $\mu$ supported on $E$. We prove estimates of the form
\begin{align}\label{goal-esti}
	\Bigg\| \sum_{j\in J} n_j |U_\phi f_j|^2 \Bigg\|_{L_t^{\frac{q}{2}}(\mathrm d\mu; L_x^{\frac{r}{2}})} \lesssim \|\{n_j\}_{j\in J}\|_{\ell^\beta}.
\end{align}
This is motivated by the recent work of Lee--Lee--Roncal \cite{LLR2025}, where the Strichartz estimates \eqref{class} over $\mathbb{R}^d\times E$ were studied.

Building on the scalar estimates of \cite{LLR2025}, we obtain orthonormal
Strichartz estimates over fractal time for the fractional Schr\"odinger,
wave and Klein--Gordon equations, in which the dimension of the time
measure enters the admissibility condition and the regularity. For the
fractional Schr\"odinger equation, the coefficient exponents are optimal.
For the Schr\"odinger equation, we prove the critical estimates at the
endpoint and, in higher dimensions, in the whole range of finite exponents
to the left of the endpoint line. Since the constants depend on the
measure only through its Frostman constant, the case of Lebesgue measure
yields global-in-time estimates; see Remark~\ref{rem-critical-scope}.

\subsection{Strichartz estimates over fractal time}\label{sec-fractal-time}
Let $E\subset\mathbb R$ be bounded, and let $\mu$ be a nonzero finite positive Borel measure supported on $E$. The sufficiency estimates are trivial for the zero measure; nonzero mass is needed for the failure statements below. Throughout this paper, we assume that $\mu$ satisfies the Frostman condition
\begin{align}\label{frac_meas}
	\Fmu:=\sup_{x_0\in\mathbb R,\,\rho>0}\rho^{-\alpha}\mu\left( (x_0 - \rho, x_0+\rho) \right) <\infty,
\end{align}
and we call $\mu$ an $\alpha$-dimensional Frostman measure, $\alpha\in(0,1]$.
The Lebesgue measure restricted to a bounded interval is a basic example with $\alpha=1$.
For some of the necessary conditions we shall also use the lower bound
\begin{align}\label{lower-regular}
	\inf_{t\in\supp\mu,\,0<\rho\le1}\rho^{-\alpha}\mu\left( (t - \rho, t+\rho) \right)>0.
\end{align}
Measures satisfying both \eqref{frac_meas} and \eqref{lower-regular} ($\alpha$-Ahlfors regular measures) exist for every $\alpha\in(0,1]$; see, e.g., \cite[Lemma~2.14]{LLR2025}, \cite{LRZZ}.

In this setting, a pair $(q,r)\in[2,\infty]^2$ is said to be $\sigma$-admissible if it satisfies
\begin{align}\label{frac_adm}
	\sigma \Big(\frac12 - \frac1r\Big) \ge \frac\alpha q,
\end{align}
and sharp $\sigma$-admissible when equality holds in \eqref{frac_adm}.
 If the inequality in \eqref{frac_adm} is strict, then $(q,r)$ is called non-sharp $\sigma$-admissible. Let $s_a(q,r)$ be given by
\begin{align}\label{frac_reg}
s_{a}(q,r)=\frac{d}{2}-\frac{d}{r}-\frac{a}{q}.
\end{align}
For finite spatial exponents, the following scalar consequences of \cite{LLR2025} will be used. Strict $r=\infty$ inputs are proved separately in Lemma~\ref{lem-scalar-OD-completion}.
\begin{thm}[Theorem~1.4, \cite{LLR2025}]\label{Class-Schro}
	Let $\alpha \in (0,1]$, $a \in \mathbb R_{>0}\setminus\{1\}$, and $2\le q\le\infty$, $2\le r<\infty$ with $\frac dr+\frac aq<d$. Suppose  $(q,r)$ is $\frac{d}{2}$-admissible and $\mu$ satisfies \eqref{frac_meas}. Then
	\begin{align}\label{class-frac}
		\| e^{it(-\Delta)^\frac{a}{2}} f \|_{L_t^q(\mathrm d\mu; L_x^r(\mathbb{R}^d))} \lesssim  \|f\|_{\dot{H}^s(\mathbb R^d)}
	\end{align}
	 holds with $s= s_{a\alpha}(q,r)$.
\end{thm}
The condition $\frac dr+\frac aq<d$ ensures that the regularity in
Theorem~\ref{Class-Schro} satisfies $s=s_{a\alpha}(q,r)>-d/2$, since
$\alpha\le1$. For $s\le- d/2$, the space $\dot H^s(\mathbb R^d)$ is not
naturally realized as a space of tempered distributions in which the
Schwartz functions are dense; see \cite[p.~4]{LLR2025}. The condition is
not used in the frequency-localized estimates
\cite[Proposition~2.9]{LLR2025}, from which Lemma~\ref{lem-scalar-annular-synthesis}
below obtains the scalar estimates at finite $r$ needed in our interpolation
arguments without this condition.

\begin{thm}[Remark~2.3, \cite{LLR2025}]\label{Class-wave} Let $d\geq2$ and $2\le q\le\infty$, $2\le r<\infty$. Suppose $(q,r)$ is $\frac{d-1}{2}$-admissible and $\mu$ satisfies \eqref{frac_meas}. Then
\begin{align}\label{class-wave}
		\| e^{it\sqrt{-\Delta}} f \|_{L_t^q(\mathrm d\mu; L_x^r(\mathbb{R}^d))} \lesssim  \|f\|_{\dot{H}^s(\mathbb R^d)}
	\end{align}
     holds with $s= s_{\alpha}(q,r)$.
\end{thm}

\subsection{Strichartz estimates for orthonormal families}
It is natural to ask whether Theorems~\ref{Class-Schro} and~\ref{Class-wave} extend to orthonormal families. If \(\{f_j\}_{j\in J}\) is an orthonormal family and
\(
\gamma_0=\sum_{j\in J}n_j|f_j\rangle\langle f_j|
\)
is the associated density matrix, then the density of its free evolution is
\[
\rho_{U_\phi(t)\gamma_0U_\phi(t)^*}(x)
=\sum_{j\in J}n_j|U_\phi(t)f_j(x)|^2.
\]
Thus estimates of the form \eqref{goal-esti} bound the densities of quantum states sampled at times in $E$, and the exponent $\beta>1$ measures the gain from orthogonality.

We briefly describe the proofs. By the duality principle of Frank and
Sabin (Proposition~\ref{prop_dual}), \eqref{goal-esti} is equivalent to a
Schatten bound for a weighted $TT^*$ operator. Since $\mu$ is singular,
convolution and Fourier arguments in time are unavailable; all temporal
estimates come from Frostman fractional integration
(Lemma~\ref{lem_HLS}). For the general dispersions, we interpolate an
analytic family containing a temporal weight, chosen so that the time
kernels at both the operator-norm and the Hilbert--Schmidt endpoints are
integrable against $\mu$ (Section~\ref{sec-sharp}). For the critical estimates for the
Schr\"odinger equation, we adapt the spectral near--far argument of
Hadama--Toyoshima \cite{HT2026} at the sharp endpoint and the direct $q=4$
Hilbert--Schmidt argument of Gong \cite{Gong2026}; for $q>4$, asymmetric
Sobolev orders give decay between initial frequency annuli, and
orthogonal block sums yield the critical strong estimate
(Section~\ref{sec-critical}).

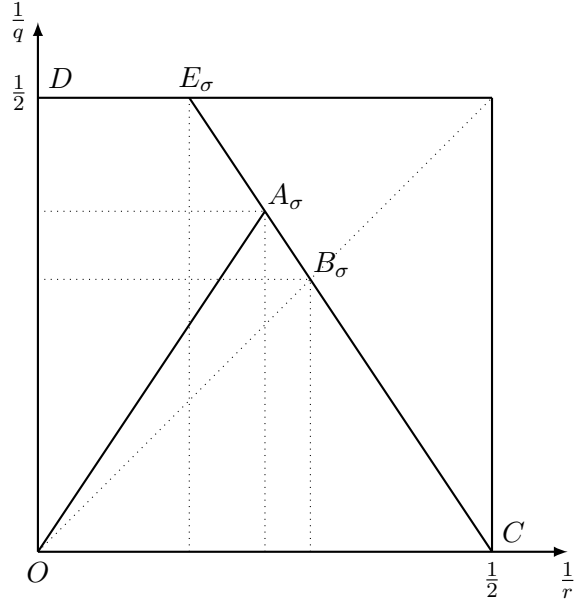
\begin{figure}[t]
\centering
\begin{tikzpicture}
\draw[thick,->, >=latex] (0,0) -- (0,7) node[left] {\(\frac{1}{q}\)};
\draw[thick,->, >=latex] (0,0) -- (7,0) node[below] {\(\frac{1}{r}\)};
\node at (0,5.7) [above left]{\(\frac{1}{2}\)};
\node at (6,0) [below]{\(\frac{1}{2}\)};
\node at (0,0) [below]{\(O\)};
\node at (0,6) [above right]{\(D\)};
\node at (2.1,5.95) [above]{\(E_\sigma\)};\node at (6,0) [above right]{\(C\)};
\node at (3.5, 3.5) [above right]{\(B_\sigma\)};
\node at (2.9, 4.4) [above right]{\(A_\sigma\)};
\draw[thick] (6,0) -- (6,6);
\draw[thick] (6,0) -- (2,6);
\draw[thick] (0,6) -- (6,6);
\draw[thick] (0,0) -- (3,9/2);
\draw[dotted] (2,0)--(2,6);
\draw[dotted] (3,0)--(3,9/2);
\draw[dotted] (0,9/2)--(3,9/2);
\draw[dotted] (18/5,0)--(18/5,18/5);
\draw[dotted] (0,18/5)--(18/5,18/5);
\draw[dotted] (0,0)--(6,6);
\end{tikzpicture}
\caption{The admissible exponent region when $\sigma>\alpha$ (schematic).}\label{fig-admissible}
\end{figure}

Describing the exponent regions in coordinates $(1/r,1/q)$, we introduce the following points in $[0,\frac{1}{2}]\times[0,\frac{1}{2}]$ (see Figure~\ref{fig-admissible}). Let $O=(0,0)$, $C=(\frac12,0)$, $D=(0,\frac12)$. For $\sigma\geq\frac{\alpha}{2}$, the points $A_\sigma$ and $B_\sigma$ are given by
\begin{equation*}
 A_\sigma = \left(\frac{2\sigma-\alpha}{2(2\sigma+\alpha)},\frac{\sigma}{2\sigma+\alpha}\right),  \quad B_\sigma = \left(\frac{\sigma}{2(\sigma+\alpha)},\frac{\sigma}{2(\sigma+\alpha)}\right).
\end{equation*}
For $\sigma\geq\alpha$, the point $E_\sigma$  is
\begin{equation*}
 E_\sigma = \left(\frac{\sigma-\alpha}{2\sigma},\frac{1}{2}\right).
\end{equation*}
Thus $[E_\sigma, C]$ is the set of sharp $\sigma$-admissible pairs, $B_\sigma\in(A_\sigma,C)$ lies on the diagonal, and $A_\sigma$ lies on the line through $O$ given by $(2\sigma-\alpha)/q=2\sigma/r$.
In all our applications $\sigma\ge\frac12\ge\frac\alpha2$, so $A_\sigma$ and $B_\sigma$ are always defined. When $2\sigma=\alpha$, the line $OA_\sigma$ is vertical; the cross-multiplied formulas below include this case. The notation $ODE_\sigma A_\sigma$ is used only when $\sigma\ge\alpha$, so that $E_\sigma$ is defined.

For points $X_j\in\R^2$, we use $[X_1,X_2]$, $[X_1,X_2)$ and
$(X_1,X_2)$ for line segments with the indicated endpoint conventions,
and $X_1X_2X_3$ or $X_1X_2X_3X_4$ for their closed convex hulls.
With $x=1/r$ and $y=1/q$,
\[
 \begin{split}
 \operatorname{int}(OA_\sigma C)
 =\{(x,y)\in(0,1/2)^2: \  
 &(2\sigma-\alpha)y<2\sigma x,\\
 &\sigma(1/2-x)>\alpha y\,\}.
 \end{split}
\]
When $\sigma\ge\alpha$, the other region is
\[
 \begin{split}
 \operatorname{int}(ODE_\sigma A_\sigma)
 =\{(x,y)\in(0,1/2)^2: \
 &(2\sigma-\alpha)y>2\sigma x,\\
 &\sigma(1/2-x)>\alpha y\, \}.
 \end{split}
\]

When $\sigma\ge\alpha$, the admissible region $ODE_\sigma C$ splits into the two parts $OA_\sigma C$ and $ODE_\sigma A_\sigma$. For $q,r\ge2$ we define $\beta_\sigma(q,r)$ by
\begin{equation}\label{beta-def}
    \frac{\sigma}{\beta_{\sigma}(q,r)}=\frac{\alpha}{q}+\frac{2\sigma}{r}.
\end{equation}
Note that
\begin{equation*}
	\beta_\sigma(q,r) =
	\left\{
	\begin{array}{lll}
	q/2,\quad &(1/r, 1/q)\in [O, A_\sigma],\\[3pt]
	\frac{2r}{r+2},\quad &(1/r, 1/q) \in [A_\sigma, C],\\[3pt]
	r/2,\quad &(1/r, 1/q)\in [O, C],
	\end{array}\right.
\end{equation*}
and $\beta_\sigma(q,r)<q/2$ in $\mathrm{int}\,OA_\sigma C$, $\beta_\sigma(q,r)>q/2$ in $\mathrm{int}\,ODE_\sigma A_\sigma$.
For the fractional Schr\"odinger equation we show that \eqref{goal-esti}
holds with $\beta=\beta_{d/2}(q,r)$ for
$(1/r,1/q)\in\textrm{int}\,OA_{d/2}C\cup(A_{d/2},C]$, and that this is sharp.
In $ODE_{d/2}A_{d/2}\setminus\{O,D\}$ we obtain \eqref{goal-esti}
for $1\le\beta<q/2$ when $q>2$, and for $\beta=1$ when $q=2$;
\eqref{goal-esti} fails for $\beta>q/2$ if $\mu$ is Ahlfors regular.

\begin{figure}[t]
\centering
\begin{tikzpicture}[scale=.75]
\fill[blue!10] (0,0)--(0,6)--(2,6)--(3,4.5)--cycle;
\draw[->,>=latex] (0,0)--(6.8,0) node[below] {$1/r$};
\draw[->,>=latex] (0,6)--(0,6.8) node[left] {$1/q$};
\draw[thick] (6,0)--(2,6);
\draw[blue,thick] (0,6)--(2,6)--(3,4.5);
\draw[dashed,thick] (0,0)--(3,4.5);
\draw[dashed,thick] (0,0)--(0,6);
\draw[dotted] (0,3)--(2,3);
\node at (0,0) [below left] {$O$};
\node at (0,5.7) [above left] {$D$};
\node at (2,6) [above] {$E_{d/2}$};
\node at (3,4.5) [right] {$A_{d/2}$};
\node at (6,0) [below] {$C$};
\node at (1.35,4.4) {$\ell^{q/2}$};
\node at (-0.4,3.5) [below] {$\frac14$};
\filldraw[draw=black,fill=white] (3,4.5) circle (.08);
\filldraw[draw=black,fill=white] (0,6) circle (.08);
\filldraw[draw=black,fill=white] (0,0) circle (.08);
\end{tikzpicture}
\caption{The Schr\"odinger equation, $d>2\alpha$ (schematic).
The shaded region and the bold finite-$r$ edges satisfy the critical
strong estimate. The dashed boundaries are excluded; $A_{d/2}$ has
Lorentz sequence data. The $r=\infty$ edge has the restricted weak
estimate of part (3).}\label{fig-quadratic-critical}
\end{figure}
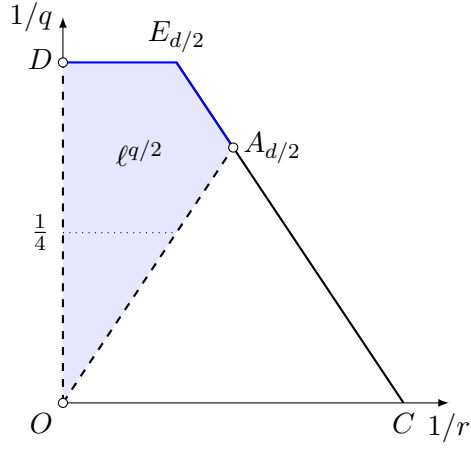

The first result concerns fractional Schr\"odinger evolution, for which $\sigma=d/2$.
\begin{thm}[The fractional Schr\"odinger equation]\label{thm-frac-schro}
Let $a\in \mathbb{R}_{>0}\setminus\{1\}$ and $d\geq1$. Suppose that $\mu$ satisfies \eqref{frac_meas} and $(q,r)$ satisfies \eqref{frac_adm} with $\sigma=d/2$ and $\frac{d}{r}+\frac{a}{q}<d$. Let $s=s_{a\alpha}(q,r)$.
    \begin{enumerate}
    \item If $(\frac{1}{r},\frac{1}{q})$ belongs to $\mathrm{int}\,OA_{d/2}C \cup (A_{d/2}, C]$, then there is a constant $C>0$ such that
    \begin{equation}\label{fra-equa-ineq}
        \Bigg\| \sum_{j\in J} n_j |e^{it(-\Delta)^\frac{a}{2}} f_j|^2 \Bigg\|_{L_t^{\frac{q}{2}}(\mathrm d\mu; L_x^{\frac{r}{2}})} \le C \| \{n_j\}_{j\in J}\|_{\ell^\beta}
    \end{equation}
    holds for all orthonormal families $\{f_j\}_{j\in J}$ in $\dot{H}^s(\mathbb{R}^d)$ with $\beta = \beta_{d/2}(q,r)$.
    By Theorem~\ref{thm-necessary}, the estimate fails for $\beta>\beta_{d/2}(q,r)$.

    \item Suppose $\alpha\le d/2$. If $(\frac{1}{r},\frac{1}{q})$ belongs to $ODE_{d/2}A_{d/2} \setminus \{O,D\}$, then \eqref{fra-equa-ineq} holds for all orthonormal families $\{f_j\}_{j\in J}$ in $\dot{H}^s(\mathbb{R}^d)$ and every $1\le\beta<\frac{q}{2}$ when $q>2$; at $q=2$ it holds with $\beta=1$. If $\mu$ satisfies \eqref{lower-regular} in addition, then, by Theorem~\ref{thm-necessary}, \eqref{fra-equa-ineq} fails for $\beta>\frac q2$.
    \end{enumerate}
\end{thm}

\begin{rem}\label{rem-range-schro}
Along $(A_{d/2},C]$ one has $2\le r<\frac{2(d+\alpha)}{d-\alpha}$
(read as $r<\infty$ if $d=\alpha$) and $\beta=\frac{2r}{r+2}$.
For $\alpha=1$ this is the Frank--Sabin range $2\le r<\frac{2(d+1)}{d-1}$.
The endpoint $r=\frac{2(d+\alpha)}{d-\alpha}$ corresponds to $A_{d/2}$,
where Theorem~\ref{thm-quadratic-critical}~(1) gives a Lorentz estimate
for $a=2$. If $\alpha\le d/2$, the segment $[A_{d/2},E_{d/2}]$ corresponds to
$\frac{2(d+\alpha)}{d-\alpha}\le r\le\frac{2d}{d-2\alpha}$.
\end{rem}

For the Schr\"odinger equation ($a=2$) the regularity vanishes along the
sharp segment, $s_{2\alpha}(q,r)=0$, so the critical endpoint arguments can
use the global energy identity for $e^{it\Delta}$ on $L^2$; for $a\ne2$
this input is not available (see Section~\ref{sec-quadratic-limitations}).
We obtain the following critical refinement. Write
\[
 q_A=\frac{2(d+\alpha)}d,\qquad r_A=\frac{2(d+\alpha)}{d-\alpha},
 \qquad P_A=q_A/2.
\]
\begin{thm}[Critical Schr\"odinger estimates]\label{thm-quadratic-critical}
Let $d\ge2$ and let $\mu$ satisfy \eqref{frac_meas}.
\begin{enumerate}
\item For every orthonormal family $(f_j)$ in $L^2(\mathbb R^d)$,
\begin{equation}\label{quadratic-A}
 \left\|\sum_j n_j|e^{it(-\Delta)}f_j|^2\right\|_
 {L_t^{P_A}(\rd\mu;L_x^{r_A/2})}
 \le C_{d,\alpha}\Fmu^{1/P_A}\|n\|_{\ell^{P_A,1}}.
\end{equation}
\item Let $d>2\alpha$, and let $2\le q<\infty$, $2\le r<\infty$ satisfy
\begin{equation}\label{quadratic-range}
 \frac d2\left(\frac12-\frac1r\right)\ge\frac\alpha q,
 \qquad r>\frac{dq}{d-\alpha}.
\end{equation}
Then, for every orthonormal family $(f_j)$ in $\dot H^s(\mathbb R^d)$
with $s=s_{2\alpha}(q,r)$,
\begin{equation}\label{quadratic-strong}
 \left\|\sum_j n_j|e^{it(-\Delta)}f_j|^2\right\|_
 {L_t^{q/2}(\rd\mu;L_x^{r/2})}
 \le C_{d,\alpha,q,r}\Fmu^{2/q}\|n\|_{\ell^{q/2}}.
\end{equation}
\item Let $d>2\alpha$ and $2<q<\infty$. Then, for every orthonormal family
$(f_j)$ in $\dot H^s(\mathbb R^d)$ with $s=s_{2\alpha}(q,\infty)$,
\begin{equation}\label{quadratic-ODweak}
 \left\|\sum_j n_j|e^{it(-\Delta)}f_j|^2\right\|_
 {L_t^{q/2,\infty}(\rd\mu;L_x^\infty)}
 \le C_{d,\alpha,q}\Fmu^{2/q}\|n\|_{\ell^{q/2,1}}.
\end{equation}
\end{enumerate}
\end{thm}
Condition \eqref{quadratic-range} describes the points of
$ODE_{d/2}A_{d/2}$ with $r<\infty$ lying strictly to the left of the line
$OA_{d/2}$; this set includes $(A_{d/2},E_{d/2}]$ and
$[D,E_{d/2}]\setminus\{D\}$ (see Figure~\ref{fig-quadratic-critical}).

\begin{rem}\label{rem-critical-scope}
For Lebesgue measure on an interval $I$ ($\alpha=1$), one has
$\langle\mu\rangle_1\le2$ independently of $I$, so letting $I\uparrow\mathbb R$
gives global-in-time estimates. Thus part (1) recovers the endpoint
estimate of Hadama--Toyoshima \cite{HT2026}, and, for $d\ge3$, part (2)
gives the critical estimate in the whole range \eqref{quadratic-range}
with $\alpha=1$. This range contains the sharp segment treated in
\cite{HT2026} and the case $q=4$ in \cite{Gong2026}; for the remaining
interior, announced in \cite{FSW2025}, our proof does not use the weak
mixed-norm interpolation identity employed there. For Ahlfors regular
$\mu$, Theorem~\ref{thm-necessary} shows that the exponent $\beta=q/2$ in
part (2) cannot be increased. The estimate with $\beta=q/2$ on $OA_{d/2}$ (at $A_{d/2}$
only $\ell^{P_A,1}$ is obtained) and the strong estimate on $r=\infty$ are
not established here; see Section~\ref{sec-quadratic-limitations}.
\end{rem}

For the wave equation, $\sigma=\frac{d-1}2$, and we obtain the following analogue.
\begin{thm}[The wave equation]\label{thm-wave} Let $d\geq2$. Suppose that $\mu$ satisfies \eqref{frac_meas} and $(q,r)$ satisfies \eqref{frac_adm} with $\sigma=\frac{d-1}2$. Let $s=s_{\alpha}(q,r)$.
\begin{enumerate}
    \item If $(\frac{1}{r},\frac{1}{q})$ belongs to $\mathrm{int}\,OA_{\frac{d-1}{2}}C \cup (A_{\frac{d-1}2}, C]$, then there is a constant $C>0$ such that
    \begin{equation}\label{wave-equa-ineq}
        \Bigg\| \sum_{j\in J} n_j |e^{it\sqrt{-\Delta}} f_j|^2 \Bigg\|_{L_t^{\frac{q}{2}}(\mathrm d\mu; L_x^{\frac{r}{2}})} \le C \| \{n_j\}_{j\in J}\|_{\ell^\beta}
    \end{equation}
    holds for all orthonormal families $\{f_j\}_{j\in J}$ in $\dot{H}^s(\mathbb{R}^d)$ with $\beta = \beta_{\frac{d-1}2}(q,r)$.
    By Theorem~\ref{thm-necessary}, the estimate fails for $\beta>\beta_{\frac{d}2}(q,r)$.

    \item Suppose $\alpha\le \frac{d-1}2$. If $(\frac{1}{r},\frac{1}{q})$ belongs to $ODE_{\frac{d-1}{2}}A_{\frac{d-1}{2}} \setminus \{O,D\}$, then \eqref{wave-equa-ineq} holds for all orthonormal families $\{f_j\}_{j\in J}$ in $\dot{H}^s(\mathbb{R}^d)$ and every $1\le\beta<\frac{q}{2}$ when $q>2$; at $q=2$ it holds with $\beta=1$. If $\mu$ satisfies \eqref{lower-regular} in addition, then, by Theorem~\ref{thm-necessary}, \eqref{wave-equa-ineq} fails for $\beta>\frac q2$.
  \end{enumerate}
\end{thm}

In the case of the Klein--Gordon equation, we obtain the following result with initial data in inhomogeneous Sobolev spaces.
\begin{thm}[The Klein--Gordon equation]\label{thm-KG}
Let $d\geq2$. Suppose that $\mu$ satisfies \eqref{frac_meas} and $(q,r)$ satisfies \eqref{frac_adm} with $\sigma=\frac{d-1}2$. Let $s\ge s_{\alpha}(q,r)$.
\begin{enumerate}
\item If $(\frac{1}{r},\frac{1}{q})$ belongs to $\mathrm{int}\,OA_{\frac{d-1}{2}}C \cup (A_{\frac{d-1}2}, C]$, then there is a constant $C>0$ such that
     \begin{equation}\label{K-G-equa-ineq}
        \Bigg\| \sum_{j\in J} n_j |e^{it\sqrt{1-\Delta}} f_j|^2 \Bigg\|_{L_t^{\frac{q}{2}}(\mathrm d\mu; L_x^{\frac{r}{2}})} \le C \| \{n_j\}_{j\in J}\|_{\ell^\beta}
    \end{equation}
holds for all orthonormal families $\{f_j\}_{j\in J}$ in ${H}^s(\mathbb{R}^d)$ with $\beta=\beta_{\frac{d-1}2}(q,r)$. By Theorem~\ref{thm-necessary}, the estimate fails for $\beta>\beta_{\frac d2}(q,r)$ when $s=s_\alpha(q,r)$.

\item Suppose $\alpha\le\frac{d-1}2$. If $(\frac{1}{r},\frac{1}{q})$ belongs to $ODE_{\frac{d-1}{2}}A_{\frac{d-1}{2}} \setminus \{O,D\}$, then \eqref{K-G-equa-ineq} holds for all orthonormal families $\{f_j\}_{j\in J}$ in ${H}^s(\mathbb{R}^d)$ and every $1\le\beta<\frac{q}{2}$ when $q>2$; at $q=2$ it holds with $\beta=1$. If $\mu$ satisfies \eqref{lower-regular} in addition, then, by Theorem~\ref{thm-necessary}, \eqref{K-G-equa-ineq} fails for $\beta>\frac q2$ when $s=s_\alpha(q,r)$.
\end{enumerate}
\end{thm}

\begin{rem}\label{rem-range-wave}
For the wave and Klein--Gordon equations, along $(A_{\frac{d-1}2},C]$ one has
$2\le r<\frac{2(d-1+\alpha)}{d-1-\alpha}$ (read as $r<\infty$ if $d-1=\alpha$)
and $\beta=\frac{2r}{r+2}$. For $\alpha=1$ and $d\ge3$ this is the
Lebesgue-time range $2\le r<\frac{2d}{d-2}$ recalled above. Unlike the
Schr\"odinger case, the regularity
$s_\alpha(q,r)=\frac{d+1}2\big(\frac12-\frac1r\big)$ is positive on the sharp
line, so the energy input used for Theorem~\ref{thm-quadratic-critical}~(1)
is not available, and we do not obtain an endpoint estimate at
$A_{\frac{d-1}2}$; see Section~\ref{sec-quadratic-limitations}. Moreover,
the exponent $\beta_{\frac{d-1}2}(q,r)$ need not be optimal: by
Theorem~\ref{thm-necessary} only $\beta\le\beta_{\frac d2}(q,r)$ is known to be
necessary (see the discussion after Theorem~\ref{thm-necessary}). If
$\alpha\le\frac{d-1}2$, the segment $[A_{\frac{d-1}2},E_{\frac{d-1}2}]$
corresponds to $\frac{2(d-1+\alpha)}{d-1-\alpha}\le r\le\frac{2(d-1)}{d-1-2\alpha}$,
where $E_{\frac{d-1}2}=D$ is excluded when $d-1=2\alpha$.
\end{rem}

\begin{rem}\label{rem-small-sigma}
The assumption $\alpha\le\sigma$ in part (2) of
Theorems~\ref{thm-frac-schro}, \ref{thm-wave} and \ref{thm-KG} fails only
when $\sigma=\frac12<\alpha$, that is, for the fractional Schr\"odinger
equation with $d=1$ and for the wave and Klein--Gordon equations with
$d=2$. Then $D$ is not admissible, the admissible region is the triangle
$OCF$ with $F=(0,\frac1{4\alpha})$, and part (1) remains valid. Let $H$ be
the convex hull of $\mathrm{int}\,OA_\sigma C$ and the segment
$\{(0,y):0<y<\frac14\}$. The interpolation argument of
Section~\ref{sec-nonsharp}, with Proposition~\ref{prop-boundary} on
$r=\infty$, gives the estimates for every $\beta<q/2$ at all points of
$(H\cap OFA_\sigma)\setminus\{O\}$; for $\alpha=1$, where $F=(0,\frac14)$,
this set is $OFA_\sigma\setminus(\{O\}\cup[F,A_\sigma])$. For
$\frac12<\alpha<1$, Proposition~\ref{prop-lowdim-band} gives further
estimates on $r=\infty$ for $4\alpha<q\le4$, and at related finite $r$,
for every $\beta$ below explicit exponents such as $B_\alpha(q)$; these
are smaller than $q/2$ except at $q=4$. The optimal $\beta$ on the rest of
$OFA_\sigma$ is not determined here.
\end{rem}

We now state the necessary conditions on $\beta$ in Theorems~\ref{thm-frac-schro}, \ref{thm-wave}, \ref{thm-KG}.
\begin{thm}[Necessary conditions]\label{thm-necessary}
Let $\mu$ satisfy \eqref{frac_meas} and let $2\le q,r\le\infty$. Suppose that \eqref{fra-equa-ineq} holds for all orthonormal families in $\dot H^s$ with $s=s_{a\alpha}(q,r)$, that \eqref{wave-equa-ineq} holds for all orthonormal families in $\dot H^s$, or that \eqref{K-G-equa-ineq} holds for all orthonormal families in $H^s$, with $s=s_\alpha(q,r)$ in the last two cases. Then
\begin{enumerate}
\item $\beta\le \beta_{d/2}(q,r)$;
\item if $\mu$ satisfies \eqref{lower-regular} in addition, then $\beta\le q/2$.
\end{enumerate}
\end{thm}
By \eqref{beta-def}, $\beta_{d/2}(q,r)<q/2$ in $\mathrm{int}\,OA_{d/2}C$,
$\beta_{d/2}(q,r)=q/2$ on $OA_{d/2}$, and $\beta_{d/2}(q,r)>q/2$ to the left
of $OA_{d/2}$. Hence, for the fractional Schr\"odinger equation,
Theorem~\ref{thm-necessary}~(1) shows that Theorem~\ref{thm-frac-schro}~(1)
is sharp for every $\mu$ satisfying \eqref{frac_meas}, while, under
\eqref{lower-regular}, part (2) shows that the range $\beta<q/2$ in
Theorem~\ref{thm-frac-schro}~(2) and the exponent $\beta=q/2$ in
Theorem~\ref{thm-quadratic-critical}~(2) cannot be improved, apart from the
endpoint $\beta=q/2$ in the former. Since the line $OA_\sigma$ becomes steeper as
$\sigma$ decreases, $ODE_{\frac{d-1}2}A_{\frac{d-1}2}$ lies to the left of
$OA_{d/2}$, so the same conclusion holds for the wave and Klein--Gordon
equations there. In $\mathrm{int}\,OA_{\frac{d-1}2}C\cup(A_{\frac{d-1}2},C]$,
however, there is a gap between the exponent $\beta_{\frac{d-1}2}(q,r)$
obtained here and the necessary condition $\beta\le\beta_{d/2}(q,r)$. The
lower-bound question in \cite{BLN2021} concerns strengthening the necessary
condition, whereas \cite[Conjecture~1.2]{BKSnote2025} proposes improving the
sufficient exponent toward $\beta_{d/2}$. We do not address this problem for
Frostman time.

\subsection*{Organization}
In Section~\ref{sec-prelim}, we recall
the function spaces used below and present preliminary results: the
fractional integral inequality for $\mu$, the duality principle,
dispersive estimates, and the scalar estimates over fractal time,
including the Klein--Gordon case and the annular synthesis for finite $r$.
Sections~\ref{sec-sharp}, \ref{sec-nonsharp} and \ref{sec-OD}
prove the sufficiency parts of Theorems~\ref{thm-frac-schro},
\ref{thm-wave} and \ref{thm-KG} for sharp admissible pairs, for non-sharp
pairs, and on the boundary $r=\infty$, respectively;
Section~\ref{sec-OD} also contains the low-dimensional estimates of
Remark~\ref{rem-small-sigma}. Section~\ref{sec-critical} proves
Theorem~\ref{thm-quadratic-critical} and discusses its endpoint
limitations. Section~\ref{sec-necessary} proves
Theorem~\ref{thm-necessary}, and Section~\ref{sec-appli} reformulates the
results in terms of density matrices.

\section{Preliminaries}\label{sec-prelim}
\subsection{Function spaces} We recall the function spaces and interpolation identities used below. Unless the time measure is displayed explicitly as $\rd t$, all mixed time norms are taken with respect to $\rd\mu$.

The spaces $\dot{H}^s(\mathbb{R}^d)$ and $H^s(\mathbb{R}^d)$ are the homogeneous and inhomogeneous Sobolev spaces with the norms
\begin{equation*}
\|f\|_{\dot{H}^s(\mathbb{R}^d)}=\||D|^sf\|_{L^2(\mathbb{R}^d)},\quad \|f\|_{{H}^s(\mathbb{R}^d)}=\|\langle D\rangle^sf\|_{L^2(\mathbb{R}^d)}.
\end{equation*}
Throughout this work, $\Psi(\xi)$ denotes either $|\xi|$ or $\langle\xi\rangle=(1+|\xi|^2)^{1/2}$.

For $1<p<\infty$ and $1\leq r\leq \infty$, the Lorentz space $L^{p,r}(\mathrm d\mu)$ is the set of $\mu$-measurable functions $f$ on $\mathbb{R}$ such that $\|f\|_{L^{p,r}(\mathrm d\mu)}<\infty$, where
\begin{equation*}
\|f\|_{L^{p,r}(\mathrm d\mu)}=
\begin{cases}
\left(\int_0^\infty (u^\frac{1}{p}f^*(u))^r\frac{\mathrm du}{u}\right)^\frac{1}{r},&\text{ if } 1\leq r<\infty,\\
\sup\limits_{u>0}u^\frac{1}{p}f^*(u),&\text{ if } r=\infty.
\end{cases}
\end{equation*}
Here the decreasing rearrangement $f^*$ of $f$ is defined by
$f^*(u)=\inf\{\lambda>0: a_f(\lambda)\leq u\}$,
 where $a_f(\lambda)=\mu\left(\{\tau\in\mathbb{R}: |f(\tau)|>\lambda\}\right)$ is the distribution function of $f$. The functional $\|\cdot\|_{L^{p,r}(\mathrm d\mu)}$ is a norm if $r\leq p$ and a quasi-norm otherwise; replacing $f^*(u)$ with $\frac{1}{u}\int_0^uf^*$ yields an equivalent norm for all $1<p<\infty$, $1\le r\le\infty$.
The Lorentz sequence spaces $\ell^{p,r}$ are defined in the same manner using the decreasing rearrangement $\{\lambda^*_j\}^\infty_{j=1}$ of the absolute values of $\{\lambda_j\}_{j=1}^\infty\in c_0$:
\begin{equation*}
\|\{\lambda_j\}^\infty_{j=1}\|_{\ell^{p,r}}=
\begin{cases}
\left(\sum\limits_{j=1}^\infty (j^\frac{1}{p}\lambda_j^*)^r\frac{1}{j}\right)^\frac{1}{r},&\text{ if } 1\leq r<\infty,\\
\sup\limits_{j\geq1}j^\frac{1}{p}\lambda_j^*,&\text{ if } r=\infty.
\end{cases}
\end{equation*}
We refer the reader to \cite{SW1971} for further details.
We have $L^{p,p}(\mathrm d\mu)=L^p(\mathrm d\mu)$, and the Lorentz spaces increase with the second exponent:
\begin{equation}\label{embedding}
    L^{p,r_1}(\mathrm d\mu)\subseteq L^{p,r_2}(\mathrm d\mu),\quad 1\leq r_1\leq r_2\leq \infty.
\end{equation}
Since $\mu$ is finite, we also have $L^{p_1,r_1}(\mathrm d\mu)\subseteq L^{p_2,r_2}(\mathrm d\mu)$ whenever $p_1>p_2$.

We use the matching-index Lions--Peetre identity
\begin{equation}\label{real-L}
 (L_t^{q_0}(\rd\mu;B_0),L_t^{q_1}(\rd\mu;B_1))_{\theta,q}
 =L_t^q(\rd\mu;(B_0,B_1)_{\theta,q}),
\end{equation}
where $1\le q_0,q_1<\infty$, $0<\theta<1$,
$1/q=(1-\theta)/q_0+\theta/q_1$, and $(B_0,B_1)$ is a compatible
Banach couple. The equality is up to equivalent norms; see
\cite[Chapter VII, Section 1]{LP1964}.
For $B_i=L_x^{r_i}$ with $r_0\ne r_1$, the inner space is
$L_x^{r,q}$, where $1/r=(1-\theta)/r_0+\theta/r_1$.
If $r_0=r_1$, the inner space is $L_x^{r_0}$ itself.
The identity requires the real interpolation second index to coincide
with the outer time exponent $q$; for weak outer spaces the corresponding
statement fails (see Section~\ref{sec-quadratic-limitations}). We also use
\begin{equation}\label{real-ell}
\big(\ell^{p_0,r_0},\ell^{p_1,r_1}\big)_{\theta,r}=\ell^{p,r}
\end{equation}
whenever $p_0, p_1\in[1,\infty)$, $p_0\neq p_1$, $\frac{1}{p}=\frac{1-\theta}{p_0}+\frac{\theta}{p_1}$, $\theta\in(0,1)$ and $r_0,r_1,r\in[1,\infty]$. Similarly, $(L^{p_0}_t(\mathrm d\mu;L^2_x),L^{p_1}_t(\mathrm d\mu;L^2_x))_{\theta,r}=L_t^{p,r}(\mathrm d\mu;L^2_x)$. See \cite{BL1976, C1974}.

\subsection{Hardy--Littlewood--Sobolev inequality for fractal measures}
\begin{lem}\label{lem_HLS}
Let $\alpha\in(0,1]$, $\lambda\in(0,\alpha)$, and suppose that $\mu$ satisfies \eqref{frac_meas}.
Let $1<p,p_1<\infty$ satisfy $2 = \frac \lambda \alpha + \frac1p+\frac1{p_1}$, and let $r_1, r_2\in[1,\infty]$ satisfy $1/r_1 + 1/r_2 \ge1$. Then
\begin{align}\label{frac_HLS}
	\iint\frac{|g_1(t) g_2(\tau)|}{|t-\tau|^{\lambda}}~\mathrm d\mu(t) \mathrm d\mu(\tau)  \le C\,\Fmu^{\lambda/\alpha} \| g_1 \|_{L^{p,r_1}(\mathrm d\mu)} \|g_2\|_{L^{p_1, r_2}(\mathrm d\mu)}.
\end{align}
\end{lem}
This follows from \cite[Corollary~3.3]{GarEdu2004}, which gives the boundedness of the fractional integral $g\mapsto\int g(\tau)|t-\tau|^{-\lambda}\mathrm d\mu(\tau)$ from $L^{p_1}(\mathrm d\mu)$ to $L^{p'}(\mathrm d\mu)$ for measures satisfying \eqref{frac_meas}, together with real interpolation and H\"older's inequality in Lorentz spaces \cite{O'N1963}. In particular, taking $p=p_1=u_\lambda := \frac{2\alpha}{2\alpha - \lambda}$, we have
\begin{equation}\label{eq: hls}
 \iint\frac{|f(t)g(\tau)|}{|t-\tau|^\lambda}
       \,\mathrm d\mu(t)\,\mathrm d\mu(\tau)
 \leq C_{\alpha,\lambda} \Fmu^{\lambda/\alpha}
       \|f\|_{L^{u_\lambda}(\mathrm d\mu)}
       \|g\|_{L^{u_\lambda}(\mathrm d\mu)}.
\end{equation}

\subsection{Schatten spaces and the duality principle}
For $r\ge 2$, we define $\widetilde r$ by
\begin{equation*}
    \frac{1}{\widetilde{r}}=\frac{1}{2}-\frac{1}{r},
\end{equation*}
and similarly $\widetilde q$, $\widetilde p$. Note that $\widetilde q=2(q/2)'$.
For $\beta\in[1,\infty)$, the Schatten space $\mathfrak S^\beta=\mathfrak S^\beta(\mathfrak H)$ on a separable Hilbert space $\mathfrak H$ is the space of all compact operators $T$ on $\mathfrak H$ whose singular values $\{\lambda_j\}_j$ (the square roots of the eigenvalues of $T^*T$) belong to $\ell^\beta$, with the norm
$\|T\|_{\mathfrak S^\beta}=\|\{\lambda_j\}\|_{\ell^\beta}$.
For $\beta=\infty$, $\mathfrak S^\infty$ denotes the space of bounded operators with the operator norm. The space $\mathfrak S^2$ is the Hilbert--Schmidt class, and $\|T\|_{\mathfrak S^2}=\|K_T\|_{L^2}$ if $T$ is an integral operator with kernel $K_T$. We refer to \cite{BS1979} for more details. We shall mostly use $\mathfrak H=L^2(\mathrm d\mu\, \mathrm dx)$ or $\mathfrak H=L^2(\mathbb R^d)$.

For a fixed flow and Sobolev order, we write
\begin{equation}\label{eq-initial-potential}
 B_{\phi,s}(V)=\Psi(D)^{-s}
 \left(\int U_\phi(t)^*V(t)U_\phi(t)\,\rd\mu(t)\right)
 \Psi(D)^{-s}.
\end{equation}
We abbreviate this to $B_s(V)$ when the flow is fixed. Initially the
integral is a sesquilinear form for bounded supported potentials and
finite-frequency test vectors, with homogeneous supports bounded away
from zero. Its bounded or Schatten-class extension is established in
each estimate below.

\begin{prop}[Duality principle]\label{prop_dual}
	Let $2< q<\infty$, $2\le p\le\infty$, $2\le r\le\infty$ and $\beta\ge1$. Suppose $T$ is a bounded operator from $L^2(\mathbb R^d)$ to $L_t^{q,p}(\mathrm d\mu; L_x^r)$. Then the following statements are equivalent:
	\begin{enumerate}
		\item For every orthonormal family $\{f_j\}_{j\in J}$ in $L^2(\mathbb R^d)$ and $\{n_j\}_{j} \in \ell^\beta$,
		\[
			\Bigg\| \sum_{j\in J}  n_j |T f_j|^2 \Bigg\|_{L_t^{\frac{q}{2}, \frac{p}{2}}(\mathrm d\mu; L_x^{\frac{r}{2}})} \le C \|\{n_j\}_{j\in J}\|_{\ell^\beta}.
		\]

		\item For every $W \in L_t^{\widetilde q, \widetilde p}L_x^{\widetilde r}$,
		\[
			\| W T T^* \overline{W} \|_{\mathfrak S^{\beta'}} \le C' \|W\|_{L_t^{\widetilde q, \widetilde p}L_x^{\widetilde r}}^2.
		\]
	\end{enumerate}
	Moreover, (2) is equivalent to the bilinear bound
\[
 \|W_1TT^*W_2\|_{\mathfrak S^{\beta'}}
 \le C''\prod_{i=1}^2\|W_i\|_{L_t^{\widetilde q,\widetilde p}(L_x^{\widetilde r})}.
\]
\end{prop}
Here $|W|^2\in L_t^{(q/2)',(p/2)'}(\mathrm d\mu;L_x^{(r/2)'})$ belongs to the norming associate space of the left-hand side of (1), which is a Banach space for $q>2$ after renorming.
The proof follows the argument of \cite[Lemma~3]{FraSa2017}. We include it for the reader's convenience.

\begin{proof}
The nonzero singular values of the two Gram operators agree, so
\[
 \|WTT^*\overline W\|_{\mathfrak S^{\beta'}}
 =\|(WT)^*(WT)\|_{\mathfrak S^{\beta'}}
 =\|T^*|W|^2T\|_{\mathfrak S^{\beta'}}.
\]
Thus statement (2) is equivalent to
\[
 \|T^*|W|^2T\|_{\mathfrak S^{\beta'}}
 \le C'\||W|^2\|_{L_t^{(q/2)',(p/2)'}(L_x^{(r/2)'})}.
\]
Let $\{n_j\}$ be nonnegative (which suffices), let $\{f_j\}$ be orthonormal, and set $\gamma=\sum_j n_j|f_j\rangle\langle f_j|$, so that $\|\gamma\|_{\mathfrak S^\beta}=\|\{n_j\}\|_{\ell^\beta}$ and $\rho_{T\gamma T^*}=\sum_jn_j|Tf_j|^2$. For $V\ge0$,
\[
\int_{\mathbb R^d\times E}\Big(\sum_jn_j|Tf_j|^2\Big)V\,\mathrm dx\,\mathrm d\mu
=\mathrm{Tr}_{L^2(\mathbb R^d)}\big(\gamma\, T^*VT\big)\le\|\gamma\|_{\mathfrak S^\beta}\|T^*VT\|_{\mathfrak S^{\beta'}}
\]
by H\"older's inequality in Schatten spaces. Taking $V=|W|^2$ and the supremum over $V$ in the unit ball of the norming associate space, (2) implies (1). Conversely, if (1) holds, then for every $\gamma\ge0$ of finite rank, $\mathrm{Tr}(\gamma T^*VT)\le C\|\gamma\|_{\mathfrak S^\beta}\|V\|_{L_t^{(q/2)',(p/2)'}(L_x^{(r/2)'})}$, and by duality $\|T^*VT\|_{\mathfrak S^{\beta'}}\le C\|V\|$ for $V\ge0$, which gives (2). The bilinear version follows from (2) by writing $W_1TT^*W_2=(W_1T)(\overline{W_2}T)^*$ and using H\"older's inequality in Schatten spaces together with $\|AB^*\|_{\mathfrak S^{\beta'}}\le\|AA^*\|_{\mathfrak S^{\beta'}}^{1/2}\|BB^*\|_{\mathfrak S^{\beta'}}^{1/2}$.
\end{proof}

\subsection{Dispersive estimates}\label{sec-disp}
Let $\chi_0\in C_c^\infty(\mathbb R)$ be even, nonnegative, at most one, nonincreasing on positive radii, and supported in $(-2,2)$ with $\chi_0\equiv 1$ on $(-1, 1)$, and let $\psi = \chi_0 - \chi_0(2\cdot)$ and $\psi_k = \psi(2^{-k} \cdot)$, $k\in\mathbb Z$, so that $\supp \psi_k \subset [2^{k-1} , 2^{k+1}]$ and $\sum_{k\in\mathbb Z}\psi_k(\rho)=1$ for $\rho>0$, $\chi_0(\rho)+\sum_{k\ge1}\psi_k(\rho)=1$.
We define $P_k$ and $P_{\le 0}$ by $\widehat{P_k g}(\xi) = \psi_k(|\xi|) \widehat{g}(\xi)$, $\widehat{P_{\le0} g}(\xi) = \chi_0(|\xi|) \widehat{g}(\xi)$.

In what follows, we consider the three cases
\begin{equation}\label{three-cases}
 (\phi,\Psi,a,\sigma)=
 \begin{cases}
 (|\xi|^a,|\xi|,a,d/2),&a\ne1,\\
 (|\xi|,|\xi|,1,(d-1)/2),&d\ge2,\\
 (\langle\xi\rangle,\langle\xi\rangle,1,(d-1)/2),&d\ge2,
 \end{cases}
\end{equation}
and for $s\in\mathbb R$, $\kappa\in\mathbb R$ and $k\in\mathbb Z$ we set
\[
\begin{aligned}
 K_k^{s,\kappa}(t,x)
 &=\int_{\mathbb R^d}\e^{i(x\cdot\xi+t\phi(\xi))}
 \Psi(\xi)^{-2s+i\kappa}\psi_k(|\xi|)\,\rd\xi,\\
 K_{\le0}^{s,\kappa}(t,x)
 &=\int_{\mathbb R^d}\e^{i(x\cdot\xi+t\phi(\xi))}
 \Psi(\xi)^{-2s+i\kappa}\chi_0(|\xi|)\,\rd\xi.
\end{aligned}
\]
\begin{lem}\label{lem-dyadic-disp}
Let $(\phi,\Psi,a,\sigma)$ be one of \eqref{three-cases}. There exist $C,M>0$, depending on $s$, such that
\begin{equation}\label{dyadic-disp}
|K_k^{s,\kappa}(t,x)|\le C(1+|\kappa|)^M\,2^{k(d-2s)}\big(1+2^{ak}|t|\big)^{-\sigma}
\end{equation}
for all $(t,x)\in\mathbb R\times\mathbb R^d$ and $\kappa\in\mathbb R$, where $k\in\mathbb Z$ in the first two cases and $k\ge1$ for the Klein--Gordon equation. In the Klein--Gordon case, we also have $|K_{\le0}^{s,\kappa}(t,x)|\le C(1+|\kappa|)^M$.
\end{lem}
\begin{proof}
In the homogeneous cases, the change of variables $\xi\to2^k\xi$ gives $K_k^{s,\kappa}(t,x)=2^{k(d-2s+i\kappa)}K_0^{s,\kappa}(2^{ak}t,2^kx)$, so it suffices to consider $k=0$. Then \eqref{dyadic-disp} follows from the standard stationary phase estimate (see, e.g., \cite{GNY2018, BLN2021}): on $\supp\psi_0$, the Hessian of $|\xi|^a$ is nondegenerate if $a\ne1$, and the Hessian of $|\xi|$ has $d-1$ nonvanishing eigenvalues. The derivatives of the amplitude $|\xi|^{-2s+i\kappa}\psi_0(|\xi|)$ grow at most polynomially in $\kappa$. For the Klein--Gordon equation and $k\ge1$, the change of variables $\xi\to2^k\xi$ gives the phase $2^kt\,(2^{-2k}+|\xi|^2)^{1/2}$, whose Hessian on $\supp\psi_0$ has $d-1$ eigenvalues comparable to $1$, uniformly in $k\ge1$. The bound for $K^{s,\kappa}_{\le0}$ is trivial.
\end{proof}

Summing \eqref{dyadic-disp} over $k$, we obtain the following.
\begin{lem}\label{lem-full-decay}
Let $(\phi,\Psi,a,\sigma)$ be one of \eqref{three-cases} and $0<d-2s<a\sigma$. Then, for every $R>0$,
\[
\Big|\sum_{k\in\mathcal K}K_k^{s,\kappa}(t,x)\Big|\le C(1+|\kappa|)^M|t|^{-\frac{d-2s}{a}}, \qquad 0<|t|\le R,
\]
uniformly over finite sets $\mathcal K\subset\mathbb Z$ (in the Klein--Gordon case, $\mathcal K\subset\mathbb N$ and the sum may include $K_{\le0}^{s,\kappa}$). In the homogeneous cases the bound holds for all $t\ne0$.
\end{lem}
\begin{proof}
With $\gamma=d-2s$, the two geometric sums give
\[
 \begin{aligned}
 \sum_k2^{k\gamma}(1+2^{ak}|t|)^{-\sigma}
 &\lesssim\sum_{2^{ak}|t|\le1}2^{k\gamma}
 +|t|^{-\sigma}\sum_{2^{ak}|t|>1}2^{k(\gamma-a\sigma)}\\
 &\lesssim |t|^{-\gamma/a},
 \end{aligned}
\]
since $0<\gamma<a\sigma$. In the Klein--Gordon case the low frequency part is bounded and $1\lesssim_R|t|^{-\gamma}$ for $|t|\le R$.
\end{proof}

We also need oscillatory integral estimates in which the real power of the weight exactly compensates the dispersion: for the weights $|\xi|^{(a-2)d/2}$ and $\Psi(\xi)^{-(d+1)/2}$ the kernels decay like $|t|^{-d/2}$ and $|t|^{-(d-1)/2}$, respectively, uniformly in the frequency cutoff.
Let $\chi\in C_c^\infty(\mathbb R^d)$ be a real-valued, nonnegative, radial function, nonincreasing in the radial variable, with $\chi\equiv1$ near the origin, and set
\[
m_\varepsilon(\xi)
:=
\chi(\varepsilon\xi)-\chi(\xi/\varepsilon),\qquad 0<\varepsilon<1.
\]
Then \(0\le m_\varepsilon\le1\), \(m_\varepsilon\) is smooth and compactly supported away from the origin, and $m_\varepsilon(\xi)\to1$ for every $\xi\neq0$ as $\varepsilon\to0$.

\begin{prop}[Proposition~6, \cite{BLN2021}]\label{prop_disp}
Let $\phi(\xi)=|\xi|^a$, $a\in\mathbb R_{>0}\setminus\{1\}$. There exist constants $C<\infty$ and $M>0$ such that
	\[
		\Big|\int_{\mathbb R^d} \e^{i(x\cdot\xi + t\phi(\xi))}  |\xi|^{\frac{(a-2)d}{2} + i\kappa}m_\varepsilon(\xi)\,\mathrm d\xi\Big| \le C (1+|\kappa|)^M |t|^{-\frac{d}{2}}
	\]
	for all $\varepsilon\in(0,1)$, $\kappa\in\mathbb R$ and $(x,t) \in \mathbb R^d \times \mathbb R$.
\end{prop}

\begin{prop}[Propositions~3 and 4, \cite{BLN2021}]\label{wave-decay}
Let $d\ge2$ and $(\phi,\Psi)=(|\xi|,|\xi|)$ or $(\langle\xi\rangle,\langle\xi\rangle)$. There exist constants $C<\infty$ and $M>0$ such that
    \begin{equation*}
        |\kappa|\,\Big|\int_{\mathbb R^d}\e^{i(x\cdot\xi+t\phi(\xi))}\Psi(\xi)^{-\frac{d+1}{2}+i\kappa}m_\varepsilon(\xi)\,\mathrm d\xi\Big|\leq C(1+|\kappa|)^M|t|^{-\frac{d-1}{2}}
    \end{equation*}
    for all $\varepsilon\in(0,1)$, $\kappa\in\mathbb{R}$ and $(x,t) \in \mathbb R^d \times \mathbb R$.
\end{prop}
In \cite{BLN2021} these estimates are stated with the high-frequency
cutoff $\chi(\varepsilon\xi)$ in place of $m_\varepsilon$. They persist
for $m_\varepsilon$. Indeed, choose $0<\varepsilon'\le\varepsilon$ so
small that $\chi(\varepsilon'\cdot)\equiv1$ on
$\supp\chi(\cdot/\varepsilon)$. Then
$\chi(\xi/\varepsilon)=\chi(\xi/\varepsilon)\chi(\varepsilon'\xi)$, so the
oscillatory integral with $\chi(\xi/\varepsilon)$ is the spatial
convolution of the one with $\chi(\varepsilon'\xi)$ and
$(2\pi)^{-d/2}\varepsilon^d\check\chi(\varepsilon\cdot)$, whose $L^1$ norm
does not depend on $\varepsilon$. Subtracting it from the integral with
$\chi(\varepsilon\xi)$ gives the stated bounds with $m_\varepsilon$.
The factor $|\kappa|$ in Proposition~\ref{wave-decay} cannot be removed: for $\kappa=0$ the bound fails logarithmically near the light cone; this is why the factor $z-1$ appears in the analytic family in Corollary~\ref{cor-wave-kg}.

\subsection{Scalar Klein--Gordon estimates}
\begin{lem}\label{K-L-Strichartz}
Let $d\ge2$, and let $2\le q\le\infty$, $2\le r<\infty$. Suppose that $(q,r)$ is $\frac{d-1}{2}$-admissible and $\mu$ satisfies \eqref{frac_meas}. Then
    \begin{equation}\label{eq:K-G-equa-fos}
       \| \e^{it\sqrt{1-\Delta}} f \|_{L_t^q(\mathrm d\mu; L_x^r)} \lesssim \| f\|_{H^s(\mathbb R^d)},\qquad s=s_\alpha(q,r).
\end{equation}
The implicit constant depends on $d,q,r$, $\Fmu$ and the diameter of $\supp\mu$.
\end{lem}
\begin{proof}
Let $R_0>0$ be such that $E\subset(-R_0,R_0)$. We write $f=P_{\le0}f+(1-P_{\le0})f$. Since $r\ge2$, Bernstein's inequality and the unitarity of $\e^{it\langle D\rangle}$ give $\|\e^{it\langle D\rangle}P_{\le0}f\|_{L^r_x}\lesssim\|f\|_{L^2}$ uniformly in $t$, and the bound for the low frequency part follows since $\mu$ is finite.
For the high frequency part, we write $\e^{it\langle D\rangle}(1-P_{\le0})=M_t\,\e^{it|D|}(1-P_{\le 0})$, where $M_t$ is the Fourier multiplier with symbol
$m_t(\xi)=\e^{it(\langle\xi\rangle-|\xi|)}(1-\chi_0(2|\xi|))$, a smooth function. The symbol equals the required phase on the support of
$1-\chi_0(|\xi|)$. Since $1-\chi_0(2|\xi|)=1$ on the support of
$1-\chi_0(|\xi|)$, this gives the stated factorization.
Moreover, $b_t=m_t-1$ is a smooth symbol of order $-1$, uniformly
for $|t|\le R_0$. Its high-frequency dyadic inverse Fourier pieces
have $L^1$ norms $O_{R_0}(2^{-k})$; its low-frequency part is smooth
and compactly supported. Thus $\sup_{|t|\le R_0}\|\check b_t\|_1<\infty$.
Since $M_t=I+(2\pi)^{-d/2}\check b_t*$, it is bounded on every $L^r$,
$1\le r\le\infty$, uniformly in $|t|\le R_0$.
Therefore, by Theorem~\ref{Class-wave},
\[
 \begin{aligned}
 \|e^{it\langle D\rangle}(1-P_{\le0})f\|_{L_t^q(\rd\mu;L_x^r)}
 &\lesssim\|e^{it|D|}(1-P_{\le0})f\|_{L_t^q(\rd\mu;L_x^r)}\\
 &\lesssim\|(1-P_{\le0})f\|_{\dot H^s}
 \lesssim\|f\|_{H^s}.
 \end{aligned}
\]
Here we used $s=s_\alpha(q,r)\ge\frac{d+1}2(\frac12-\frac1r)\ge0$ for $\frac{d-1}2$-admissible pairs.
\end{proof}

\subsection{Scalar estimates for finite spatial exponents}
\begin{lem}[Scalar estimates by annular synthesis]\label{lem-scalar-annular-synthesis}
Let $a>0$, $a\ne1$, $0<\alpha\le1$, and let $\mu$ satisfy
\eqref{frac_meas}. Suppose that $2\le q\le\infty$, $2\le r<\infty$,
and
\[
 \frac d2\left(\frac12-\frac1r\right)\ge\frac\alpha q.
\]
Set $s=d/2-d/r-a\alpha/q$. Then the operator initially defined on
functions with compact Fourier support away from the origin by
\[
 \mathcal U_{a,s}g(t)=e^{it|D|^a}|D|^{-s}g
\]
extends to a bounded operator
\[
 \begin{gathered}
 \mathcal U_{a,s}:L^2(\mathbb R^d)
 \longrightarrow L_t^q(\rd\mu;L_x^r(\mathbb R^d)),\\
 \|\mathcal U_{a,s}g\|_{L_t^q(\rd\mu;L_x^r)}
 \le C_{d,a,\alpha,q,r}\Fmu^{1/q}\|g\|_2.
 \end{gathered}
\]
Here $\Fmu^{1/q}$ is interpreted as $1$ when $q=\infty$.
The condition $d/r+a/q<d$ of Theorem~\ref{Class-Schro} is not assumed;
in particular, $s\le-d/2$ is allowed. In that case the conclusion
concerns the extension from annular data, which agrees with the Fourier
definition of $e^{it|D|^a}|D|^{-s}g$ as a tempered distribution (see the
end of the proof).
\end{lem}
\begin{proof}
We use the frequency-localized estimate in
\cite[Proposition~2.9]{LLR2025}, which requires only \eqref{frac_adm}
with $\sigma=d/2$ and $(q,r)\ne(2,\infty)$; in particular, it covers
$q=2$ and does not require $d/r+a/q<d$. For a fixed smooth annular cutoff $Q_k$ at
frequency $2^k$, the estimate is
\begin{equation}\label{eq:scalar-localized-repair}
 \|e^{it|D|^a}Q_k h\|_{L_t^q(\rd\mu;L_x^r)}
 \le C\Fmu^{1/q}2^{ks}\|h\|_2,
 \qquad k\in\mathbb Z.
\end{equation}
The factor $\Fmu^{1/q}$ is explicit in the cited estimate; equivalently it follows by normalizing $\mu/\Fmu$. The same result holds for any fixed smooth
annular cutoff: its support can be covered by finitely many
rescaled cutoff regions in that proposition, with smooth
subordinate cutoffs, and the constants are independent of $k$.

Let $P_k$ be a homogeneous Littlewood--Paley partition, and let
$Q_k$ be a smooth cutoff equal to one on the support of $P_k$.
Put $v_k=e^{it|D|^a}|D|^{-s}P_kg$. Applying
\eqref{eq:scalar-localized-repair} to $h=|D|^{-s}P_kg$ gives
\[
 \|v_k\|_{L_t^q(\rd\mu;L_x^r)}
 \le C\Fmu^{1/q}2^{ks}\||D|^{-s}P_kg\|_2
 \le C_s\Fmu^{1/q}\|P_kg\|_2.
\]
For every finite set $\mathcal K\subset\mathbb Z$, the spatial
Littlewood--Paley inequality, followed by Minkowski's inequality
in $L_t^{q/2}(\rd\mu;L_x^{r/2})$, yields
\begin{align*}
 \left\|\sum_{k\in\mathcal K}v_k\right\|_{L_t^q(\rd\mu;L_x^r)}
 &\le C_r
 \left\|\left(\sum_{k\in\mathcal K}|v_k|^2\right)^{1/2}
 \right\|_{L_t^q(\rd\mu;L_x^r)}\\
 &\le C_r
 \left(\sum_{k\in\mathcal K}
 \|v_k\|_{L_t^q(\rd\mu;L_x^r)}^2\right)^{1/2}\\
 &\le C\Fmu^{1/q}
 \left(\sum_{k\in\mathcal K}\|P_kg\|_2^2\right)^{1/2}
 \le C\Fmu^{1/q}\|g\|_2.
\end{align*}
The second line follows from Minkowski's inequality in
$L_t^{q/2}L_x^{r/2}$, since $q/2,r/2\ge1$. The constants do not depend
on the size or location of $\mathcal K$. Applying the same
calculation to frequency tails proves convergence in the target
space, as the $\ell^2(L_x^2)$ norm of these tails tends to zero.

Finally, $s<d/2$ for all the pairs in the statement. For every
$g\in L^2$ and every Schwartz test function $\varphi$,
\[
 \int |\xi|^{-s}|\widehat g(\xi)|
       |\widehat\varphi(\xi)|\,d\xi
 \le \|g\|_2\,
       \||\xi|^{-s}\widehat\varphi(\xi)\|_2<\infty.
\]
Thus $|\xi|^{-s}\widehat g$ defines a tempered distribution.
The same majorant applies after multiplication by
$e^{it|\xi|^a}$, uniformly in $t$. The annular limits therefore
agree with this canonical Fourier definition as distributions;
for $\mu$-almost every $t$ they agree with the target-space limit.
This identifies the extension and completes the proof.
\end{proof}

\begin{rem}
On the sharp line, $d/r+a/q=d/2+(a-2\alpha)/q$.
When $d>2\alpha$, the additional hypothesis in
Theorem~\ref{Class-Schro} holds at $E_{d/2}$ precisely when
$a<d+2\alpha$. It holds at $A_{d/2}$ precisely
when $a<d+3\alpha$.
Lemma~\ref{lem-scalar-annular-synthesis} supplies all finite-$r$
scalar inputs needed by the sharp-line interpolation, for every
$a>0$, $a\ne1$. One may therefore keep the original condition
$d/r+a/q<d$ at the final theorem's target, without adding an
artificial restriction on $a$: internal interpolation points
outside that condition are interpreted through the $L^2$
operator in the lemma. The lemma concerns scalar estimates only;
it gives no orthonormal summation over frequency annuli.
\end{rem}

\section{Sharp admissible pairs}\label{sec-sharp}
We prove the sufficiency statements for sharp admissible pairs with
finite spatial exponent:
\[
 \frac\alpha q=\sigma\left(\frac12-\frac1r\right),\qquad
 2\le q\le\infty,\quad2\le r<\infty.
\]
The argument for $(A_\sigma,C]$ also applies when $\sigma<\alpha$, in which
case $E_\sigma$ is not defined. We first treat fractional Schr\"odinger
 evolution and then the wave and Klein--Gordon equations.
The boundary $r=\infty$ is considered in Section~\ref{sec-OD}.

\subsection{The fractional Schr\"odinger equation}\label{sec-sharp-frac}
We set
\[
\mathcal U_sf(x,t)=U_a(t)|D|^{-s}f,
\quad
U_a(t)=\e^{it|D|^a}, \quad a\in\mathbb{R}_{>0}\setminus\{1\}.
\]
This is the operator $\mathcal U_{a,s}$ of
Lemma~\ref{lem-scalar-annular-synthesis}. For finite $r$, that lemma supplies the
bounded map $\mathcal U_s:L^2\to L_t^q(\rd\mu;L_x^r)$ at every
admissible pair, with $s=s_{a\alpha}(q,r)$, including internal
interpolation points outside the extra hypothesis of
Theorem~\ref{Class-Schro}. Theorem~\ref{thm-frac-schro} assumes this
condition only at the target point.
\begin{prop}\label{prop-frac-schro}
Let $2\le q\le\infty$, $2\le r<\infty$, and suppose
\[
\frac{\alpha}{q}
=
\frac d2\left(\frac12-\frac1r\right).
\]
If \(\widetilde r>1+\frac d\alpha\), then for any \(W_1,W_2
\in
L_t^{\widetilde q}
\left(\mathrm d\mu;L_x^{\widetilde r}\right)\), we have
\begin{equation}\label{eq: Schatten-goal}
    \left\|
W_1 \mathcal U_s \mathcal U_s^*W_2
\right\|_{\mathfrak S^{\widetilde r}}
\le
C\Fmu^{2/q}
\|W_1\|_{L_t^{\widetilde q}(\mathrm d\mu;L_x^{\widetilde r})}
\|W_2\|_{L_t^{\widetilde q}(\mathrm d\mu;L_x^{\widetilde r})},
\end{equation}
where \(s=s_{a \alpha}(q,r)=
\frac d2-\frac dr-\frac{a\alpha}{q}.\)
\end{prop}
Since $\beta=\frac{2r}{r+2}$ satisfies $\beta'=\frac{2r}{r-2}=\widetilde r$, Propositions~\ref{prop-frac-schro} and \ref{prop_dual} give \eqref{fra-equa-ineq} with $\beta=\frac{2r}{r+2}=\beta_{d/2}(q,r)$. Note that, for $2\le r\le\infty$, the condition $\widetilde r>1+\frac d\alpha$ is equivalent to
\[
	2\le r < \frac{2(d+\alpha)}{d-\alpha},
\]
with the upper bound read as infinity when $d=\alpha$.
This corresponds to $(A_{d/2}, C]$. (When $r=2$, i.e. at $C$, the estimate is trivial.) Therefore, we obtain Theorem~\ref{thm-frac-schro} for $(1/r, 1/q)\in (A_{d/2}, C]$.

\begin{rem}\label{rem-varying-order}
When the Sobolev order changes along an interpolation segment, we
interpolate the initial-space operators \eqref{eq-initial-potential} with
$s(z)$ affine in $z$ and the usual mixed-norm powers of the potential. On
each boundary line the imaginary powers of $\Psi(D)$ act as separate
unitary operators on the initial $L^2$ space, on the left and on the
right, so they preserve the strong or weak Schatten bound available at
that endpoint. The finite-frequency construction is given in
Lemma~\ref{newq:varying-order}; the same argument applies to the
Lorentz sequence input used in Section~\ref{sec-nonsharp}.
\end{rem}

For the segment $[E_{d/2}, A_{d/2}]$ (assuming $\alpha\le d/2$), let $P\in(E_{d/2},A_{d/2})$, or $P=E_{d/2}$ if $d\neq2\alpha$. By Lemma~\ref{lem-scalar-annular-synthesis} and the triangle inequality, we have the trivial estimate
\begin{align}\label{frac-schro-trivial}
	\Big\| \sum_{j\in J} n_j | \mathcal U_{s_P} f_j |^2 \Big\|_{L_t^{q_P/2}(\mathrm d\mu; L_x^{r_P/2})} \lesssim \| \{n_j\}_{j\in J} \|_{\ell^1}
\end{align}
at $P=(1/r_P,1/q_P)$, with $s_P=s_{a\alpha}(q_P,r_P)$. Let $Q\in[P, A_{d/2}]$. Interpolating (see Remark~\ref{rem-varying-order}) \eqref{frac-schro-trivial} with the estimate at points $X\in(A_{d/2},C]$ close to $A_{d/2}$, we obtain \eqref{fra-equa-ineq} at $Q$ with $\beta$ given by $\frac1\beta=\frac{1-\theta}{1}+\frac{\theta}{\beta_X}$, where $Q=(1-\theta)P+\theta X$. Since $\frac1\beta$ and $\frac2q$ are affine along the segment, $1\ge\frac2{q_P}$ and $\frac1{\beta_X}>\frac2{q_X}$ for $X\in(A_{d/2},C)$, we have $\beta<q_Q/2$. Moreover, $\frac1\beta-\frac2{q_Q}=(1-\theta)(1-\frac2{q_P})+\theta(\frac1{\beta_X}-\frac2{q_X})\to0$ as $P\to E_{d/2}$ and $X\to A_{d/2}$, so every $\beta<q_Q/2$ is obtained. This proves Theorem~\ref{thm-frac-schro} (2) on $[E_{d/2}, A_{d/2}]\setminus\{D\}$.
It remains to prove Proposition~\ref{prop-frac-schro}.

\begin{proof}[Proof of Proposition~\ref{prop-frac-schro}]
We may assume $r>2$. Choose \(\nu\) such that
\[
    0<\nu<\alpha,
\quad
0<h:=d-\nu(\widetilde r-2)<\alpha,
\]
that is,
\(\max\{0,\frac{d-\alpha}{\widetilde r-2}\}<\nu<
\min\{\alpha,\frac d{\widetilde r-2}\}.\)
Since $d\ge\alpha$, such $\nu$ exists if and only if \(\widetilde r>1+d/\alpha\).

For $0\le \operatorname{Re}z\le1$ and $0<\delta<R$, $0<\varepsilon<1$, define
\[
T_{z,\varepsilon}^{\delta,R}F(t)
=
\int_{\delta<|t-\tau|<R}
|t-\tau|^{\nu(\widetilde r z/2-1)}
U_a(t-\tau)
|D|^{bz}
m_\varepsilon(D)F(\tau)
\,\mathrm d\mu(\tau),
\]
where \(b={(a-2)d}/{2}\). Since $|\xi|^{bz}m_\varepsilon(\xi)$ is smooth and compactly supported away from the origin and the time integration is restricted to $\delta<|t-\tau|<R$, this is an analytic family of bounded operators on $L^2(\mathrm d\mu\,\mathrm dx)$ of admissible growth. When \(z=2/\widetilde{r}\), we have
\[
T_{2/\widetilde r,\varepsilon}^{\delta,R}F(t)
=
\int_{\delta<|t-\tau|<R}
U_a(t-\tau)
|D|^{-2s}
m_\varepsilon(D)F(\tau)
\,\mathrm d\mu(\tau),
\]
because $b\cdot\frac2{\widetilde r}=\frac{(a-2)d}{\widetilde r}=-2s_{a\alpha}(q,r)$ by the sharp admissibility $\frac\alpha q=\frac d{2\widetilde r}$.

We first prove the estimate when \(W_1, W_2\) are bounded simple functions with bounded
support. Let $z=iy$.  Since \(U_a(t)\) is unitary on \(L^2_x\),
\(
\|
U_a(t-\tau)|D|^{iby}m_\varepsilon(D)
\|_{L_x^2\to L_x^2}
\le 1.
\)
Therefore, for any $F,G\in L_t^2L_x^2$, we have
\[
\begin{aligned}
 \left|\langle W_1T_{iy,\varepsilon}^{\delta,R}W_2F,G\rangle\right|
 &\le\iint |t-\tau|^{-\nu}F_1(t)F_2(\tau)\,\rd\mu(t)\rd\mu(\tau),\\
 F_1(t)&=\|W_1(t)\|_\infty\|G(t)\|_2,\qquad
 F_2(\tau)=\|W_2(\tau)\|_\infty\|F(\tau)\|_2.
\end{aligned}
\]
Choose $v_0=2\alpha/(\alpha-\nu)$ and
$u_\nu=2\alpha/(2\alpha-\nu)$ as in \eqref{eq: hls}. Then
\[
 \frac1{u_\nu}=\frac1{v_0}+\frac12.
\]
Thus H\"older places each product in $L_t^{u_\nu}$, and applying
\eqref{eq: hls} with $\lambda=\nu$ gives
\begin{equation}\label{eq: Schatten-infty}
    \left\|
W_1T_{iy,\varepsilon}^{\delta,R}W_2
\right\|_{\mathfrak S^\infty}
\le
C \Fmu^{\nu/\alpha}
\|W_1\|_{L_t^{v_0}(\mathrm d\mu; L_x^\infty)}
\|W_2\|_{L_t^{v_0}(\mathrm d\mu; L_x^\infty)},
\end{equation}
where \(v_0=2\alpha/{(\alpha-\nu)}\), and
the constant is independent of $y,\varepsilon,\delta,R$.

Let $z=1+iy$. By Proposition~\ref{prop_disp},
\[
\left|
\int_{\mathbb R^d}
\e^{i(x\cdot\xi+(t-\tau)|\xi|^a)}
|\xi|^{b+iby}
m_\varepsilon(\xi)\,\mathrm d\xi
\right|
\le
C(1+|y|)^M|t-\tau|^{-d/2}
\]
uniformly in \(\varepsilon\).
Consequently, the kernel of \(T_{1+iy,\varepsilon}^{\delta,R}\) satisfies
\begin{equation}\label{eq: Hilbert--Schmidt ker}
        |K_{1+iy}(t,x;\tau,x')|
 \leq C(1+|y|)^M|t-\tau|^{\nu(\widetilde{r}/2-1)-d/2}
 =C(1+|y|)^{M}|t-\tau|^{-h/2}.
\end{equation}
Set $v_1=4\alpha/(2\alpha-h)$. For the HLS exponent
$u_h=2\alpha/(2\alpha-h)$ one has $u_h=v_1/2$.
The squared spatial $L^2$ norms therefore belong to $L_t^{u_h}$.
Using \eqref{eq: Hilbert--Schmidt ker} and \eqref{eq: hls} with
$\lambda=h$, we obtain
\begin{equation*}
 \begin{aligned}
&\left\|
W_1T_{1+iy,\varepsilon}^{\delta,R}W_2
\right\|_{\mathfrak S^2}^2\\
&\quad\le
C(1+|y|)^{2M}
\iint |t-\tau|^{-h}
\|W_1(t)\|_2^2
\|W_2(\tau)\|_2^2
\,\mathrm d\mu(t)\,\mathrm d\mu(\tau)
\\
&\quad\le
C(1+|y|)^{2M}
\Fmu^{h/\alpha}
\|W_1\|_{L_t^{v_1}L_x^2}^2
\|W_2\|_{L_t^{v_1}L_x^2}^2,
\end{aligned}
\end{equation*}
where \(v_1={4\alpha}/{(2\alpha-h)}\), that is,
\begin{equation}\label{eq: Schatten-2}
\left\|
W_1T_{1+iy,\varepsilon}^{\delta,R}W_2
\right\|_{\mathfrak S^2}
\le
C(1+|y|)^{M}
\Fmu^{h/(2\alpha)}
\|W_1\|_{L_t^{v_1}L_x^2}
\|W_2\|_{L_t^{v_1}L_x^2}.
\end{equation}
Let \(\theta=2/{\widetilde r}\). Stein's interpolation theorem, applied to the analytic family $z\mapsto W_{1,z}T_{z,\varepsilon}^{\delta,R}W_{2,z}$ with the usual choice of $W_{j,z}$ adapted to mixed norms (cf. \cite{FraSa2017, BLN2021}), together with \eqref{eq: Schatten-infty} and \eqref{eq: Schatten-2} gives
\begin{equation}\label{eq: Schatten-r-tilde-trunc}
    \left\|W_1
T_{2/\widetilde r,\varepsilon}^{\delta,R}W_2
\right\|_{\mathfrak S^{\widetilde r}}
\le
C \Fmu^{2/q}
\|W_1\|_{L_t^{\widetilde q}L_x^{\widetilde r}}
\|W_2\|_{L_t^{\widetilde q}L_x^{\widetilde r}}.
\end{equation}
Indeed, $\frac{1-\theta}\infty+\frac\theta2=\frac1{\widetilde r}$ and, using $h=d-\nu(\widetilde r-2)$ and $\frac1q=\frac{d}{2\alpha\widetilde r}$,
\[
\begin{aligned}
 \frac{1-\theta}{v_0}+\frac\theta{v_1}
 &=\frac{(\widetilde r-2)(\alpha-\nu)+2\alpha-h}{2\alpha\widetilde r}\\
 &=\frac12-\frac d{2\alpha\widetilde r}
 =\frac12-\frac1q=\frac1{\widetilde q},\\
 (1-\theta)\frac\nu\alpha+\theta\frac h{2\alpha}
 &=\frac d{\alpha\widetilde r}=\frac2q.
\end{aligned}
\]

For fixed $\varepsilon$, the multiplier \(|\xi|^{-2s}m_\varepsilon(\xi)\) is smooth and compactly supported away from the origin, so the kernel of $U_a(t-\tau)|D|^{-2s}m_\varepsilon(D)$ is bounded. The Frostman condition implies that $\mu$ has no atoms. Thus
\(
\mathbf 1_{\{\delta<|t-\tau|<R\}}
\to1
\)
for \(\mu\times\mu\)-almost every \((t,\tau)\)
 as
\(\delta\to0\) and \(R\to\infty\). Since \(W_1,W_2\) are bounded with bounded support and $\mu$ is finite, dominated convergence gives
\[
W_1T_{2/\widetilde r,\varepsilon}^{\delta,R}W_2
\longrightarrow
W_1\mathcal U_sm_\varepsilon(D)\mathcal U_s^*W_2
\]
in the weak operator topology. By the lower semicontinuity of the Schatten norm with respect to weak operator convergence, we have
\[
\bigl\|
W_1 \mathcal U_sm_\varepsilon(D) \mathcal U_s^*W_2
\bigr\|_{\mathfrak S^{\widetilde r}}
\le
C\Fmu^{2/q}
\|W_1\|_{L_t^{\widetilde q}(\mathrm d\mu;L_x^{\widetilde r})}
\|W_2\|_{L_t^{\widetilde q}(\mathrm d\mu;L_x^{\widetilde r})}.
\]
Since \(m_\varepsilon(\xi)\to1\)
for every \(\xi\ne0\) and \(0\le m_\varepsilon\le1\), we have
\(
m_\varepsilon(D)\to I\) strongly on \(L_x^2
\). By the boundedness of \(\mathcal U_s:L_x^2\to
L_t^q(\mathrm d\mu;L_x^r)\) and H\"older's inequality, \(W_1 \mathcal U_s:L_x^2\to L_t^2L_x^2\) and \( \mathcal U_s^*W_2: L_t^2L_x^2\to L_x^2\) are bounded. It follows that
\(
W_1 \mathcal U_sm_\varepsilon(D) \mathcal U_s^*W_2
\to
W_1 \mathcal U_s \mathcal U_s^*W_2
\)
strongly on \(L_t^2L_x^2\). Applying lower semicontinuity once again, we obtain \eqref{eq: Schatten-goal}
for bounded simple functions \(W_1, W_2\) with bounded support. Finally, the general case follows by a standard density argument.
\end{proof}

\begin{rem}
The time truncation $\delta<|t-\tau|<R$ avoids the singularity at $t=\tau$ and ensures that the operators are well defined throughout the interpolation argument. The frequency cutoff $m_\varepsilon$ removes the possible singular behavior of $|\xi|^{bz}$ near $\xi=0$ and at infinity, so that $|\xi|^{bz}m_\varepsilon(\xi)$ is a smooth multiplier compactly supported away from the origin.
\end{rem}

\subsection{The wave and Klein--Gordon equations}\label{sec-sharp-wave}
For the wave and the Klein--Gordon equations, we set
\begin{align*}
\mathcal U_s^{\textrm{wave}} f(x,t)&=\e^{it|D|} |D|^{-s} f(x),
\\
\mathcal U_s^{\textrm{KG}} f(x,t) &=  \e^{it\langle D\rangle} \langle D\rangle^{-s} f (x),
\end{align*}
and write $\mathcal U_s^{(\cdot)}$ with $\cdot \in \{\mathrm{wave}, \mathrm{KG}\}$. By Theorem~\ref{Class-wave} and Lemma~\ref{K-L-Strichartz}, $\mathcal U_s^{(\cdot)} : L_x^2(\mathbb R^d) \to L_t^q(\mathrm d \mu; L_x^r(\mathbb R^d))$ is bounded for $s=s_\alpha(q,r)$ and $\frac{d-1}2$-admissible $(q,r)$ with $2\le q\le\infty$, $2\le r<\infty$.

\begin{cor}\label{cor-wave-kg}
Let $d\ge2$ and let $2\le q\le \infty$, $2\le r<\infty$ satisfy
\[
\frac{\alpha}{q}
=
\frac {d-1}2\left(\frac12-\frac1r\right).
\]
If \(\widetilde r>1+\frac {d-1}\alpha\), then for any \(W_1,W_2
\in
L_t^{\widetilde q}
\left(\mathrm d\mu;L_x^{\widetilde r}\right)\), we have
\begin{equation}\label{eq: Schatten-goal-wave-kg}
    \left\|
W_1 \mathcal U_s^{(\cdot)} (\mathcal U_s^{(\cdot)})^*W_2
\right\|_{\mathfrak S^{\widetilde r}}
\le
C\Fmu^{2/q}
\|W_1\|_{L_t^{\widetilde q}(\mathrm d\mu;L_x^{\widetilde r})}
\|W_2\|_{L_t^{\widetilde q}(\mathrm d\mu;L_x^{\widetilde r})},
\end{equation}
where \(s=s_{\alpha}(q,r)=
\frac d2-\frac dr-\frac{\alpha}{q}.\)
\end{cor}
\begin{proof}
The case $r=2$ is the energy estimate, so assume $r>2$.
The proof is the same as that of Proposition~\ref{prop-frac-schro}, with the following changes.
Choose \(\nu\) such that
\(
    0<\nu<\alpha\) and
\(
0<h:=d-1-\nu(\widetilde r-2)<\alpha,
\)
which is possible if and only if \(\widetilde r>1+\frac{d-1}\alpha\), since $d-1\ge\alpha$ when $d\ge2$.
We set $b=-\frac{d+1}2$; then $b\cdot\frac2{\widetilde r}=-\frac{d+1}{\widetilde r}=-2s_\alpha(q,r)$ by the sharp admissibility $\frac\alpha q=\frac{d-1}{2\widetilde r}$. We consider the analytic family
\[
\begin{split}
\widetilde T_{z,\varepsilon}^{\delta,R}F(t)
&=
(z-1)\int_{\delta<|t-\tau|<R}
|t-\tau|^{\nu(\widetilde r z/2-1)}\\
&\qquad\times\e^{i(t-\tau)\phi(D)}
\Psi(D)^{bz}
m_\varepsilon(D)F(\tau)
\,\mathrm d\mu(\tau),
\end{split}
\]
with $(\phi,\Psi)=(|\xi|,|\xi|)$ or $(\langle\xi\rangle,\langle\xi\rangle)$.
For $z=iy$, the factor $|z-1|\le1+|y|$ only affects the admissible growth, and \eqref{eq: Schatten-infty} holds with the additional factor $1+|y|$.
For $z=1+iy$, Proposition~\ref{wave-decay} with $\kappa=by$ gives
\[
\begin{aligned}
 &|z-1|\left|\int_{\mathbb R^d}
 e^{i(x\cdot\xi+(t-\tau)\phi(\xi))}
 \Psi(\xi)^{b+iby}m_\varepsilon(\xi)\,\rd\xi\right|\\
 &\qquad=\frac{|by|}{|b|}\left|\int\cdots\right|
 \le C(1+|y|)^M|t-\tau|^{-(d-1)/2}
\end{aligned}
\]
uniformly in $\varepsilon$; the factor $z-1$ compensates the factor $|\kappa|$ in Proposition~\ref{wave-decay}. Then, with $v_0,v_1$ defined as before (with the new $h$) and $\theta=2/\widetilde r$, we have
\[
\begin{aligned}
\frac{1-\theta}{v_0}+\frac{\theta}{v_1}
&=\frac12-\frac{d-1}{2\alpha\widetilde r}
=\frac12-\frac1q=\frac1{\widetilde q},\\
(1-\theta)\frac\nu\alpha+\theta\frac h{2\alpha}
&=\frac{d-1}{\alpha\widetilde r}=\frac2q,
\end{aligned}
\]
and interpolation gives \eqref{eq: Schatten-r-tilde-trunc} for $\widetilde T_{2/\widetilde r,\varepsilon}^{\delta,R}$. Since $2/\widetilde r-1\neq0$, we may divide by this factor. The truncations are removed as before, using the boundedness of $\mathcal U_s^{(\cdot)}:L^2_x\to L_t^q(\mathrm d\mu;L_x^r)$. In the Klein--Gordon case, the multiplier $\langle\xi\rangle^{bz}$ has no singularity at the origin, but we keep $m_\varepsilon$ for uniformity.
\end{proof}

Using Corollary~\ref{cor-wave-kg}, we prove Theorems~\ref{thm-wave}, \ref{thm-KG} for sharp admissible pairs.
For $2\le r<\infty$, the condition $\widetilde r>1+\frac {d-1}\alpha$ is equivalent to
\(
	2\le r < \frac{2(d-1+\alpha)}{d-1-\alpha}
\)
(read as $r<\infty$ if $d-1=\alpha$),
which corresponds to $(A_{\frac{d-1}2}, C]$.
By Proposition~\ref{prop_dual}, we obtain Theorems~\ref{thm-wave}, \ref{thm-KG} for $(1/r, 1/q)\in (A_{\frac{d-1}2}, C]$ with $\beta=\frac{2r}{r+2}$ and $s=s_\alpha(q,r)$. For the Klein--Gordon equation write $s_0=s_\alpha(q,r)$.
If $s\ge s_0$, set $C_s=\langle D\rangle^{s_0-s}$, so
$\|C_s\|_{L^2\to L^2}\le1$. The initial-space potential operators satisfy
$B_s(V)=C_sB_{s_0}(V)C_s$. The Schatten ideal inequality preserves
the dual bound. Applying duality to the orthonormal vectors
$g_j=\langle D\rangle^sf_j$ proves the estimate for the original
$H^s$-orthonormal family $(f_j)$. This argument applies to all the
Klein--Gordon estimates below as well.
For $[E_{\frac{d-1}2}, A_{\frac{d-1}2}]\setminus\{D\}$ (assuming $\alpha\le\frac{d-1}2$), we argue as in the case of the fractional Schr\"odinger equation, using the trivial estimate
\begin{align}\label{frac-wave-kg-trivial}
	\Bigg\| \sum_{j\in J} n_j | \mathcal U_{s_P}^{(\cdot)} f_j |^2 \Bigg\|_{L_t^{q_P/2}(\mathrm d\mu; L_x^{r_P/2})} \lesssim \| \{n_j\}_{j\in J} \|_{\ell^1},\quad s_P=s_\alpha(q_P,r_P),
\end{align}
at points $P=(1/r_P,1/q_P)\in(E_{\frac{d-1}2},A_{\frac{d-1}2})$ (or $P=E_{\frac{d-1}{2}}$, where $(q_P,r_P)=(2,\frac{2(d-1)}{d-1-2\alpha})$, if $d-1\neq2\alpha$), which follows from Theorem~\ref{Class-wave}, Lemma~\ref{K-L-Strichartz} and the triangle inequality.

\section{Non-sharp admissible pairs}\label{sec-nonsharp}

In this section, we prove the results for non-sharp admissible pairs $(q,r)$ off the segment $(O,D)$, which is treated in Section~\ref{sec-OD}.
Let $(\phi,\Psi,a,\sigma)$ be one of \eqref{three-cases}, let $U_\phi(t)=\e^{it\phi(D)}$, and set
\[
s(q,r)=\frac{d}{2}-\frac{d}{r}-\frac{a\alpha}{q},
\]
so that $s=s_{a\alpha}(q,r)$ for the fractional Schr\"odinger equation and $s=s_\alpha(q,r)$ for the wave and Klein--Gordon equations. Note that $s(q,r)$ is affine in $(1/r,1/q)$.

\begin{prop}\label{non-sharp} Suppose that $\mu$ satisfies \eqref{frac_meas} and that the estimate
\begin{equation}\label{sharp}
    \Bigg\|\sum_{j\in J}n_j|U_\phi\Psi(D)^{-s(q,r)}f_j|^2\Bigg\|_{L_t^{\frac{q}{2}}(\mathrm d\mu; L_x^{\frac{r}{2}})}\lesssim \|\{n_j\}_{j\in J}\|_{\ell^\beta}
\end{equation}
holds for all orthonormal families $\{f_j\}_{j\in J}$ in $L^2(\mathbb{R}^d)$ with $\beta=\beta_{\sigma}(q,r)$ whenever $(\frac{1}{r},\frac{1}{q})\in (A_\sigma,C]$. Then \eqref{sharp} holds with $\beta=\beta_{\sigma}(q,r)$ for all $(\frac{1}{r},\frac{1}{q})\in \mathrm{int}\,OA_\sigma C$.
\end{prop}

\begin{proof} We use the following orthonormal Sobolev inequality, which follows from Lieb's inequality \cite[Theorem~3]{EHL1983} by decomposing the coefficients into dyadic level sets: for $1\le p<\infty$,
\begin{equation}\label{Sobolev ine}
    \Bigg\|\sum_{j\in J}n_j|\Psi(D)^{-\frac{d}{2p'}}f_j|^2\Bigg\|_{L^p(\mathbb{R}^d)}\lesssim\|\{n_j\}_{j\in J}\|_{\ell^{p,1}}
\end{equation}
for all orthonormal families $\{f_j\}_{j\in J}$ in $L^2(\mathbb{R}^d)$.
Since $\{U_\phi(t) f_j\}_{j\in J}$ is an orthonormal family in $L^2{(\mathbb{R}^d)}$ for each fixed $t$ and $s(\infty,r)=\frac{d}{2}-\frac dr=\frac{d}{2(r/2)'}$, \eqref{Sobolev ine} with $p = \frac{r}{2}$ implies that, for $(\frac{1}{r},\frac{1}{q})\in (O,C)$,
\begin{equation}\label{OC}
    \Bigg\|\sum_{j\in J}n_j|U_\phi\Psi(D)^{-s(\infty,r)}f_j|^2\Bigg\|_{L_t^{\infty}(\mathrm d\mu; L_x^{\frac{r}{2}})}\lesssim \|\{n_j\}_{j\in J}\|_{\ell^{\frac{r}{2},1}}.
\end{equation}
Note that $\beta_\sigma(\infty,r)=r/2$. Every point of $\mathrm{int}\,OA_\sigma C$ lies on a segment joining a point of $(A_\sigma,C)$ and a point of $(O,C)$. We interpolate \eqref{sharp} on $(A_\sigma,C)$ with \eqref{OC}, using that $1/\beta_\sigma(q,r)$ and $s(q,r)$ are affine. Since one endpoint bound is of weak Schatten type, complex interpolation (see Remark~\ref{rem-varying-order}) gives the estimate with $\ell^{\beta_\sigma(q,r),u}$ on the right-hand side for some $u\ge1$; since $\ell^{\beta,1}\subset\ell^{\beta,u}$, we obtain for all $(\frac{1}{r},\frac{1}{q})\in \text{int}\, OA_\sigma C$
\begin{equation}\label{beta1}
    \Bigg\|\sum_{j\in J}n_j|U_\phi\Psi(D)^{-s(q,r)}f_j|^2\Bigg\|_{L_t^{\frac{q}{2}}(\mathrm d\mu; L_x^{\frac{r}{2}})}\lesssim \|\{n_j\}_{j\in J}\|_{\ell^{\beta_{\sigma}(q,r),1}}.
\end{equation}

To upgrade \eqref{beta1} to the strong-type estimate, we fix $(\frac{1}{r},\frac{1}{q})\in \text{int}\, OA_\sigma C$ and choose two distinct points $(\frac{1}{r_i},\frac{1}{q_i})\in \text{int}\, OA_\sigma C$, $i=0,1$, such that $(\frac1r,\frac1q)$ lies on the open segment joining them and $s(q_i,r_i)=s(q,r)$. Then $q_0\neq q_1$, $r_0\ne r_1$, and also $\beta_\sigma(q_0,r_0)\neq\beta_\sigma(q_1,r_1)$, since the level lines $dx+a\alpha y=\mathrm{const}$ of $s$ and $2\sigma x+\alpha y=\mathrm{const}$ of $\sigma/\beta_\sigma$, with $x=1/r$, $y=1/q$, are not parallel ($2a\sigma\neq d$ in all three cases). For fixed $s$, the map $\{n_j\}\mapsto\sum_jn_j|U_\phi\Psi(D)^{-s}f_j|^2$ is linear. Hence real interpolation with second index $q/2$, using \eqref{real-L} and \eqref{real-ell}, gives
\begin{equation*}
    \Bigg\|\sum_{j\in J}n_j|U_\phi\Psi(D)^{-s(q,r)}f_j|^2\Bigg\|_{L_t^{\frac{q}{2}}(\mathrm d\mu; L_x^{\frac{r}{2}, \frac{q}{2}})}\lesssim \|\{n_j\}_{j\in J}\|_{\ell^{\beta_{\sigma}(q,r), \frac{q}{2}}}.
\end{equation*}
If $(\frac1r,\frac1q)\in\mathrm{int}\,OA_\sigma B_\sigma$, then
  \begin{equation*}
      \beta_{\sigma}(q,r)< \frac{q}{2}\le \frac{r}{2}.
  \end{equation*}
Indeed, $\beta_{\sigma}(q,r)<\frac{q}{2}$ holds below the line through $O$ and $A_\sigma$, and $q\leq r$ holds above the diagonal $\frac{1}{q}=\frac{1}{r}$, which passes through $B_\sigma$. Thus the embeddings $L_x^{\frac{r}{2}, \frac{q}{2}}\subseteq L_x^{\frac{r}{2}}$ and $\ell^{\beta_\sigma(q,r)}\subseteq \ell^{\beta_\sigma(q,r), \frac{q}{2}}$ show that \eqref{sharp} holds with $\beta=\beta_\sigma(q,r)$ for all $(\frac{1}{r},\frac{1}{q})\in \text{int}\,OA_\sigma B_\sigma$.

At $C$, i.e. $(q,r)=(\infty,2)$, we have $s(\infty,2)=0$ and $\beta_\sigma(\infty,2)=1$, and the triangle inequality together with the unitarity of $U_\phi(t)$ gives
\begin{equation*}
\Bigg\|\sum_{j\in J}n_j|U_\phi f_j|^2\Bigg\|_{L_t^{\infty}(\mathrm d\mu; L_x^{1})}\le \|\{n_j\}_{j\in J}\|_{\ell^1}.
\end{equation*}
Every point of $\mathrm{int}\,OA_\sigma C$ lies on a segment joining $C$ and a point of $\mathrm{int}\,OA_\sigma B_\sigma$. Interpolating (see Remark~\ref{rem-varying-order}), we obtain \eqref{sharp} with $\beta=\beta_\sigma(q,r)$ for all $(\frac{1}{r},\frac{1}{q})\in\text{int}\,OA_\sigma C$.
\end{proof}

\begin{proof}[Proof of Theorems~\ref{thm-frac-schro}, \ref{thm-wave} and \ref{thm-KG} for non-sharp pairs off $(O,D)$]
For Theorems~\ref{thm-frac-schro}, \ref{thm-wave} and \ref{thm-KG}, put $\sigma=\frac d2$ in the fractional Schr\"odinger case and $\sigma=\frac{d-1}2$ in the wave and Klein--Gordon cases. By Proposition~\ref{prop-frac-schro}, Corollary~\ref{cor-wave-kg} and Proposition~\ref{prop_dual}, the hypothesis of Proposition~\ref{non-sharp} is satisfied. (For the Klein--Gordon equation we take $s=s_\alpha(q,r)$; the case $s>s_\alpha(q,r)$ follows as in Section~\ref{sec-sharp-wave}.) Hence part (1) of the theorems holds on $\mathrm{int}\,OA_\sigma C$.

For part (2), assume $\alpha\le\sigma$. For $(\frac1r,\frac12)\in(D,E_\sigma]$ (with $E_\sigma$ excluded if $E_\sigma=D$, i.e. $\sigma=\alpha$), the triangle inequality and Lemma~\ref{lem-scalar-annular-synthesis}, Theorem~\ref{Class-wave} or Lemma~\ref{K-L-Strichartz} give
    \begin{equation}\label{trivial-DE}
    \Bigg\|\sum_{j\in J}n_j\big|U_\phi\Psi(D)^{-s(2,r)}f_j\big|^2\Bigg\|_{L_t^{1}(\mathrm d\mu; L_x^{\frac{r}{2}})}\lesssim\|\{n_j\}_{j\in J}\|_{\ell^1},
\end{equation}
that is, \eqref{sharp} with $\beta=1=q/2$. On $(O,D)$, we have \eqref{sharp} for every $\beta<q/2$ by Proposition~\ref{prop-boundary} below, and on $[E_\sigma,A_\sigma]\setminus\{D\}$ by Section~\ref{sec-sharp}. Let $P\in ODE_\sigma A_\sigma\setminus\{O,D\}$ be a point not on these segments. For $X\in\mathrm{int}\,OA_\sigma C$ close to the line through $O$ and $A_\sigma$, the line through $X$ and $P$ meets $(D,E_\sigma]\cup(O,D)\cup[E_\sigma,A_\sigma]$ at a point $Y$ with $P=(1-\theta)Y+\theta X$, $\theta\in(0,1)$. Interpolating (see Remark~\ref{rem-varying-order}) the estimates at $X$ (with $\beta_\sigma(X)$) and $Y$ (with $\beta_Y\le q_Y/2$ arbitrarily close to $q_Y/2$), we obtain \eqref{sharp} at $P$ with $\frac1\beta=\frac{1-\theta}{\beta_Y}+\frac{\theta}{\beta_\sigma(X)}$. Since
\[
\frac1\beta-\frac2{q_P}=(1-\theta)\Big(\frac1{\beta_Y}-\frac2{q_Y}\Big)+\theta\Big(\frac1{\beta_\sigma(X)}-\frac2{q_X}\Big)
\]
is positive and can be made arbitrarily small, every $\beta<q_P/2$ is obtained. The same argument, with $Y\in(O,(0,\frac14))$, gives the statement in Remark~\ref{rem-small-sigma} when $\sigma<\alpha$.
\end{proof}

\section{The boundary \texorpdfstring{$r=\infty$}{r=infinity}}\label{sec-OD}

In this section, we deal with the boundary case $r=\infty$, i.e. $1/r =0$, $1/q \in (0, 1/2)$, and with the low-dimensional cases of Remark~\ref{rem-small-sigma}.

\subsection{Estimates on \texorpdfstring{$r=\infty$}{r=infinity}}\label{sec-OD-main}
In the proofs of Proposition~\ref{prop-boundary} and Lemma~\ref{lem-three-bounds} below, we normalize $\Fmu=1$ (note that $\Fmu>0$ since $\mu\ne0$); replacing $\mu$ by $\mu/\Fmu$ and rescaling the density time norm restores the factor $\Fmu^{2/q}$ in every conclusion.
Let $(\phi,\Psi,a,\sigma)$ be one of \eqref{three-cases}, $U_\phi(t)=\e^{it\phi(D)}$, and
\[
s(q)=s(q,\infty)=\frac d2-\frac{a\alpha}q .
\]
\begin{lem}[Scalar input on $r=\infty$]\label{lem-scalar-OD-completion}
Let $2<q<\infty$ and $2\alpha/q<\sigma$ in one of the cases
\eqref{three-cases}. Then, at $s(q)=d/2-a\alpha/q$,
\[
 \|U_\phi(t)\Psi(D)^{-s(q)}g\|_{L_t^q(\rd\mu;L_x^\infty)}
 \le C\Fmu^{1/q}\|g\|_2.
\]
The constant is uniform in finite annular cutoffs. In the
Klein--Gordon case it may depend on the bounded time support.
\end{lem}
\begin{proof}
Put $\gamma=d-2s(q)=2a\alpha/q$ and $\lambda=\gamma/a$.
The strict inequality $0<\gamma<a\sigma$ permits absolute summation
of the annular bounds \eqref{dyadic-disp}; see Lemma~\ref{lem-full-decay}.
For nonnegative finite frequency cutoffs $m_N\le1$ obtained from
the partition, the kernel of
$T_NT_N^*$, where $T_N=U_\phi\Psi(D)^{-s(q)}m_N(D)^{1/2}$,
is uniformly bounded by $C|t-\tau|^{-\lambda}$.
The bounded low-frequency KG kernel is absorbed on the bounded time
interval. Since $\lambda=2\alpha/q<\alpha$ and
$\lambda/\alpha+2/q'=2$, Lemma~\ref{lem_HLS} gives
\[
 \|T_N^*F\|_2^2\le C\Fmu^{2/q}
        \|F\|_{L_t^{q'}(\rd\mu;L_x^1)}^2.
\]
Since the truncated outputs are continuous in $x$ and vanish at
infinity, their $L_x^\infty$ norms are attained by testing against
$L_x^1$ functions, and time measurability follows from suprema over
rational points; hence the $TT^*$ argument applies on finite-annulus
data. Uniform bounds applied to initial-data tails give the
$L^2$ completion. In the homogeneous cases $s(q)<d/2$, and the
Fourier distribution argument in
Lemma~\ref{lem-scalar-annular-synthesis} identifies this completion
with $U_\phi\Psi(D)^{-s(q)}g$. As in Lemma~\ref{lem-scalar-annular-synthesis},
the condition $d/r+a/q<d$ is not needed at auxiliary interpolation points.
\end{proof}

\begin{prop}\label{prop-boundary}
Assume that $\mu$ satisfies \eqref{frac_meas}. Let $2<q<\infty$, and assume in addition $q>4$ if $\sigma<\alpha$. Then the estimate
\begin{equation}\label{boundary}
    \Bigg\|\sum_{j\in J}n_j|U_\phi\Psi(D)^{-s(q)}f_j|^2\Bigg\|_{L_t^{\frac{q}{2}}(\mathrm d\mu; L_x^\infty)}\lesssim \|\{n_j\}_{j\in J}\|_{\ell^\beta}
\end{equation}
   holds for all orthonormal families $\{f_j\}_{j\in J}$ in $L^2(\mathbb{R}^d)$ and every $\beta<\frac{q}{2}$.
\end{prop}

Theorems~\ref{thm-frac-schro}, \ref{thm-wave}, \ref{thm-KG} on $(O,D)$ (and the corresponding part of Remark~\ref{rem-small-sigma}) follow immediately from Proposition~\ref{prop-boundary}, since $s(q)=s_{a\alpha}(q,\infty)$ for the fractional Schr\"odinger equation and $s(q)=s_{\alpha}(q,\infty)$ for the wave and Klein--Gordon equations. (For the Klein--Gordon equation with $s>s_\alpha(q,\infty)$ we argue as in Section~\ref{sec-sharp-wave}.)

Since $r=\infty$, we have $\widetilde r=2$, and $\beta<q/2$ if and only if $\beta'>\widetilde q/2$. By Proposition~\ref{prop_dual} and the inclusions $\mathfrak S^{p_1}\subset\mathfrak S^{p_2}$, $p_1\le p_2$, it suffices to show that, for every $p>\widetilde q/2$ sufficiently close to $\widetilde q/2$,
\begin{equation}\label{goal-bound}
    \|W_1U_\phi\Psi(D)^{-2s(q)}U_\phi^*W_2\|_{\mathfrak S^{p}}\lesssim\prod_{i=1}^2\|W_i\|_{L_t^{\widetilde q}(\rd\mu;L_x^2)}.
\end{equation}
For $\kappa\in\mathbb R$ we write
\[
\begin{aligned}
 \mathcal T_q^\kappa(W_1,W_2)
 &=W_1U_\phi\Psi(D)^{-2s(q)+i\kappa}U_\phi^*W_2,\\
 \mathcal T_{q,k}^\kappa(W_1,W_2)
 &=W_1U_\phi\Psi(D)^{-2s(q)+i\kappa}P_kU_\phi^*W_2,
\end{aligned}
\]
and $\mathcal T_{q,\le0}^\kappa$ is defined with $P_{\le0}$ in place of $P_k$. The kernel of $\mathcal T_{q,k}^{\kappa}$ is $W_1(t,x)K_k^{s(q),\kappa}(t-\tau,x-y)W_2(\tau,y)$, where $K_k^{s,\kappa}$ is as in Section~\ref{sec-disp}. In what follows, $W_1,W_2$ are bounded simple functions with bounded support, and all bounds are first proved with $\Psi(D)^{-2s(q)+i\kappa}\sum_{k\in\mathcal K}P_k$ in place of $\Psi(D)^{-2s(q)+i\kappa}$, uniformly over finite sets $\mathcal K\subset\mathbb Z$. The general case then follows by lower semicontinuity of Schatten norms and density.

\begin{lem}\label{lem-three-bounds}
Let $\mu$ satisfy \eqref{frac_meas}. For each fixed choice of exponents
and flow parameters, there are constants $C,M>0$ such that the following
hold for all $\kappa\in\mathbb R$, uniformly in the frequency cutoffs.
In the Klein--Gordon case they may depend on the bounded time support.
\begin{enumerate}
\item If $q>2$ and $\frac{2\alpha}q<\sigma$, then
\begin{equation}\label{infty-bound}
        \|\mathcal T_q^\kappa(W_1,W_2)\|_{\mathfrak S^{\infty}}\le C(1+|\kappa|)^M\|W_1\|_{L_t^{\widetilde{q},\infty}L_x^2}\|W_2\|_{L_t^{\widetilde{q},\infty}L_x^2}.
\end{equation}
\item If $q>4$ and $\frac{2\alpha}q<\sigma$, then
\begin{equation}\label{2-bound}
        \|\mathcal T_q^\kappa(W_1,W_2)\|_{\mathfrak S^{2}}\le C(1+|\kappa|)^M\|W_1\|_{L_t^{\widetilde{q},4}L_x^2}\|W_2\|_{L_t^{\widetilde{q},4}L_x^2}.
\end{equation}
\item If $q>4$ and $\frac{\widetilde q}2<p\le2$, then
  \begin{equation}\label{intermediate-bound}
        \|\mathcal T_q^\kappa(W_1,W_2)\|_{\mathfrak S^{p}}\le C(1+|\kappa|)^M\|W_1\|_{L_t^{\widetilde{q},2}L_x^2}\|W_2\|_{L_t^{\widetilde{q},2}L_x^2}.
\end{equation}
\end{enumerate}
\end{lem}

We first show how Lemma~\ref{lem-three-bounds} implies \eqref{goal-bound}. Note that the condition $\frac{2\alpha}{q}<\sigma$ holds for all $q>2$ if $\sigma\ge\alpha$, and for all $q>4$ in any case, since $\sigma\ge\frac12\ge\frac\alpha2$. We use Stein's interpolation for the analytic family $z\mapsto\mathcal T_{q(z)}^{\kappa(z)}(W_{1,z},W_{2,z})$, where $\frac1{q(z)}=\frac{1-z}{q_0}+\frac{z}{q_1}$, so that $-2s(q(z))$ is affine in $z$ and its imaginary part on $\mathrm{Re}\,z\in\{0,1\}$ is linear in $\mathrm{Im}\,z$, and complex interpolation of Lorentz spaces \cite{BL1976}.

\noindent\textbf{Case 1: $2<q\le4$ (so that $\sigma\ge\alpha$).}
Let $0<\eta<\frac12-\frac1q$ and set $\frac{1}{q_0}=\frac{1}{4}-\eta$, $\frac{1}{q_1}=\frac{1}{2}-\eta$, so that $q_1<q<q_0$, $q_0>4$ and $q_1>2$. Let $\theta\in(0,1)$ be such that $\frac1q=\frac{1-\theta}{q_0}+\frac\theta{q_1}$. Then $\frac1{\widetilde q_0}=\frac14+\eta$, $\frac1{\widetilde q_1}=\eta$ and $\frac1{\widetilde q}=\frac{1-\theta}4+\eta$.
We interpolate \eqref{2-bound} at $q_0$ and \eqref{infty-bound} at $q_1$, where in the latter we replace $L^{\widetilde q_1,\infty}$ by the smaller space $L^{\widetilde q_1,1/\eta}$. This gives the bound in $\mathfrak S^p$ with $\frac1p=\frac{1-\theta}2$ for $W_1,W_2\in L_t^{\widetilde q,u}(\mathrm d\mu;L_x^2)$, where $\frac1u=\frac{1-\theta}4+\theta\eta\le\frac1{\widetilde q}$. Since $L^{\widetilde q}\subseteq L^{\widetilde q,u}$, we obtain \eqref{goal-bound} with
\[
\frac2{\widetilde q}-\frac1p=2\eta,
\]
which can be made arbitrarily small.

\noindent\textbf{Case 2: $q>4$.} Then $\frac14<\frac1{\widetilde q}<\frac12$. Let $X=\frac2{\widetilde q}-\frac12\in(0,\frac12)$, and choose $0<\eta<X/2$, $0<\omega<4\eta(1-2X)$, and $\theta=2X-\omega\in(0,1)$. Set $\frac1{\widetilde q_0}=\frac14+\eta$ and define $\widetilde q_1$ by $\frac1{\widetilde q}=\frac{1-\theta}{\widetilde q_0}+\frac\theta{\widetilde q_1}$. A direct computation gives
\[
\frac1{\widetilde q_1}=\frac{X-\omega/4-\eta(1-\theta)}{2X-\omega},
\]
and the choice of $\eta,\omega$ ensures $\frac14<\frac1{\widetilde q_1}<\frac12$, i.e. $4<q_1<\infty$. We also have $q_0>4$. We interpolate \eqref{2-bound} at $q_0$ and \eqref{intermediate-bound} at $q_1$ with $p_1\in(\frac{\widetilde q_1}2,2]$. This gives the bound in $\mathfrak S^p$, $\frac1p=\frac{1-\theta}2+\frac\theta{p_1}$, for $W_1,W_2\in L_t^{\widetilde q,u}(\mathrm d\mu;L_x^2)$ with
\[
\frac1u=\frac{1-\theta}4+\frac\theta2=\frac{1+\theta}4=\frac1{\widetilde q}-\frac\omega4<\frac1{\widetilde q}.
\]
Thus \eqref{goal-bound} holds for this $p$. Since $\frac2{\widetilde q}=\frac{1-\theta}2+2\eta(1-\theta)+\frac{2\theta}{\widetilde q_1}$, letting $p_1\to\widetilde q_1/2$ we obtain all $p$ with $\frac1p<\frac2{\widetilde q}-2\eta(1-\theta)$. Since $\eta$ can be taken arbitrarily small, this proves \eqref{goal-bound} for all $p>\widetilde q/2$.

It remains to prove Lemma~\ref{lem-three-bounds}.

\begin{proof}[Proof of \eqref{infty-bound}]
Let $\gamma=d-2s(q)=\frac{2a\alpha}q$. By the assumption, $0<\gamma<a\sigma$, so Lemma~\ref{lem-full-decay} gives that the kernel $K$ of $U_\phi(t-\tau)\Psi(D)^{-2s(q)+i\kappa}$ (with frequency truncation) satisfies $|K(t-\tau,x)|\lesssim(1+|\kappa|)^M|t-\tau|^{-2\alpha/q}$ for $t,\tau\in E$. For $F\in L_t^2(\rd\mu;L_x^2)$, set $H(\tau)=\|W_2F(\tau)\|_1$ and
$J_\lambda H(t)=\int|t-\tau|^{-\lambda}H(\tau)\,\rd\mu(\tau)$.
Spatial kernel domination and Lorentz H\"older give
\begin{align*}
 \|\mathcal T_q^\kappa(W_1,W_2)F\|_{L_t^2L_x^2}
 &\lesssim(1+|\kappa|)^M
 \|W_1\|_{L_t^{\widetilde q,\infty}L_x^2}
 \|J_{2\alpha/q}H\|_{L_t^{q,2}}\\
 &\lesssim(1+|\kappa|)^M
 \|W_1\|_{L_t^{\widetilde q,\infty}L_x^2}\|H\|_{L_t^{q',2}}\\
 &\lesssim(1+|\kappa|)^M
 \prod_{i=1}^2\|W_i\|_{L_t^{\widetilde q,\infty}L_x^2}
 \|F\|_{L_t^2L_x^2}.
\end{align*}
Here we used $\frac12=\frac1{\widetilde q}+\frac1q$, $\frac1{q'}=\frac1{\widetilde q}+\frac12$, and the fact that the fractional integral of order $\alpha-\lambda$ with respect to $\mu$ maps $L^{q',2}(\mathrm d\mu)$ to $L^{q,2}(\mathrm d\mu)$ when $\lambda=\frac{2\alpha}q$, which follows from Lemma~\ref{lem_HLS} by duality.
\end{proof}

\begin{proof}[Proof of \eqref{2-bound}]
As above, the kernel satisfies $|K(t-\tau,x)|\lesssim(1+|\kappa|)^M|t-\tau|^{-2\alpha/q}$. Since $\lambda=\frac{4\alpha}q<\alpha$, Lemma~\ref{lem_HLS} with $p=p_1=\frac{\widetilde q}2$ and $r_1=r_2=2$ gives
\begin{align*}
 &\|\mathcal T_q^\kappa(W_1,W_2)\|_{\mathfrak S^2}^2\\
 &\quad=\iiiint |W_1(t,x)|^2|K(t-\tau,x-y)|^2|W_2(\tau,y)|^2
       \,\rd x\rd y\rd\mu(t)\rd\mu(\tau)\\
 &\quad\lesssim(1+|\kappa|)^{2M}
       \iint\frac{\|W_1(t)\|_2^2\|W_2(\tau)\|_2^2}
                    {|t-\tau|^{4\alpha/q}}\,\rd\mu(t)\rd\mu(\tau)\\
 &\quad\lesssim(1+|\kappa|)^{2M}
       \prod_{i=1}^2\|W_i\|_{L_t^{\widetilde q,4}L_x^2}^2.
\end{align*}
Here $2=(4\alpha/q)/\alpha+4/\widetilde q$, and squaring changes
the time Lorentz indices $(\widetilde q,4)$ to $(\widetilde q/2,2)$:
\[
 \|\|W(t)\|_2^2\|_{L_t^{\widetilde q/2,2}}
 =\|W\|_{L_t^{\widetilde q,4}L_x^2}^2.
\]
\end{proof}

For the proof of \eqref{intermediate-bound}, we need two lemmas.
\begin{lem}\label{lem-young}
Let $\lambda>0$, $\sigma\ge\frac12$, and $\kappa_\lambda(t)=(1+\lambda|t|)^{-2\sigma}$. If $p_1,p_2\in[1,\infty]$ and $1<\frac1{p_1}+\frac1{p_2}\le2$, then
\[
\begin{aligned}
 &\iint|G_1(t)G_2(\tau)|\kappa_\lambda(t-\tau)
      \,\rd\mu(t)\rd\mu(\tau)\\
 &\qquad\le C(\Fmu\lambda^{-\alpha})^{2-1/p_1-1/p_2}
      \|G_1\|_{L^{p_1}(\rd\mu)}\|G_2\|_{L^{p_2}(\rd\mu)}.
\end{aligned}
\]
\end{lem}
\begin{proof}
Let $\frac1m=2-\frac1{p_1}-\frac1{p_2}\in[0,1)$. If $m<\infty$, then $2\sigma m>1\ge\alpha$, and by \eqref{frac_meas}
\[
\begin{aligned}
 \sup_t\int\kappa_\lambda(t-\tau)^m\,\rd\mu(\tau)
 &\le\sum_{j\ge0}2^{-2\sigma m(j-1)}
       \mu\{\tau:|t-\tau|<2^j\lambda^{-1}\}\\
 &\lesssim\Fmu\lambda^{-\alpha}
       \sum_{j\ge0}2^{j(\alpha-2\sigma m)}
 \lesssim\Fmu\lambda^{-\alpha}.
\end{aligned}
\]
Thus, including the trivial case $m=\infty$,
\[
 B:=\sup_t\|\kappa_\lambda(t-\cdot)\|_{L^m(\rd\mu)}
 \lesssim(\Fmu\lambda^{-\alpha})^{1/m}.
\] By H\"older's inequality, the bilinear form is bounded by $B\|G_1\|_{L^1}\|G_2\|_{L^{m'}}$ and, by symmetry, by $B\|G_1\|_{L^{m'}}\|G_2\|_{L^1}$. Bilinear interpolation between these two bounds gives the claim on the line $\frac1{p_1}+\frac1{p_2}=1+\frac1{m'}$.
\end{proof}

\begin{lem}[{\cite[Section~3.13, Exercise~5(b)]{BL1976}}]\label{lem-bilinear}
Let $(A_0,A_1)$, $(B_0,B_1)$, $(C_0,C_1)$ be compatible couples of Banach spaces and let $\mathcal T$ be a bilinear operator which is bounded from $A_0\times B_0$ to $C_0$, from $A_0\times B_1$ to $C_1$, and from $A_1\times B_0$ to $C_1$. Then $\mathcal T$ is bounded from $(A_0,A_1)_{\theta_0,\rho_1}\times(B_0,B_1)_{\theta_1,\rho_2}$ to $(C_0,C_1)_{\theta,1}$, provided $\theta_0,\theta_1\in(0,1)$, $\theta=\theta_0+\theta_1<1$, and $\frac1{\rho_1}+\frac1{\rho_2}\ge1$.
\end{lem}
For a Banach space $X$, let $\ell^\nu_\rho(X)$ denote the space of $X$-valued sequences $\{x_k\}_{k\in\mathbb Z}$ with $\|\{2^{k\nu}\|x_k\|_X\}\|_{\ell^\rho}<\infty$. We also use the identity $(\ell^{\nu_0}_\infty(X),\ell^{\nu_1}_\infty(X))_{\theta,1}=\ell^{\nu}_1(X)$, $\nu=(1-\theta)\nu_0+\theta\nu_1$, valid for $\nu_0\neq\nu_1$ \cite[Theorem~5.6.1]{BL1976}.

\begin{proof}[Proof of \eqref{intermediate-bound}]
Let $s=s(q)$. Fix $k$, and let $\widetilde\psi_k=\widetilde\psi(2^{-k}\cdot)$ with $\widetilde\psi\in C_c^\infty(0,\infty)$ satisfying $\widetilde\psi\psi=\psi$.

\emph{Trace class bound.} We factor
$\mathcal T_{q,k}^\kappa(W_1,W_2)=A_kB_k^*$, where
\[
 \begin{aligned}
 A_k&=W_1U_\phi\Psi(D)^{-s+i\kappa}\widetilde\psi_k(|D|),\\
 B_k&=\overline{W_2}U_\phi\Psi(D)^{-s}\psi_k(|D|).
 \end{aligned}
\]
These operators map $L^2(\mathbb R^d)$ to $L^2(\rd\mu\,\rd x)$.
By Plancherel's theorem,
\[
\begin{aligned}
 \|A_k\|_{\mathfrak S^2}^2
 &=c_d\iint|W_1(t,x)|^2\,\rd x\rd\mu(t)
       \int\Psi(\xi)^{-2s}\widetilde\psi_k(|\xi|)^2\,\rd\xi\\
 &\lesssim2^{k(d-2s)}\|W_1\|_{L_t^2L_x^2}^2,
\end{aligned}
\]
and similarly for $B_k$. Hence
\begin{equation}\label{local-1-bound}
 \begin{split}
 \|\mathcal T_{q,k}^\kappa(W_1,W_2)\|_{\mathfrak S^1}
 &\le\|A_k\|_{\mathfrak S^2}\|B_k\|_{\mathfrak S^2}\\
 &\lesssim2^{k(d-2s)}\prod_{i=1}^2\|W_i\|_{L_t^2L_x^2}.
 \end{split}
\end{equation}

\emph{Hilbert--Schmidt bound.} Let $\widetilde q_1,\widetilde q_2\in[2,\infty]$ with $\frac1{\widetilde q_1}+\frac1{\widetilde q_2}>\frac12$. Set $G_i(t)=\|W_i(t)\|_2^2$, $p_i=\widetilde q_i/2$, and $C_k=(1+|\kappa|)^{2M}2^{2k(d-2s)}$. Applying Lemma~\ref{lem-dyadic-disp} and Lemma~\ref{lem-young} with $\lambda=2^{ak}$ gives
\begin{equation}\label{local-2-bound}
 \begin{split}
 \|\mathcal T_{q,k}^\kappa(W_1,W_2)\|_{\mathfrak S^2}^2
 &\lesssim C_k\iint G_1(t)G_2(\tau)
 (1+2^{ak}|t-\tau|)^{-2\sigma}\,\rd\mu(t)\rd\mu(\tau)\\
 &\lesssim C_k2^{-2ak\alpha(1-1/\widetilde q_1-1/\widetilde q_2)}
 \prod_{i=1}^2\|W_i\|_{L_t^{\widetilde q_i}L_x^2}^2.
 \end{split}
\end{equation}

\emph{Interpolation.} Let $1<p\le2$ and $\theta_p=2-\frac2p\in(0,1]$, and let
\[
\Delta(p)=\Big\{\Big(\frac1{w_1},\frac1{w_2}\Big):\ \frac1p-\frac12\le\frac1{w_i}\le\frac12,\ \frac1{w_1}+\frac1{w_2}>\frac1p\Big\}.
\]
For $(\frac1{w_1},\frac1{w_2})\in\Delta(p)$, define $\widetilde q_i$ by $\frac1{w_i}=\frac{1-\theta_p}2+\frac{\theta_p}{\widetilde q_i}$. Then $\frac1{\widetilde q_i}\in[0,\frac12]$ and $\frac1{\widetilde q_1}+\frac1{\widetilde q_2}>\frac12$. Bilinear complex interpolation between \eqref{local-1-bound} and \eqref{local-2-bound} gives
\begin{equation}\label{goal-bound2}
    \|\mathcal T_{q,k}^\kappa(W_1,W_2)\|_{\mathfrak S^{p}}\lesssim (1+|\kappa|)^M2^{-k\nu(w_1,w_2)}\|W_1\|_{L_t^{w_1}L_x^2}\|W_2\|_{L_t^{w_2}L_x^2},
\end{equation}
where, since $\theta_p(1-\frac1{\widetilde q_1}-\frac1{\widetilde q_2})=1-\frac1{w_1}-\frac1{w_2}$,
\[
\nu(w_1,w_2)=a\alpha\Big(1-\frac1{w_1}-\frac1{w_2}\Big)-(d-2s)=a\alpha\Big(1-\frac1{w_1}-\frac1{w_2}-\frac2q\Big).
\]
In particular $\nu(\widetilde q,\widetilde q)=0$. Now let $q>4$ and $\frac{\widetilde q}2<p\le2$. Then $(\frac1{\widetilde q},\frac1{\widetilde q})$ lies in the interior of $\Delta(p)$, because $\frac1p-\frac12<\frac2{\widetilde q}-\frac12<\frac1{\widetilde q}<\frac12$ and $\frac2{\widetilde q}>\frac1p$. Choose $\epsilon>0$ so small that $(\frac1{w_i},\frac1{w_j})\in\Delta(p)$ for $i,j\in\{0,1\}$, where $\frac1{w_0}=\frac1{\widetilde q}+\epsilon$ and $\frac1{w_1}=\frac1{\widetilde q}-2\epsilon$. Then $\nu(w_0,w_0)=-2a\alpha\epsilon$ and $\nu(w_0,w_1)=\nu(w_1,w_0)=a\alpha\epsilon$. Thus the bilinear operator $\mathcal T=\{\mathcal T_{q,k}^\kappa\}_k$ is bounded from
\[
\begin{aligned}
 L_t^{w_0}L_x^2\times L_t^{w_0}L_x^2
 &\longrightarrow\ell_\infty^{-2a\alpha\epsilon}(\mathfrak S^p),\\
 L_t^{w_0}L_x^2\times L_t^{w_1}L_x^2
 &\longrightarrow\ell_\infty^{a\alpha\epsilon}(\mathfrak S^p),\\
 L_t^{w_1}L_x^2\times L_t^{w_0}L_x^2
 &\longrightarrow\ell_\infty^{a\alpha\epsilon}(\mathfrak S^p).
\end{aligned}
\]
Here we use the identity stated after \eqref{real-ell}, with inner space $L_x^2$:
\[
 (L_t^{w_0}(\rd\mu;L_x^2),L_t^{w_1}(\rd\mu;L_x^2))_{1/3,2}
 =L_t^{\widetilde q,2}(\rd\mu;L_x^2),
\]
because
\[
 \frac23\left(\frac1{\widetilde q}+\epsilon\right)
 +\frac13\left(\frac1{\widetilde q}-2\epsilon\right)
 =\frac1{\widetilde q}.
\]
Apply Lemma~\ref{lem-bilinear} with $\theta_0=\theta_1=1/3$ and
$\rho_1=\rho_2=2$. The output weight is zero, since
\[
 \frac13(-2a\alpha\epsilon)+\frac23(a\alpha\epsilon)=0.
\]
Consequently $\mathcal T$ maps the product of the two Lorentz spaces
above into $\ell_1^0(\mathfrak S^p)$.
That is,
\[
\sum_k\|\mathcal T_{q,k}^\kappa(W_1,W_2)\|_{\mathfrak S^p}\lesssim(1+|\kappa|)^M\|W_1\|_{L_t^{\widetilde q,2}(\mathrm d\mu;L_x^2)}\|W_2\|_{L_t^{\widetilde q,2}(\mathrm d\mu;L_x^2)}.
\]
In the Klein--Gordon case, the sum is over $k\ge1$, and the low frequency part satisfies, by the same argument as for \eqref{local-1-bound},
$\|\mathcal T_{q,\le0}^\kappa(W_1,W_2)\|_{\mathfrak S^p}\le\|\mathcal T_{q,\le0}^\kappa(W_1,W_2)\|_{\mathfrak S^1}\lesssim\|W_1\|_{L_t^2(L_x^2)}\|W_2\|_{L_t^2(L_x^2)}$, which is bounded by the right-hand side since $\mu$ is finite and $\widetilde q>2$. The triangle inequality gives \eqref{intermediate-bound}.
\end{proof}

\subsection{Low-dimensional non-sharp estimates}
\begin{prop}[Low-dimensional non-sharp estimates]\label{prop-lowdim-band}
Let $1/2<\alpha<1$, and let $\mu$ satisfy \eqref{frac_meas}.
Consider either the fractional Schr\"odinger equation in dimension
$d=1$, with $a>0$, $a\ne1$, or the wave or Klein--Gordon equation
in dimension $d=2$, with $a=1$. Use the corresponding $\phi$ and
$\Psi$ from \eqref{three-cases}, and put
\[
 s(q,r)=\frac d2-\frac dr-\frac{a\alpha}{q}.
\]
For every orthonormal family $(g_j)_j$ in $L^2(\mathbb R^d)$,
the following assertions hold.
\begin{enumerate}
\item If $4\alpha<q\le4$ and $r=\infty$, then
\begin{equation}\label{eq:lowdim-OD-band}
 \left\|\sum_j n_j
 |U_\phi(t)\Psi(D)^{-s(q,\infty)}g_j|^2\right\|_
 {L_t^{q/2}(d\mu;L_x^\infty)}
 \le C\Fmu^{2/q}\|n\|_{\ell^\beta}
\end{equation}
for every
\[
 1\le\beta<
 B_\alpha(q):=\frac{2(1-\alpha)}{1-2\alpha+4\alpha/q}.
\]
\item If $2<r<\infty$ and
\[
 4\alpha<q(1-2/r)\le4,
\]
then
\begin{equation}\label{eq:lowdim-finite-band}
 \left\|\sum_j n_j
 |U_\phi(t)\Psi(D)^{-s(q,r)}g_j|^2\right\|_
 {L_t^{q/2}(d\mu;L_x^{r/2})}
 \le C\Fmu^{2/q}\|n\|_{\ell^\beta}
\end{equation}
for every
\[
 1\le\beta<
 \frac{2(1-\alpha)}{1-2\alpha+2/r+4\alpha/q}.
\]
\end{enumerate}
The constant may depend on the fixed exponents and, in the
Klein--Gordon case, on the diameter of the bounded time support,
but is independent of the number of functions and all frequency
cutoffs. Homogeneous expressions are defined first on annular
data and then by their canonical $L^2$ extensions. In the
Klein--Gordon case the same assertions hold with any regularity
larger than $s(q,r)$, by the Schatten ideal property.
\end{prop}

\begin{proof}
Here $\sigma=1/2$. For $q>4\alpha$,
\[
 \frac{2\alpha}q<\frac12=\sigma,\qquad
 s(q,\infty)=\frac d2-\frac{a\alpha}q.
\]
Lemma~\ref{lem-scalar-OD-completion} therefore gives
\begin{equation}\label{eq:scalar-lowdim-OD}
 \|U_\phi(t)\Psi(D)^{-s(q,\infty)}g\|_{L_t^q(\rd\mu;L_x^\infty)}
 \le C\Fmu^{1/q}\|g\|_2.
\end{equation}
The input is the annular $L^2$ completion from that lemma, so it also
applies when the additional hypothesis in Theorem~\ref{Class-Schro}
fails at an auxiliary scalar pair.

Fix first $4\alpha<q\le4$. Choose
\[
 4\alpha<q_0<q,\qquad q_1>4,
 \qquad \beta_1=2<q_1/2,
\]
and let $\theta\in(0,1)$ satisfy
$1/q=(1-\theta)/q_0+\theta/q_1$.
The scalar bound \eqref{eq:scalar-lowdim-OD} and the triangle
inequality give the orthonormal estimate with $\beta_0=1$ at
$(q_0,\infty)$. There is a direct Hilbert--Schmidt input at $(q_1,\infty)$:
since $q_1>4$, $4\alpha/q_1<\alpha$, and the squared full-kernel
bound above followed by Lemma~\ref{lem_HLS} gives
\[
 \|B_{s(q_1,\infty)}(V)\|_{\Schatten^2}
 \le C\Fmu^{2/q_1}\|V\|_{L_t^{(q_1/2)'}(\rd\mu;L_x^1)}.
\]
For a bounded supported potential and a nonnegative frequency cutoff
$m_N$, let $B_{s,N}(V)$ denote the corresponding initial-space operator,
where $s=s(q_1,\infty)$. Expanding its Fourier-kernel norm gives, up to
fixed Fourier-normalization constants,
\[
 \begin{aligned}
 \|B_{s,N}(V)\|_{\Schatten^2}^2
 &\lesssim\iiiint |V(t,x)V(\tau,y)|
 |K_{s,N}(t-\tau,x-y)|^2\,\rd x\rd y\rd\mu(t)\rd\mu(\tau)\\
 &\lesssim\iint |t-\tau|^{-4\alpha/q_1}
 \|V(t)\|_1\|V(\tau)\|_1\,\rd\mu(t)\rd\mu(\tau).
 \end{aligned}
\]
The full kernel bound in Lemma~\ref{lem-full-decay} is uniform in $N$.
For $p_1=(q_1/2)'$,
\[
 \frac{4\alpha/q_1}{\alpha}+\frac2{p_1}
 =\frac4{q_1}+2\left(1-\frac2{q_1}\right)=2.
\]
Lemma~\ref{lem_HLS} bounds the last expression by
$C\Fmu^{4/q_1}\|V\|_{L_t^{p_1}L_x^1}^2$.
The same annular completion and Schatten weak limit identify the
uncut operator, so the displayed HS bound follows after taking square
roots.
Trace duality gives coefficient exponent $\beta_1=2$.

We interpolate the initial-space operators \eqref{eq-initial-potential}.
Write $P_i=q_i/2$ and $P=q/2$.
The two input estimates give
\[
 \|B_{s(q_i,\infty)}(V)\|_{\mathfrak S^{\beta_i'}}
 \le C_i\Fmu^{2/q_i}
 \|V\|_{L_t^{P_i'}(d\mu;L_x^1)},
 \qquad i=0,1,
\]
with $\beta_0'=\infty$. These bounds hold for complex $V$ as
well: factor $V=|V|^{1/2}J|V|^{1/2}$ with $J$ a contraction,
then use the Schatten ideal inequality.

For a bounded simple rectangular potential $V$ with compact support, put
$h(t)=\|V(t,\cdot)\|_1$ and $v=V/h$ on $\{h>0\}$, zero
otherwise. Define
\[
 \begin{gathered}
 s(z)=(1-z)s(q_0,\infty)+zs(q_1,\infty),\\
 \frac1{p(z)}=\frac{1-z}{P_0'}+\frac z{P_1'},\qquad
 V_z=h^{P'/p(z)}v.
 \end{gathered}
\]
Regularize the initial-space analytic family by
\[
 \mathcal B_N(z)=m_N(D)^{1/2}\Psi(D)^{-s(z)}
 \left(\int U_\phi^*V_zU_\phi\,d\mu\right)
 \Psi(D)^{-s(z)}m_N(D)^{1/2}.
\]
On the boundary $\Re z=j$, $j\in\{0,1\}$, the imaginary factors in $s(z)$
are unitary Fourier multipliers on the initial $L^2$ space and
preserve every Schatten norm. Thus no $L_x^\infty$ or $L_x^1$
bound for imaginary powers is needed. Moreover
\[
 \|V_{j+iy}\|_{L_t^{P_j'}L_x^1}
 =\|h\|_{L^{P'}}^{P'/P_j'},
 \qquad
 P'\left(\frac{1-\theta}{P_0'}+\frac\theta{P_1'}\right)=1.
\]
Stein interpolation, followed by annular cutoff removal on dense
initial vectors, gives
\[
 \|B_{s(q,\infty)}(V)\|_{\mathfrak S^{\beta'}}
 \le C\Fmu^{2/q}\|V\|_{L_t^{P'}L_x^1},
 \qquad
 \frac1\beta=1-\theta+\frac\theta{\beta_1}.
\]
Since both endpoint bounds are strong Schatten bounds with strong
potential spaces, density duality gives \eqref{eq:lowdim-OD-band} with
strong time and sequence norms; in particular, no weak mixed-norm
interpolation is used.

Letting $q_0\downarrow4\alpha$ and $q_1\downarrow4$, with
$\beta_1=2$, in the exponent calculation gives the scalar
weight
\[
 b=1-\theta=\frac{4\alpha/q-\alpha}{1-\alpha},
 \qquad
 \frac1{\beta_*}=b+\frac{1-b}{2}
 =\frac{1-2\alpha+4\alpha/q}{2(1-\alpha)}.
\]
Every $\beta<\beta_*=B_\alpha(q)$ is obtained by choosing
input exponents sufficiently close to these
limits, followed if necessary by the inclusion
$\ell^\beta\subset\ell^{\beta_{\mathrm{obtained}}}$.

For part (2), put $\vartheta=1-2/r\in(0,1)$ and $q_b=\vartheta q$.
Interpolate the just proved strong estimate at $(q_b,\infty)$
with the energy estimate at $(\infty,2)$, where $s=0$ and
$\beta=1$. The interpolation has weight $\vartheta$ for the OD input.
The spatial potential endpoints here are $L_x^\infty$ at energy
and $L_x^1$ at OD. At the target put $p=(q/2)'$ and
$\rho=(r/2)'=1/\vartheta$, normalize $\|V\|_{L_t^pL_x^\rho}=1$,
and write $h(\tau)=\|V(\tau)\|_\rho$, $v=V/h$.
We use the mixed potential family
\[
 V_z(\tau,x)=\operatorname{phase}(V)
       |v(\tau,x)|^{\rho z}
       h(\tau)^{p[(1-z)+z/P_b']},\qquad P_b=q_b/2,
\]
set to zero on $V=0$. At $z=\vartheta$ it equals $V$, since
$1/p=1-\vartheta+\vartheta/P_b'$. On the two boundaries its spatial norms are
$h^p$ and $h^{p/P_b'}$, giving respectively unit norms in
$L_t^1L_x^\infty$ and $L_t^{P_b'}L_x^1$.
Use this family in the same initial-space regularized sandwich;
imaginary regularity powers act only on the initial $L^2$ space. Both endpoint potential
spaces and Schatten spaces are strong, so the output remains
$L_t^{q/2}L_x^{r/2}$ and the coefficient space remains strong
$\ell^\beta$. The affine regularity and constants are
\[
 \vartheta\,s(q_b,\infty)=\frac d2-\frac dr-\frac{a\alpha}{q},
 \qquad \vartheta\frac2{q_b}=\frac2q.
\]
Finally,
\[
 \frac1{\beta_*}
 =1-\vartheta+\frac{\vartheta}{B_\alpha(q_b)}
 =\frac{1-2\alpha+2/r+4\alpha/q}{2(1-\alpha)}.
\]
This proves the stated strict range.

For Klein--Gordon and $s\ge s_0=s(q,r)$, let
$C=\langle D\rangle^{s_0-s}$, so $\|C\|\le1$. Then
$B_s(V)=CB_{s_0}(V)C$ and the Schatten ideal property gives the
same estimate.
\end{proof}

\begin{rem}
The scalar endpoint $q=4\alpha$ and the coefficient endpoint
$\beta=B_\alpha(q)$ are not obtained. For $4\alpha<q<4$,
$B_\alpha(q)<q/2$; equality holds at $q=4$, where the result
still has the strict restriction $\beta<2$.
\end{rem}

\section{Critical Schr\"odinger estimates}\label{sec-critical}
Throughout this section $d\ge2$ and $0<\alpha\le1$.
Write $M=\Fmu$ for the Frostman constant, distinct from the total mass
$\mu(\R)$. We use $U(t)=e^{it\Delta}$; all time-kernel estimates
depend only on $|t-\tau|$ and apply to the opposite sign as well. Note that $M>0$ since $\mu\ne0$. We prove
Theorem~\ref{thm-quadratic-critical}. For part (1) put
\[
 p=1+d/\alpha>2,\qquad P_A=p',\qquad R_A=(p/2)'.
\]
The proof uses only the two inputs
\begin{equation}\label{critical-energy}
 \|U(t)\|_{2\to2}=1,\qquad |K_{U(t)U(\tau)^*}(x,y)|
 \le D_d|t-\tau|^{-d/2}\quad(t\ne\tau).
\end{equation}
By Lorentz sequence duality it suffices to prove
\begin{equation}\label{critical-weakS}
 \left\|\int U(t)^*V(t)U(t)\,\rd\mu(t)\right\|_{\Schatten^{p,\infty}}
 \le C_{d,\alpha}M^{1/P_A}\|V\|_{L_t^p(\rd\mu;L_x^{p/2})}.
\end{equation}
Here $\Schatten^{p,u}$ denotes the compact operators whose singular
values belong to $\ell^{p,u}$. Initially $V$ is bounded, positive and
simple, with bounded time support and spatial support of finite measure.
The proof follows the spectral construction of \cite{HT2026}, with
Schur estimates replacing time convolution.

\subsection{The Lorentz endpoint at \texorpdfstring{$A_{d/2}$}{A}}

\paragraph{The time-tail estimate.}
For $\ell>0$, dyadic shells and \eqref{frac_meas} give
\begin{align}\label{ht:tail}
 \sup_t\int_{|t-\tau|\ge\ell}|t-\tau|^{-d}\,\rd\mu(\tau)
 &\le \sum_{m\ge0}(2^m\ell)^{-d}
            \mu((t-2^{m+1}\ell,t+2^{m+1}\ell))\notag\\
 &\le \frac{2^\alpha}{1-2^{\alpha-d}}M\ell^{\alpha-d}.
\end{align}
The same bound holds with the two variables interchanged.
Consequently, Schur's test on $L^2(\rd\mu)$ gives
\begin{equation}\label{ht:schur}
 \begin{split}
 &\left|\iint \ind_{\{|j_\ell(t)-j_\ell(\tau)|\ge3\}}
 |t-\tau|^{-d}g_1(t)g_2(\tau)\,\rd\mu(t)\rd\mu(\tau)\right|\\
 &\qquad\le C_{d,\alpha}M\ell^{\alpha-d}\|g_1\|_2\|g_2\|_2,
 \end{split}
\end{equation}
where $I_{\ell,j}=[j\ell,(j+1)\ell)$ and $j_\ell(t)$ is its interval
index. Separation of the indices implies $|t-\tau|\ge2\ell$.
No doubling property of $\mu$ is used; the geometric series in
\eqref{ht:tail} converges because $d>\alpha$.

\paragraph{Compactness.}
Set $f=\sqrt V$, $T=fU:L_x^2\to\mathcal K$, and
$\mathcal K=L^2(\rd\mu\,\rd x)$. Energy and finiteness of the time
support give
\[
 \|Th\|_{\mathcal K}^2\le \|V\|_\infty\mu(\text{time support})\|h\|_2^2.
\]
Thus $A=TT^*\ge0$ is well defined and bounded, without invoking any
endpoint Strichartz bound. Let $R_j=\ind_{I_{\ell,j}}$ on $\mathcal K$.
Since $\mu(I_{\ell,j})\le M\ell^\alpha$,
\[
 \|R_jT\|^2\le M\ell^\alpha\|V\|_\infty,\qquad
 \|R_iAR_j\|\le M\ell^\alpha\|V\|_\infty.
\]
For each fixed offset $h\in\{-2,-1,0,1,2\}$, orthogonality of the
input and output interval projections gives
\[
 \left\|\sum_jR_{j+h}AR_j\right\|
 \le M\ell^\alpha\|V\|_\infty.
\]
Hence the fixed-length near part has operator norm at most
$5M\ell^\alpha\|V\|_\infty$, which tends to zero as $\ell\downarrow0$.
The separated part has Hilbert--Schmidt norm squared bounded by
\[
\begin{aligned}
 &D_d^2\iint_{|j_\ell(t)-j_\ell(\tau)|\ge3}
 |t-\tau|^{-d}\|V(t)\|_1\|V(\tau)\|_1\,\rd\mu(t)\rd\mu(\tau)\\
 &\qquad\le C_{d,\alpha}M\ell^{\alpha-d}
      \|V\|_{L_t^2(\rd\mu;L_x^1)}^2<\infty.
\end{aligned}
\]
The kernel identity for separated pieces follows from the explicit
Schr\"odinger kernel and Fubini; spatial support of finite measure and
positive time separation make these pieces Hilbert--Schmidt.
The near part is defined as an operator block sum, not by an absolutely
integrable kernel on the time diagonal. Also, \eqref{frac_meas} implies that
$\mu$ is atomless. These observations justify the decomposition.
It follows that $A$ is an operator-norm limit of Hilbert--Schmidt
operators and is compact. Its spectral bands therefore have finite rank
away from zero. Note that \eqref{critical-weakS} was not used here.

\paragraph{The near--far decomposition.}
First assume $V$ is dyadic simple:
\[
 V=\sum_k2^k\ind_{E_k},
\]
with finitely many levels and pairwise disjoint measurable $E_k$.
For $\varepsilon\in(0,1)$ and $\lambda>0$ choose
\begin{equation}\label{ht:ell}
 \ell(k,l)^\alpha
 =\frac{\varepsilon\lambda}{M}
     2^{-(k+l)/2}2^{|k-l|/4}.
\end{equation}
Thus the time-cell length is
\[
 \ell(k,l)=\left(\frac{\varepsilon\lambda}{M}
  2^{-(k+l)/2}2^{|k-l|/4}\right)^{1/\alpha}.
\]
This choice matches the local energy factor $M2^k\ell(k,l)^\alpha$
with the spectral scale $\lambda$.
Write
\[
 Q_{k,l,j}=\ind_{E_k}\ind_{I_{\ell(k,l),j}}.
\]
For each $(k,l)$, partition both time variables at this same length and
define $F$ using $|i-j|\ge3$, and $N$ using $|i-j|\le2$.
Then $A=F+N$. For fixed $V,\varepsilon,\lambda$ these are finite sums:
there are finitely many levels and finitely many active time cells.

Let $e_k(t)=|E_k(t)|$, where $|\cdot|$ denotes spatial Lebesgue measure.
Disjoint level supports in both input and output variables imply
orthogonality of distinct $(k,l)$ in the Hilbert--Schmidt calculation.
Using \eqref{ht:schur},
\begin{align*}
 \|F\|_{\Schatten^2}^2
 &\le C_d\sum_{k,l}2^{k+l}
 \iint_{|j_{\ell(k,l)}(t)-j_{\ell(k,l)}(\tau)|\ge3}
 |t-\tau|^{-d}e_k(t)e_l(\tau)\,\rd\mu(t)\rd\mu(\tau)\\
 &\le C_{d,\alpha}M\sum_{k,l}2^{k+l}
       \ell(k,l)^{\alpha-d}\|e_k\|_{L^2(\rd\mu)}\|e_l\|_{L^2(\rd\mu)}.
\end{align*}
Since $p-2=(d-\alpha)/\alpha$, \eqref{ht:ell} yields the identity
\[
 M2^{k+l}\ell(k,l)^{\alpha-d}
 =M^{p-1}(\varepsilon\lambda)^{-(p-2)}
   2^{p(k+l)/2}2^{-(p-2)|k-l|/4}.
\]
Set $b_k=\|2^{pk/2}e_k\|_{L^2(\rd\mu)}$. Because $p>2$,
the discrete sequence $2^{-(p-2)|n|/4}$ is summable; the discrete
$\ell^2$ convolution inequality gives
\[
 \sum_{k,l}2^{-(p-2)|k-l|/4}b_kb_l\le C_p\sum_k b_k^2.
\]
Furthermore,
\[
 \sum_k b_k^2
 \le\int\left(\sum_k2^{pk/2}e_k(t)\right)^2\rd\mu(t)
 =\|V\|_{L_t^p(\rd\mu;L_x^{p/2})}^p.
\]
We have proved
\begin{equation}\label{ht:far}
 \|F\|_{\Schatten^2}^2
 \le C_{d,\alpha}M^{p-1}(\varepsilon\lambda)^{-(p-2)}
       \|V\|_{L_t^p(\rd\mu;L_x^{p/2})}^p.
\end{equation}
No logarithmic loss occurs in \eqref{ht:far}.

\paragraph{Spectral compression.}
Put $P_m=\ind_{[2^m\lambda,2^{m+1}\lambda)}(A)$, $a_m=P_mAP_m$,
and $N_m=\rank P_m$. For $Q=\ind_{E_k}\ind_{I_{\ell,j}}$, energy gives
\[
 QAQ\le M2^k\ell^\alpha Q.
\]
Since $P_m\le(2^m\lambda)^{-1}A$ and $P_m\le I$,
\begin{equation}\label{ht:compression}
 QP_mQ\le\min\left\{1,\frac{M2^k\ell^\alpha}{2^m\lambda}\right\}Q.
\end{equation}
For $n\ge0$, $|h|\le2$, define the extraction map
\[
 L_{n,h}(B)=\sum_{l,j}Q_{l+n,l,j}\,B\,Q_{l,l+n,j+h}.
\]
The left projections in this sum are pairwise orthogonal, and so are
the right projections. Indeed distinct levels are disjoint, and for
a fixed level the cells form a partition. The lengths may vary with
$l$ without affecting this assertion.

For every $s\ge1$,
\begin{equation}\label{ht:pinch}
 \|L_{n,h}(B)\|_{\Schatten^s}\le\|B\|_{\Schatten^s},
 \qquad \|L_{n,h}(B)\|\le\|B\|.
\end{equation}
To verify this, index the corresponding left and right projections by
$i$ and take independent signs $\epsilon_i$. Both signed projection
sums have operator norm at most one, and
$L_{n,h}(B)=\mathbb E[(\sum_i\epsilon_iR_i)B
(\sum_i\epsilon_iS_i)]$. The ideal property and triangle inequality
prove \eqref{ht:pinch}. This proof is independent of the time measure.

For $k=l+n$, the length choice gives
\[
 \frac{M2^{l+n}\ell(l+n,l)^\alpha}{2^m\lambda}
     =\varepsilon2^{3n/4-m},\qquad
 \frac{M2^l\ell(l+n,l)^\alpha}{2^m\lambda}
     =\varepsilon2^{-n/4-m}.
\]
Denote the left and right projection square sums by
\[
 A(u)^2=\sum_{l,j}\|P_mQ_{l+n,l,j}P_0u\|^2,\qquad
 B(v)^2=\sum_{l,j}\|P_mQ_{l,l+n,j+h}P_0v\|^2.
\]
Cauchy--Schwarz and \eqref{ht:compression} give
\begin{align*}
 |\langle u,P_0L_{n,h}(a_m)P_0v\rangle|
 &\le\|a_m\|A(u)B(v)\\
 &\le2^{m+1}\lambda
 \left[\min(1,\varepsilon2^{3n/4-m})
 \min(1,\varepsilon2^{-n/4-m})\right]^{1/2}\\
 &\qquad\times\|u\|\|v\|.
\end{align*}
Consequently,
\begin{equation}\label{ht:near-block}
 \begin{split}
 \|P_0L_{n,h}(a_m)P_0\|&\le2\lambda w(n,m),\\
 w(n,m)&=2^m\left[\min(1,\varepsilon2^{3n/4-m})
                     \min(1,\varepsilon2^{-n/4-m})\right]^{1/2}.
 \end{split}
\end{equation}

\paragraph{The eigenvalue count.}
Fix $\lambda_0>0$. Compactness shows that
\[
 H_0=\sup_{\zeta\ge\lambda_0}\zeta^p
        \rank\ind_{[\zeta,2\zeta)}(A)<\infty:
\]
the ranks are bounded by the number of eigenvalues above $\lambda_0$,
and a nonempty band has $\zeta\le\|A\|$. If $H_0=0$ there is nothing
to prove above $\lambda_0$. Otherwise choose $\lambda_*\ge\lambda_0$
with
\begin{equation}\label{ht:good}
 H_0\le2\lambda_*^pN_{0,\lambda_*}.
\end{equation}
In all remaining formulas set $\lambda=\lambda_*$ and abbreviate
$N_0=N_{0,\lambda_*}>0$. Then for every $m\ge1$,
\begin{equation}\label{ht:rank-decay}
 N_m\le2^{1-pm}N_0.
\end{equation}
This uses only \eqref{ht:good}, not an already known global weak-Schatten
estimate.

For $m\le0$, \eqref{ht:near-block} gives
\[
 |\Tr(P_0L_{n,h}(a_m)P_0)|\le2\lambda_*N_0w(n,m).
\]
For $m\ge1$ we also use a second bound. By Schatten H\"older,
\eqref{ht:pinch}, and \eqref{ht:rank-decay},
\begin{align*}
 |\Tr(P_0L_{n,h}(a_m)P_0)|
 &\le N_0^{1-2/p}\|L_{n,h}(a_m)\|_{\Schatten^{p/2}}\\
 &\le N_0^{1-2/p}\,2^{m+1}\lambda_*N_m^{2/p}
 \le C_p\lambda_*N_0\,2^{-m}.
\end{align*}
Here we used $p/2>1$.
Combining the two estimates,
\begin{equation}\label{ht:near-trace}
 |\Tr(P_0L_{n,h}(a_m)P_0)|
 \le C_p\lambda_*N_0
 \begin{cases}w(n,m),&m\le0,\\
 \min(w(n,m),2^{-m}),&m\ge1.\end{cases}
\end{equation}
The factor $2^{-m}$ is needed to sum over the high spectral bands
$m\ge1$.

\paragraph{Summation over levels.}
For all $n\ge0$, $m\in\mathbb Z$,
\[
 w(n,m)\le\varepsilon^{1/2}2^{m/2-n/8}.
\]
Thus
\[
 \sum_{n\ge0,m\le0}w(n,m)\le C\varepsilon^{1/2}.
\]
For $m\ge1$, $\min(a,b)\le\sqrt{ab}$ implies
\[
 \min(w(n,m),2^{-m})
 \le\varepsilon^{1/4}2^{-m/4-n/16},
\]
and hence
\begin{equation}\label{ht:numeric}
 \sum_{n\ge0,m\le0}w(n,m)
 +\sum_{n\ge0,m\ge1}\min(w(n,m),2^{-m})
 \le C\varepsilon^{1/4}.
\end{equation}
The constant in \eqref{ht:numeric} does not depend on $\alpha$; the
dimension of $\mu$ enters only through the choice \eqref{ht:ell}.

Since $A$ is compact, its spectral-band expansion
$A=\sum_ma_m$ converges in operator norm. Extraction maps are bounded
by \eqref{ht:pinch}; compression by $P_0$ has finite rank. Thus one may
expand each extracted trace into the spectral bands, and
\eqref{ht:near-trace}--\eqref{ht:numeric} justify absolute convergence of the
resulting sum. For the near part,
\[
 N=\sum_{|h|\le2}L_{0,h}(A)
       +2\operatorname{Re}\sum_{n\ge1,|h|\le2}L_{n,h}(A).
\]
Using the five offsets and \eqref{ht:numeric},
\[
 |\Tr(P_0NP_0)|\le C_p\varepsilon^{1/4}\lambda_*N_0.
\]
Choose $\varepsilon>0$ once, depending only on $p$, so that the
right-hand constant is at most $\tfrac12\lambda_*N_0$.
Since $\Tr(P_0AP_0)\ge\lambda_*N_0$ and $A=F+N$,
\[
 \tfrac12\lambda_*N_0\le |\Tr(P_0FP_0)|
 \le\|F\|_{\Schatten^2}N_0^{1/2}.
\]
By \eqref{ht:far},
\[
 N_0\le4\lambda_*^{-2}\|F\|_{\Schatten^2}^2
 \le C_{d,\alpha}M^{p-1}\lambda_*^{-p}
      \|V\|_{L_t^p(\rd\mu;L_x^{p/2})}^p.
\]
Combining this with \eqref{ht:good}, then letting $\lambda_0\downarrow0$,
gives the same bound for every dyadic band starting at any $\zeta>0$.
Summing the disjoint bands starting at $2^j\lambda$, $j\ge0$, yields
\begin{equation}\label{ht:count}
 \rank\ind_{[\lambda,\infty)}(A)
 \le C_{d,\alpha}M^{p-1}\lambda^{-p}
      \|V\|_{L_t^p(\rd\mu;L_x^{p/2})}^p.
\end{equation}
The nonzero eigenvalues of $TT^*$ and $T^*T$ agree. Therefore
\eqref{ht:count} proves \eqref{critical-weakS} for dyadic simple positive potentials.

\paragraph{Lorentz duality.}
For an arbitrary bounded positive simple $V$ of the prescribed support,
round its finitely many positive values up to powers of two, obtaining
$V\le\widetilde V\le2V$. Positivity gives
$0\le U^*VU\le U^*\widetilde VU$ on the initial-data Hilbert space.
The min--max principle and \eqref{ht:count} give \eqref{critical-weakS} for $V$.

The passage to general potentials does not use an endpoint scalar
estimate. On compact operators the
functional
\[
 \|B\|_*:=\sup_{N\ge1}N^{1/p-1}\sum_{j=1}^N s_j(B)
\]
is an equivalent Banach norm for $\Schatten^{p,\infty}$:
$\|B\|_{p,\infty}\le\|B\|_*\le p'\|B\|_{p,\infty}$.
The triangle inequality follows from Ky Fan's inequality. Split a
complex simple $V$ into the positive and negative parts of its real
and imaginary parts. Each part has mixed norm at most $\|V\|$, so
\eqref{critical-weakS} and the triangle inequality for $\|\cdot\|_*$ give a
bounded linear map on simple potentials. Because $p,p/2<\infty$,
these potentials with bounded time support and spatial support of
finite measure are dense in $L_t^p(\rd\mu;L_x^{p/2})$. Completion
therefore gives \eqref{critical-weakS} for all potentials (up to a fixed constant).
The extension is compact and agrees with the quadratic-form integral:
the rank-one case of the next argument bounds this form, and norm
approximation identifies it.

For a finite orthonormal family and $n_j\ge0$, let
$\gamma=\sum_jn_j|h_j\rangle\langle h_j|$.
For positive simple $V$, the trace identity gives
\[
 \iint V(t,x)\sum_jn_j|U(t)h_j(x)|^2\,\rd x\rd\mu(t)
 =\Tr(\gamma U^*VU).
\]
The singular-value trace inequality and
$s_j(U^*VU)\le\|U^*VU\|_{p,\infty}j^{-1/p}$ imply
\begin{align*}
 |\Tr(\gamma U^*VU)|
 &\le\|U^*VU\|_{p,\infty}\sum_{j\ge1}n_j^*j^{-1/p}\\
 &=\|U^*VU\|_{p,\infty}\|n\|_{\ell^{p',1}}.
\end{align*}
Taking the supremum over positive simple $V$ in the unit ball of
$L_t^p(\rd\mu;L_x^{p/2})$ gives \eqref{quadratic-A} by ordinary mixed-norm
duality. Such test functions suffice by truncation and approximation,
even before knowing that the density belongs to the target space.
For complex $n_j$ use pointwise domination by the density with $|n_j|$.
Finite truncations converge for $n\in\ell^{p',1}$, since $p'>1$ and
the second Lorentz index is finite. This completes the proof of part (1) of Theorem~\ref{thm-quadratic-critical}.

\subsection{The sharp segment and finite spatial exponents}
Assume now $d>2\alpha$. Write $A=A_{d/2}$, $E=E_{d/2}$,
$P_E=1$, $R_E=d/(d-2\alpha)$. At an interior point
$(1/r,1/q)=(1-\theta)A+\theta E$, put $P=q/2$ and $R=r/2$.
Then
\begin{equation}\label{critical-segment-exponents}
 \frac1P=\frac{1-\theta}{P_A}+\theta,\qquad
 \frac1R=\frac{1-\theta}{R_A}+\frac\theta{R_E}.
\end{equation}
The regularity is zero along this segment. We shall prove the refinement
\begin{equation}\label{critical-refined}
 \left\|\sum_j n_j|U(t)f_j|^2\right\|_{L_t^P(\rd\mu;L_x^{R,P})}
 \le C M^{1/P}\|n\|_{\ell^{P,P}}\qquad(0<\theta<1).
\end{equation}
\paragraph{The two endpoint estimates.}
Fix an orthonormal family $\Phi=(f_j)$ and define on finite sequences
\[
 T_\Phi n=\sum_j n_j|e^{it\Delta}f_j|^2.
\]
This is a linear map in $n$, with constants uniform over $\Phi$.
Part (1), proved above, supplies
\begin{equation}\label{seg:Ainput}
 T_\Phi:\ell^{P_A,1}\longrightarrow
 Y_A:=L_t^{P_A}(\rd\mu;L_x^{R_A}),\quad
 \|T_\Phi\|\le C_{d,\alpha}M^{1/P_A}.
\end{equation}
The strong scalar endpoint used at $E$ is
\begin{equation}\label{seg:scalar}
 \|e^{it\Delta}f\|_{L_t^2(\rd\mu;L_x^{2R_E})}
 \le C_{d,\alpha}M^{1/2}\|f\|_2.
\end{equation}
This is \cite[Theorem 1.4]{LLR2025} with $a=2$, whose hypotheses hold:
\begin{enumerate}
\item $q_E=2$ and $r_E=2d/(d-2\alpha)<\infty$, since $d>2\alpha$.
\item $\frac d2(\frac12-\frac1{r_E})=\frac\alpha2$, so the pair is
sharp admissible.
\item $s_{2\alpha}(2,r_E)=0$, so the input Hilbert space is $L^2$.
\item The additional scalar hypothesis holds:
\[
 \frac d{r_E}+\frac2{q_E}
 =\frac d2-\alpha+1<d
 \quad(d\ge2,\ \alpha>0).
\]
\item Only the upper growth condition \eqref{frac_meas} is required.
\end{enumerate}
By the triangle inequality and \eqref{seg:scalar},
\begin{align}\label{seg:Einput}
 \|T_\Phi n\|_{Y_E}
 &\le\sum_j|n_j|\,
    \||e^{it\Delta}f_j|^2\|_{L_t^1(\rd\mu;L_x^{R_E})}\notag\\
 &\le C_{d,\alpha}M\sum_j|n_j|,
 \qquad Y_E=L_t^1(\rd\mu;L_x^{R_E}).
\end{align}
This proves \eqref{quadratic-strong} at $E$ itself.

\paragraph{Interpolation of the sequence and mixed norms.}
Apply \eqref{real-L} with $(q_0,q_1,q)=(P_A,1,P)$ and
$(B_0,B_1)=(L_x^{R_A},L_x^{R_E})$. This gives
\begin{equation}\label{seg:LP}
 (Y_A,Y_E)_{\theta,P}
 =L_t^P(\rd\mu;(L_x^{R_A},L_x^{R_E})_{\theta,P}).
\end{equation}
The outer time exponent and the real interpolation second index are
both $P$. The Lions--Peetre identity applies to a positive measure
without a time-translation assumption; see
\cite[Chapter VII, Section 1, Theorem 1.1]{LP1964} and
\cite[Lemma 5.4]{HT2026}. For distinct spatial Lebesgue exponents,
the inner interpolation space is $L_x^{R,P}$.

For sequence Lorentz spaces, the corresponding $K$-method identity is
\begin{equation}\label{seg:seq}
 (\ell^{b_0,u_0},\ell^{b_1,u_1})_{\theta,u}=\ell^{b,u}
 \quad(b_0\ne b_1).
\end{equation}
For our domain this can also be obtained by reiteration: choose $b>P_A$,
represent $\ell^{P_A,1}$ as $(\ell^1,\ell^b)_{\tau,1}$, and interpolate
this space with $\ell^1$. The new interpolation parameter is
$(1-\theta)\tau$, and its secondary index is the selected $u$.
The resulting reciprocal main exponent is
$(1-\theta)/P_A+\theta$. In particular, the second index $1$ of
$\ell^{P_A,1}$ is not retained after real interpolation.

Apply real interpolation to \eqref{seg:Ainput} and \eqref{seg:Einput}, choosing
the interpolation index $u=P$.
On the domain, \eqref{seg:seq} gives (note that $P_A\ne1$)
\begin{equation}\label{seg:domain}
 (\ell^{P_A,1},\ell^{1,1})_{\theta,P}
 =\ell^{P,P}=\ell^P.
\end{equation}
On the output, first apply \eqref{seg:LP} with the time measure $\mu$:
\begin{align}\label{seg:output}
 (Y_A,Y_E)_{\theta,P}
 &=L_t^P(\rd\mu;(L_x^{R_A},L_x^{R_E})_{\theta,P})\notag\\
 &=L_t^P(\rd\mu;L_x^{R,P}).
\end{align}
The spatial exponents are distinct, since
\[
 R_E-R_A=\frac{2\alpha^2}{(d-2\alpha)(d-\alpha)}>0.
\]
Thus the time exponent is strong $L^{P,P}$, the spatial secondary
exponent is $P$, and the sequence secondary exponent is also $P$.
Real interpolation of the operator norm gives
\begin{align*}
 \|T_\Phi n\|_{L_t^P(\rd\mu;L_x^{R,P})}
 &\le C_{d,\alpha,\theta}
 (M^{1/P_A})^{1-\theta}M^\theta\|n\|_{\ell^{P,P}}\\
 &=C_{d,\alpha,\theta}M^{1/P}\|n\|_{\ell^P},
\end{align*}
where the last exponent follows directly from \eqref{critical-segment-exponents}.
This proves \eqref{critical-refined}.

For consistency of the operator interpolation, both endpoint maps are
extensions of the same map on finite sequences. Since
$\ell^1\subset\ell^{P_A,1}$ and finite truncations converge in both
endpoint domains, their values agree on the intersection. Both target
spaces embed in locally integrable functions on bounded spatial sets
and finite-measure time sets, so equality of limits can be tested in
this common ambient space. The constants are uniform in the chosen
orthonormal family and in the number of coefficients.

For the spatial Lorentz embedding we note that
\begin{equation}\label{seg:sign}
 \frac1P-\frac1R
 =(1-\theta)\frac\alpha{d+\alpha}
       +\theta\frac{2\alpha}{d}>0.
\end{equation}
Thus $P<R$, and
\[
 L_x^{R,P}\hookrightarrow L_x^{R,R}=L_x^R.
\]
Applying this embedding to \eqref{critical-refined} proves \eqref{quadratic-strong} for
$0<\theta<1$. At $\theta=1$ use \eqref{seg:Einput} directly. Since we
do not control the interpolation constants as $\theta\to0$, the strong
estimate is obtained on $(A,E]$; at $A$ itself we have only the Lorentz
estimate \eqref{quadratic-A}.

For infinite sequences in $\ell^P$, apply the same estimate to differences
of finite truncations. Their outputs are Cauchy in the target Banach
space. For nonnegative coefficients the limit equals the pointwise
sum by monotone convergence; for arbitrary coefficients apply the same
argument to $|n_j|$ and use domination. This gives the stated extension.

\paragraph{Spatial Sobolev transfer.}
Put $q_A=2(d+\alpha)/d$ and $r_*(q)=2dq/(dq-4\alpha)$.
For $2\le q<q_A$ the sharp estimate just proved transfers to
$r_*(q)\le r<\infty$ at $s(q,r)=d(1/r_*(q)-1/r)$.
The following computation also tracks the stronger spatial Lorentz norm.

Write $r_0=r_*(q)$ and $s=d(1/r_0-1/r)$.
For nonnegative coefficients, introduce the finite-dimensional Hilbert
space $H_J=\ell^2(J)$ and the field
\[
 v(t,x)=(\sqrt{n_j}\,U(t)g_j(x))_{j\in J}.
\]
Since the Riesz potential $I_s=|D|^{-s}$ has a positive scalar kernel
when $0<s<d$,
\[
 \|I_sv(t,x)\|_{H_J}\le I_s(\|v(t,\cdot)\|_{H_J})(x).
\]
For finite $r>r_0$, we have $0<s<d/r_0$ and $r_0>1$. The scalar
strong Riesz-potential inequality therefore gives
\begin{equation}\label{seg:vector-sobolev}
 \|I_sv(t)\|_{L_x^r(H_J)}\le C_{d,s,r_0}\|v(t)\|_{L_x^{r_0}(H_J)}.
\end{equation}
For $s=0$ this is just the identity. Apply the sharp-segment strong
estimate at $(q,r_0)$, and square the field norm:
\begin{align*}
 \left\|\sum_jn_j|U(t)I_sg_j|^2\right\|_{L_t^{q/2}L_x^{r/2}}
 &=\|I_sv\|_{L_t^qL_x^r(H_J)}^2\\
 &\le C\|v\|_{L_t^qL_x^{r_0}(H_J)}^2
 \le C M^{2/q}\|n\|_{\ell^{q/2}}.
\end{align*}
This proves \eqref{quadratic-strong}. Complex coefficients are handled
by domination with $|n_j|$. The Hilbert-space dimension does not enter
the constant. Here $(g_j)$ is orthonormal in $L^2$; equivalently,
$(I_sg_j)$ is orthonormal in $\dot H^s$.

The spatial Lorentz indices are preserved by the transfer. For
$2<q<q_A$, the preceding sharp estimate has spatial density norm
$L^{r_0/2,q/2}$; its square-root field norm is $L^{r_0,q}(H_J)$.
Scalar Lorentz Riesz-potential bounds preserve the secondary index:
\[
 I_s:L^{r_0,q}\longrightarrow L^{r,q}.
\]
These bounds follow by real interpolation of strong Riesz-potential
bounds at nearby reciprocal exponents, with the same $s$.
Positivity gives their Hilbert-valued version. Squaring returns
$L^{r/2,q/2}$ for the density, while the sequence space remains
$\ell^{q/2,q/2}$. As $q<r_0\le r$, the spatial Lorentz embedding is in
the required direction. At $q=2$, the scalar input at $E$ gives the
strong spatial norm directly; no spatial secondary index $1$ is claimed.
At $q=q_A$, this spatial transfer alone retains sequence index $1$.
The strong estimate for $r>r_A$ is supplied instead by the interpolated
inputs below.

\subsection{Asymmetric Hilbert--Schmidt estimates}
\label{newq:section}

We extend the direct $q=4$ Hilbert--Schmidt method of
\cite[Section 3]{Gong2026} to this measure and to asymmetric Sobolev
orders. Orthogonal initial-space frequency blocks then give the critical
strong estimate for every $q>4$ in the strict finite-$r$ region.

\begin{thm}[Schatten estimates for $q\ge4$]\label{newq:qge4}
Suppose
\[
 4\le q<\infty,\qquad \frac{dq}{d-\alpha}<r<\infty,
 \qquad s=\frac d2-\frac dr-\frac{2\alpha}{q}.
\]
With $p=(q/2)'$ and $c=(r/2)'$,
\begin{equation}\label{newq:dual}
 \begin{split}
 \|B_s(V)\|_{\Schatten^p}
 &\le C_{d,\alpha,q,r}M^{2/q}\|V\|_{L_t^p(\rd\mu;L_x^c)},\\
 B_s(V)&=|D|^{-s}\int U(t)^*V(t)U(t)\,\rd\mu(t)|D|^{-s}.
 \end{split}
\end{equation}
By Schatten duality, for every $L^2$-orthonormal family $(g_j)$,
\begin{equation}\label{newq:density}
 \left\|\sum_j n_j|U(t)|D|^{-s}g_j|^2\right\|_
 {L_t^{q/2}(\rd\mu;L_x^{r/2})}
 \le C_{d,\alpha,q,r}M^{2/q}\|n\|_{\ell^{q/2}}.
\end{equation}
This is the $q\ge4$ density estimate used in
\eqref{quadratic-strong}.
The integral operator is defined by extension of its form on
finite-annulus Fourier data. All frequency restrictions have the same
constant.
\end{thm}

\begin{lem}\label{newq:mixed-hs}
Let
\[
 1<c_4<\frac{2d}{d+\alpha},\qquad
 s_4=\frac d{2c_4}-\frac\alpha2,\qquad \gamma_4=d-2s_4.
\]
For $0\le\delta\le\delta_0$, where $\delta_0>0$ is sufficiently small,
the two mixed-order operators
\[
 C_{s_4,\pm\delta}(V)
 =|D|^{-(s_4\pm\delta)}
   \int U(t)^*V(t)U(t)\rd\mu(t)|D|^{-(s_4\mp\delta)}
\]
satisfy
\begin{equation}\label{newq:mixed-bound}
 \|C_{s_4,\pm\delta}(V)\|_{\Schatten^2}
 \le C M^{1/2}\|V\|_{L_t^2(\rd\mu;L_x^{c_4})}.
\end{equation}
The choice of $\delta_0$ and the constant may be uniform as $c_4$ varies
over a small compact neighborhood of a fixed permitted exponent.
\end{lem}

\begin{proof}
Choose
\begin{equation}\label{newq:delta-margin}
 0<\delta_0<\min\{s_4,d/2-s_4,(d/4)(2/c_4-1)\}.
\end{equation}
Put $s_\pm=s_4\pm\delta$ and $\gamma_\pm=\gamma_4\mp2\delta$.
Choose a nonnegative smooth annular partition
$\sum_{j\in\mathbb Z}\psi(2^j\xi)=1$ for $\xi\ne0$.
The kernel for $U(\tau)|D|^{-2s_\pm}$ has dyadic pieces
\[
 I_{j,\pm}(\tau,z)=2^{-j\gamma_\pm}\int
 e^{i(2^{-j}z\cdot\eta-2^{-2j}\tau|\eta|^2)}
 \psi(\eta)|\eta|^{-2s_\pm}\rd\eta.
\]
On the stationary indices $2^j|z|\asymp|\tau|$, their number is $O(1)$
and the trivial and nondegenerate-Hessian bounds give
\[
 S_\pm:=\sum_{\substack{j\text{ stationary}}}|I_{j,\pm}|
 \lesssim\min\{(|z|/|\tau|)^{\gamma_\pm},
 |z|^{-2s_\pm}|\tau|^{-(\gamma_\pm-d/2)}\}.
\]
Outside a fixed enlargement of that region, repeated integration by
parts gives
\[
 |I_{j,\pm}|\lesssim_N2^{-j\gamma_\pm}
       \min\{1,(2^{-j}|z|)^{-N},(2^{-2j}|\tau|)^{-N}\}.
\]
Geometric summation with $N>\max(\gamma_+,\gamma_-)$ consequently yields
\[
 H_\pm:=\sum_{\substack{j\text{ nonstationary}}}|I_{j,\pm}|
 \lesssim\min\{|z|^{-\gamma_\pm},|\tau|^{-\gamma_\pm/2}\}.
\]
All these bounds dominate every finite dyadic frequency restriction.

The product of the two absolute kernel sums is bounded near
$|z|^2\le|\tau|$ by $C|\tau|^{-\gamma_4}$. In the far region its four terms
are bounded by
\begin{equation}\label{newq:mixed-majorant}
 \begin{split}
 &|z|^{-2\gamma_4},\qquad |z|^{-4s_4}|\tau|^{-(\gamma_4-2s_4)},\\
 &|z|^{-(d+4\delta)}|\tau|^{-(\gamma_4-d/2-2\delta)},\qquad
 |z|^{-(d-4\delta)}|\tau|^{-(\gamma_4-d/2+2\delta)}.
 \end{split}
\end{equation}
They are respectively $H_+H_-$, $S_+S_-$, $S_+H_-$ and $S_-H_+$.
No cancellation or restriction on the difference of the two frequency
indices is needed.

For each positive majorant split $|\tau|$ into
$2^k\le|\tau|<2^{k+1}$. Let $Q_k(F,G)$ be its bilinear form against
$F(t,x)G(t',y)$. For spatial exponents $a,b$, write
$u=1/a$, $v=1/b$, and $1/h=2-u-v$. Spatial Young after scaling by
$2^{k/2}$, and then Frostman Schur on $L^2(\rd\mu)$, give
\begin{equation}\label{newq:shell-bound}
 \begin{gathered}
 |Q_k(F,G)|\le CM2^{-k\nu(a,b)}
 \|F\|_{L_t^2(\rd\mu;L_x^a)}
 \|G\|_{L_t^2(\rd\mu;L_x^b)},\\
 \nu(a,b)=\frac d2(u+v)-\frac d{c_4}.
 \end{gathered}
\end{equation}
Indeed the spatial factor is
$2^{k[-\gamma_4+d-d(u+v)/2]}$, and the time-shell indicator has both row and
column integrals at most $CM2^{\alpha k}$. The exponent reduces to
\eqref{newq:shell-bound} because $-\gamma_4+d+\alpha=d/c_4$.
The near spatial kernel is a rescaled ball indicator. For the far
kernels, spatial Young requires
\[
 2\gamma_4h>d,\qquad4s_4h>d,\qquad(d-4\delta)h>d;
\]
the $d+4\delta$ condition is weaker than the last one.
Using $2s_4=d/c_4-\alpha$ and $\gamma_4=d-d/c_4+\alpha$, these become
\[
 u+v>2/c_4-2\alpha/d,\qquad
 u+v>2-2/c_4+2\alpha/d,\qquad
 u+v>1+4\delta/d,
\]
together with $1<u+v<2$. They hold in an open neighborhood of
$(u,v)=(1/c_4,1/c_4)$ by \eqref{newq:delta-margin}.

Choose $1/a_0=1/c_4-\varepsilon$ and
$1/a_1=1/c_4+2\varepsilon$ inside this neighborhood. The sequence
$Q=(Q_k)_k$ has the three bilinear bounds with output weights
$-d\varepsilon$ at $(a_0,a_0)$ and $d\varepsilon/2$ at
$(a_0,a_1)$ and $(a_1,a_0)$. Apply the three-bound bilinear real
interpolation theorem \cite[Section 3.13, Exercise 5(b)]{BL1976}
with input parameters $(1/3,2),(1/3,2)$ and output parameters $(2/3,1)$.
The output space is $\ell^1$, since the interpolated weight is zero.
The matching-outer-index Lions--Peetre identity gives the input spaces
$L_t^2(\rd\mu;L_x^{c_4,2})$. Since $c_4<2$,
$L_x^{c_4}\subset L_x^{c_4,2}$, so
\[
 \sum_k|Q_k(F,G)|\le CM
 \|F\|_{L_t^2(\rd\mu;L_x^{c_4})}
 \|G\|_{L_t^2(\rd\mu;L_x^{c_4})}.
\]
Time retains second index $2$, space passes through second index $2$,
and the shell sum has second index $1$.

For a finite set $J$, put $m_J(\xi)^2=\sum_{j\in J}\psi(2^j\xi)$,
and let $C^J_{s_4,\delta}(V)$ be the form with $m_J$ on both initial
sides. Write $K^J_{2u}$ for the finite annular kernel sum of
$U(\tau)|D|^{-2u}m_J(D)^2$. For bounded supported $V$, the
squared Fourier-kernel norm is, up to a positive dimensional constant,
\[
 \begin{aligned}
 \|C^J_{s_4,\delta}(V)\|_{\Schatten^2}^2
 =c_d\iiiint &\overline{V(t,x)}V(t',y)
 K^J_{2(s_4+\delta)}(t-t',x-y)\\
 &\cdot\overline{K^J_{2(s_4-\delta)}(t-t',x-y)}
 \,\rd x\rd y\rd\mu(t)\rd\mu(t').
 \end{aligned}
\]
Taking absolute values, apply the positive majorants with $F=G=|V|$.
This proves
\eqref{newq:mixed-bound}, also for complex $V$ and either sign of
$\delta$. Initial weighted operators may be restricted to a bounded
time interval containing the potential support; their constants are
independent of this interval. The time diagonal is null because the
Frostman condition makes $\mu$ non-atomic. Weak limits in the Hilbert
space $\Schatten^2$ are uniquely identified by their stabilized forms on
finite-annulus Fourier data. Density in $L_t^2L_x^{c_4}$ completes the
operator construction. In particular the case $\delta=0$ is a global,
cutoff-uniform Hilbert--Schmidt estimate.
\end{proof}

For $q=4$ the identity $2s=d/c-\alpha$ and the case $\delta=0$ of
Lemma~\ref{newq:mixed-hs} prove \eqref{newq:dual}. Its density formulation
follows by Schatten duality as detailed at the end of the proof below.

\subsection{Orthogonal frequency blocks for \texorpdfstring{$q>4$}{q>4}}
\begin{proof}[Proof of Theorem~\ref{newq:qge4} for $q>4$]
Put
\[
 p=\frac q{q-2}\in(1,2),\qquad c=\frac r{r-2},\qquad
 \vartheta=2-2/p=4/q.
\]
Since $r>dq/(d-\alpha)>q$, one has $1<c<p<2$.
Let $\mathsf P_k$ be the orthogonal Fourier projection onto
$2^k\le|\xi|<2^{k+1}$, and put $N_k=2^k$. Thus
\[
 \mathsf P_k\mathsf P_l=0\quad(k\ne l),\qquad
 \sum_{k\in\mathbb Z}\mathsf P_k=I
\]
in the strong operator topology. These projections are distinct from
the smooth Littlewood--Paley multipliers $P_k$ used above.
For bounded supported $V$ define the initial-space frequency blocks
\[
 B_{kl}(V)=\mathsf P_k|D|^{-s}\int U(t)^*V(t)U(t)\rd\mu(t)|D|^{-s}\mathsf P_l.
\]
Factoring through multiplication by $\sqrt{|V|}$ and its complex phase
gives
\begin{equation}\label{newq:block-trace}
 \|B_{kl}(V)\|_{\Schatten^1}
 \lesssim(N_kN_l)^{d/2-s}\|V\|_{L_t^1(\rd\mu;L_x^1)}.
\end{equation}
Indeed, by Plancherel's theorem,
\[
 \|\sqrt{|V|}U|D|^{-s}\mathsf P_k\|_{\Schatten^2}^2
 =C\int|V|\rd\mu\rd x
       \int_{2^k\le|\xi|<2^{k+1}}|\xi|^{-2s}\rd\xi
 \lesssim N_k^{d-2s}\|V\|_{L_t^1L_x^1}.
\]

For any permitted $c_4$ put $s_4=d/(2c_4)-\alpha/2$.
Choose the sign of $\delta$ opposite to $k-l$ in
Lemma~\ref{newq:mixed-hs}. Bounded frequency multipliers on its two sides
yield
\begin{equation}\label{newq:block-hs}
 \|B_{kl}(V)\|_{\Schatten^2}
 \le CM^{1/2}2^{-\delta|k-l|}(N_kN_l)^{s_4-s}
       \|V\|_{L_t^2(\rd\mu;L_x^{c_4})}.
\end{equation}
For example, when $k\ge l$ use
\[
 B_{kl}(V)=\mathsf P_k|D|^{s_4-\delta-s}
 C_{s_4,-\delta}(V)|D|^{s_4+\delta-s}\mathsf P_l.
\]

Choose nearby spatial potential exponents $c_i$ around $c$ and set
\[
 \frac1{c_i}=1-\vartheta+\frac\vartheta{c_{4,i}}.
\]
At the central target this gives
\[
 \frac1{c_4}=1-\frac q{2r},\qquad
 r_4=4r/q>4d/(d-\alpha),\qquad c_4=r_4/(r_4-2).
\]
Thus a common $\delta>0$ in \eqref{newq:block-hs} is available for all
these neighboring exponents. Complex interpolation of the linear map
$V\mapsto B_{kl}(V)$, at the fixed order $s$, between
\eqref{newq:block-trace} and \eqref{newq:block-hs} gives
\begin{equation}\label{newq:block-p}
 \|B_{kl}(V)\|_{\Schatten^p}
 \le CM^{2/q}2^{-\vartheta\delta|k-l|}(N_kN_l)^{a_i}
       \|V\|_{L_t^p(\rd\mu;L_x^{c_i})},
 \quad a_i=\frac d2\left(\frac1{c_i}-\frac1c\right).
\end{equation}
The exponent is obtained by the substitution
\[
 \begin{aligned}
 a_i&=(1-\vartheta)(d/2-s)
       +\vartheta\big(d/(2c_{4,i})-\alpha/2-s\big)\\
 &=\frac d{2c_i}-\frac{\vartheta\alpha}2-s
  =\frac d2\left(\frac1{c_i}-\frac1c\right),
 \end{aligned}
 \qquad 2s=d/c-\vartheta\alpha.
\] Both mixed Lebesgue input norms and the
Schatten output are strong in this complex interpolation step.

Now choose $1/c_0=1/c-\varepsilon$ and $1/c_1=1/c+\varepsilon$, with
$d\varepsilon/2<\vartheta\delta/2$. Fix $n\in\mathbb Z$ and consider
$L_n(V)=(B_{l+n,l}(V))_{l\in\mathbb Z}$. Since
$(N_{l+n}N_l)^{a_i}=2^{2la_i}2^{na_i}$, the residual factor $2^{na_i}$
satisfies
\[
 2^{-\vartheta\delta|n|}2^{na_i}
 \le2^{-(\vartheta\delta-d\varepsilon/2)|n|}
 \le2^{-\kappa|n|},\qquad\kappa=\vartheta\delta/2.
\]
For the weighted norm
$\|(A_l)\|_{\ell_\gamma^\infty(\Schatten^p)}
=\sup_l2^{l\gamma}\|A_l\|_{\Schatten^p}$,
\eqref{newq:block-p} gives
\[
 L_n:L_t^p(\rd\mu;L_x^{c_i})\longrightarrow
 \ell_{\gamma_i}^\infty(\Schatten^p),\qquad
 \gamma_0=d\varepsilon,\quad\gamma_1=-d\varepsilon,
\]
with norm at most $CM^{2/q}2^{-\kappa|n|}$.

Apply linear real interpolation at parameter $1/2$ and second index
$p$. The matching-index and weighted-sequence formulas give
\[
 (L_t^pL_x^{c_0},L_t^pL_x^{c_1})_{1/2,p}
 =L_t^p(\rd\mu;L_x^{c,p}),
\]
\[
 (\ell_{d\varepsilon}^\infty(\Schatten^p),
  \ell_{-d\varepsilon}^\infty(\Schatten^p))_{1/2,p}
 =\ell^p(\Schatten^p).
\]
Because $c<p$, $L_x^c=L_x^{c,c}\subset L_x^{c,p}$, and consequently
\begin{equation}\label{newq:gap-sum}
 \left(\sum_l\|B_{l+n,l}(V)\|_{\Schatten^p}^p\right)^{1/p}
 \le CM^{2/q}2^{-\kappa|n|}\|V\|_{L_t^p(\rd\mu;L_x^c)}.
\end{equation}
Time has second index $p$, the intermediate space second index is $p$,
and the frequency sequence second index is $p$.

At a fixed $n$, the blocks have mutually orthogonal domains
$\operatorname{Ran}\mathsf P_l$ and mutually orthogonal ranges
$\operatorname{Ran}\mathsf P_{l+n}$. Therefore
\[
 \left\|\sum_{l\in K}B_{l+n,l}\right\|_{\Schatten^p}^p
 =\sum_{l\in K}\|B_{l+n,l}\|_{\Schatten^p}^p
\]
for every finite $K$: cross terms in the adjoint-times-operator product
vanish, and the remaining positive terms act on disjoint domains.
Thus $\mathcal B_n(V)=\sum_lB_{l+n,l}(V)$ converges in $\Schatten^p$ with the
bound \eqref{newq:gap-sum}. Summing $n$ absolutely in $\Schatten^p$
defines $B_s(V)$ and gives \eqref{newq:dual}.

On finite-annulus Fourier test functions this operator has the original
integral form, since only finitely many block pairs contribute. Density
extends $V$ to $L_t^pL_x^c$. Every finite frequency restriction is
$Q B_s(V)Q$ and hence has the same norm bound. Rank-one duality gives
the full scalar $L_t^q(\rd\mu;L_x^r)$ extension from finite-annulus
data. For finite orthonormal families apply Schatten H\"older to
$\gamma=\sum_jn_j|f_j\rangle\langle f_j|$, with
$\|\gamma\|_{\Schatten^{q/2}}=\|n\|_{\ell^{q/2}}$.
Uniform cutoff convergence identifies the density forms; this does not
require the cutoff vectors to be orthonormal. Positive monotone convergence and
domination by absolute coefficients handle infinite families and complex
coefficients. This proves \eqref{newq:density}; the same argument applies
to the $q=4$ Hilbert--Schmidt case.
\end{proof}

\begin{rem}\label{newq:boundaries}
The proof uses strict inequalities at three places. At $r=dq/(d-\alpha)$,
$c_4=2d/(d+\alpha)$ and the stationary far kernel loses its open
integrability neighborhood. At $r=\infty$, $c=1$ and the neighboring
$c_i$ cross the Banach boundary. At $q=\infty$, $\vartheta=0$ and the
frequency-gap decay disappears. None of these boundaries is included
by taking a limit. The temporal trace/Hilbert--Schmidt interpolation
in the proof of Lemma~\ref{lem-three-bounds}~(3) requires
$\frac1{w_1}+\frac1{w_2}>\frac1p$, whereas the critical target lies on the
boundary $\frac1{w_1}=\frac1{w_2}=\frac1{2p}$; this is why that argument gives
only $\beta<q/2$. The present argument instead uses asymmetric smoothing,
linear real interpolation and orthogonal initial-space blocks.
\end{rem}

\subsection{Interpolation across finite spatial exponents}
\label{newq:lowq-section}

We record the interpolation argument explicitly, since the Sobolev order
changes with $(q,r)$ and the density map itself is then not fixed.

\begin{lem}\label{newq:varying-order}
Suppose two strong density estimates hold at finite spatial exponents
$r_i>2$, with $P_i=q_i/2\ge1$, $R_i=r_i/2>1$, coefficient exponent
$P_i$, and
\[
 s_i=\frac d2-\frac d{2R_i}-\frac\alpha{P_i},\qquad i=0,1.
\]
Equivalently, writing $a_i=P_i'$ and $b_i=R_i'$, assume
\[
 \|B_{s_i}(V)\|_{\Schatten^{a_i}}
 \le C_iM^{1/P_i}\|V\|_{L_t^{a_i}(\rd\mu;L_x^{b_i})}.
\]
For $P_i=1$, this is interpreted as an operator-norm estimate with
$a_i=\infty$. If $0<\theta<1$ and
\[
 \frac1P=\frac{1-\theta}{P_0}+\frac\theta{P_1},\qquad
 \frac1R=\frac{1-\theta}{R_0}+\frac\theta{R_1},\qquad
 s=(1-\theta)s_0+\theta s_1,
\]
then, whenever $1<P<\infty$, the strong estimate at $(q,r)=(2P,2R)$
holds with coefficient exponent $P$ and constant
$C_0^{1-\theta}C_1^\theta M^{1/P}$.
\end{lem}

\begin{proof}
Write $a=P'$, $b=R'$. Normalize a bounded simple rectangular potential
$V$ with finite support by $\|V\|_{L_t^aL_x^b}=1$ and put
$F(t)=\|V(t)\|_{L_x^b}$. Such potentials are dense at the finite
interior target. On the strip define
\[
 \frac1{a(z)}=\frac{1-z}{a_0}+\frac z{a_1},\qquad
 \frac1{b(z)}=\frac{1-z}{b_0}+\frac z{b_1},\qquad
 s(z)=(1-z)s_0+zs_1,
\]
and, on $V\ne0$, set
\[
 V_z(t,x)=\frac{V(t,x)}{|V(t,x)|}
 F(t)^{a/a(z)-b/b(z)}|V(t,x)|^{b/b(z)}.
\]
Set $V_z=0$ on $V=0$. The positive values of $F$ and $|V|$ lie in
finite sets, so the powers define analytic functions without growth
problems. At $z=\theta$, $V_\theta=V$, and on either boundary
\[
 \|V_{i+iy}(t)\|_{L_x^{b_i}}=F(t)^{a/a_i},\qquad
 \|V_{i+iy}\|_{L_t^{a_i}L_x^{b_i}}=1.
\]
When $a_i=\infty$, the first expression is $1$ on the nonzero time
support and the same conclusion follows.

Choose a smooth multiplier $Q_N(D)$ supported in a finite annular union,
with $\|Q_N\|_{2\to2}\le1$, $Q_N=1$ on $N^{-1}\le|\xi|\le N$,
and $Q_N\to I$ strongly. Consider
\[
 B_N(z)=Q_N|D|^{-s(z)}
       \int U(t)^*V_z(t)U(t)\rd\mu(t)|D|^{-s(z)}Q_N.
\]
Finite frequency and potential supports make this a trace-class analytic
family on each closed substrip. On $z=i+iy$, the imaginary multiplier
$|D|^{-iy(s_1-s_0)}$ is unitary on the initial $L^2$ space and commutes
with $Q_N$ and $U(t)$. It can be extracted on both sides:
\[
 B_N(i+iy)=|D|^{-iy(s_1-s_0)}
 Q_NB_{s_i}(V_{i+iy})Q_N|D|^{-iy(s_1-s_0)}.
\]
The two imaginary powers have the same sign; unitary invariance of
singular values applies to separate left and right unitaries. Thus
\[
 \|B_N(i+iy)\|_{\Schatten^{a_i}}\le C_iM^{1/P_i}
\]
uniformly in $N,y$. Complex interpolation of the analytic family into
strong Schatten classes gives
\[
 \|B_N(\theta)\|_{\Schatten^a}
 \le C_0^{1-\theta}C_1^\theta M^{1/P}
       \|V\|_{L_t^a(\rd\mu;L_x^b)}.
\]
This is the strong mixed-norm complex interpolation: time, space and
Schatten second indices are respectively $a,b,a$. After trace duality
the sequence, density-time and density-space second indices are
$P,P,R$.

Since $1<a<\infty$, the uniform $\Schatten^a$ bounds have weakly
convergent subsequences. On finite-annulus Fourier data their matrix
entries stabilize to the target integral form, which identifies a unique
limit with the same bound. Potential density extends the operator.
Rank-one duality identifies the scalar target norm, and Schatten trace
duality gives the estimate for every original orthonormal family.
Uniform frequency restrictions and strong scalar approximation remove
the cutoffs without assuming that projected vectors are orthonormal.
\end{proof}

\begin{thm}\label{newq:full-finite-interior}
Suppose $a=2$, $d\ge2$ and $d>2\alpha$. Let $2<q,r<\infty$ satisfy
\begin{equation}\label{newq:strict-region}
 \frac d r+\frac{2\alpha}q<\frac d2,
 \qquad \frac1q>\frac d{d-\alpha}\frac1r.
\end{equation}
Then \eqref{newq:density} holds at
$s=d/2-d/r-2\alpha/q$ with coefficient exponent $\beta=q/2$.
Thus the strict finite-$r$ part of $ODE_{d/2}A_{d/2}$ is covered.
\end{thm}

\begin{proof}
For $q\ge4$, the second inequality in \eqref{newq:strict-region} is
$r>dq/(d-\alpha)$, and Theorem~\ref{newq:qge4} applies.
For $2<q<q_A$, the first inequality is $r>r_*(q)$, and
the finite spatial Sobolev transfer proved above applies. It remains to consider
$q_A\le q<4$.

Write $x=1/r$, $y=1/q$, and set
\[
 y_A=\frac d{2(d+\alpha)},\quad
 x_A=\frac{d-\alpha}{2(d+\alpha)},\quad
 y_4=\frac14,\quad x_4=\frac{d-\alpha}{4d},\quad
 k=\frac{d-\alpha}d,
\]
and $h(y)=1/2-(2\alpha/d)y$. In the remaining range
$y_4<y\le y_A$, the strict region is exactly $0<x<ky$.
Choose $y_0>y_A$ sufficiently close to $y_A$, and let
\[
 \theta=\frac{y_0-y}{y_0-y_4},\qquad
 X(y_0)=(1-\theta)h(y_0)+\theta x_4.
\]
As $y_0\downarrow y_A$, $X(y_0)\to ky$, because $A$ and
$(x_4,y_4)$ lie on $x=ky$. Hence choose $y_0$ with $X(y_0)>x$ and put
$\lambda=x/X(y_0)\in(0,1)$. The two points
\[
 (x_0,y_0)=(\lambda h(y_0),y_0),\qquad
 (x_1,y_4)=(\lambda x_4,y_4)
\]
have finite spatial exponents and are respectively in the low-$q$
strong region and the strict $q=4$ strong segment. Moreover
$(x,y)=(1-\theta)(x_0,y_0)+\theta(x_1,y_4)$ exactly. Both endpoint
sequence exponents equal their density time exponents. Apply
Lemma~\ref{newq:varying-order}; the order
$s=d/2-dx-2\alpha y$ is affine, so the interpolated order is the target order
and the strong sequence exponent $q/2$. This includes $q=q_A$ with
$r>r_A$; the point $A$ itself is not used as an endpoint.
\end{proof}

\begin{rem}\label{newq:finite-scope}
The sharp segment and its spatial Sobolev transfer also include the
finite-$r$ admissible boundary for $2\le q<q_A$. Together with the
strict-region argument, they prove part (2) of
Theorem~\ref{thm-quadratic-critical}. At $q=q_A$ the strict condition
$r>r_A$ is supplied by the interpolated strong inputs, while the point
$A$ itself retains the Lorentz sequence estimate.
\end{rem}

\subsection{The boundary \texorpdfstring{$r=\infty$}{r=infinity}}
For every $1<P=q/2<\infty$, let $s_\infty=d/2-\alpha/P$.
We prove part (3) of Theorem~\ref{thm-quadratic-critical}:
\begin{equation}\label{od:ODweak}
 \left\|\sum_j n_j|U(t)|D|^{-s_\infty}g_j|^2\right\|_
 {L_t^{P,\infty}(\rd\mu;L_x^\infty)}
 \le C_{d,\alpha,P}M^{1/P}\|n\|_{\ell^{P,1}}.
\end{equation}
Here both critical Lorentz losses are explicit: time has secondary
index $\infty$ and the coefficients have secondary index $1$.
The following proof uses only $d/2>\alpha$ and the upper growth bound.

Take smooth annular spatial projections $\Pi_k$ forming a homogeneous
partition of unity. The annular kernel for $U(t)I_s\Pi_k$ composed with
its adjoint satisfies
\[
 |K_k(t-\tau,x-y)|\le C2^{k(d-2s)}
       (1+2^{2k}|t-\tau|)^{-d/2}.
\]
Growth and dyadic time shells give
\[
 \sup_t\int(1+2^{2k}|t-\tau|)^{-d/2}\,\rd\mu(\tau)
 \le C_{d,\alpha}M2^{-2\alpha k},
\]
because $d/2>\alpha$. By Schur and the scalar $TT^*$ argument,
\begin{equation}\label{od:annular-scalar}
 \|U I_s\Pi_k g\|_{L_t^2(\rd\mu;L_x^\infty)}^2
 \le C M2^{k(d-2s-2\alpha)}\|g\|_2^2.
\end{equation}
Since the annular outputs are continuous in $x$, their $L_x^\infty$ norms
are attained by testing against $L_x^1$ functions, and time measurability
follows from suprema over rational points; hence the $TT^*$ argument
applies.

For a set $J$ of $L$ orthonormal vectors, Bessel's inequality also gives
\begin{equation}\label{od:annular-bessel}
 \sup_{t,x}\sum_{j\in J}|U I_s\Pi_k g_j(t,x)|^2
 \le C2^{k(d-2s)}.
\end{equation}
First work with a finite annular truncation and with coefficients
$n_j=\ind_J(j)$. Split the Hilbert-valued field into $k\le m$ and
$k>m$. For $s=s_\infty$, Minkowski and
\eqref{od:annular-bessel}--\eqref{od:annular-scalar} imply
\[
 \begin{gathered}
 \|v_{\le m}\|_{L_t^\infty L_x^\infty(H_J)}
 \le C2^{\alpha m/P},\\
 \|v_{>m}\|_{L_t^2(\rd\mu;L_x^\infty(H_J))}
 \le C(ML)^{1/2}2^{-\alpha m(1-1/P)}.
 \end{gathered}
\]
Given $z>0$, choose $m$ so that the first bound is at most
$\sqrt z/2$ and $2^{2\alpha m/P}$ is comparable to $z$.
Chebyshev's inequality gives
\[
 \begin{aligned}
 \mu\{t:\|v(t)\|_{L_x^\infty(H_J)}^2>z\}
 &\le CMLz^{-1}2^{-2\alpha m(1-1/P)}\\
 &=CMLz^{-1}\big(2^{2\alpha m/P}\big)^{-(P-1)}
 \lesssim MLz^{-P}.
 \end{aligned}
\]
Thus the restricted bound is $C M^{1/P}L^{1/P}$.
The finite truncation may be removed as follows: for $q>2$, summing
the annular kernels gives $|K_s(t,x)|\lesssim |t|^{-\alpha/P}$.
The Frostman fractional integral inequality, whose exponent is strictly
less than $\alpha$, yields the scalar bounded map
$UI_s:L^2\to L_t^{2P}(\rd\mu;L_x^\infty)$.
Consequently annular truncations converge in that space for each
initial vector; for a finite family their densities converge in
$L_t^P L_x^\infty$. The restricted weak bound survives this limit.

Finally use the dyadic layer decomposition of $|n_j|$. The equivalent
Banach norm on $L^{P,\infty}$, valid for $P>1$, permits summing the
restricted bounds, and
\[
 \sum_{a\in\mathbb Z}2^a
       \#\{j:|n_j|>2^a\}^{1/P}\simeq\|n\|_{\ell^{P,1}}.
\]
This proves \eqref{od:ODweak}. No strong $\ell^P\to L_t^P L_x^\infty$
estimate follows from this proof. At $P=1$ the high-frequency sum has
zero decay, so the argument does not establish the point $D$.

\begin{rem}\label{newq:OD-scope}
Estimate \eqref{od:ODweak} leaves open strong type at $\beta=P$.
The low/high-frequency proof gives only
$\mu\{t:\|\rho(t)\|_\infty>z\}\le CMLz^{-P}$ for an $L$-element
orthonormal family. Integrating this distributional bound to estimate
the strong $L^P$ norm leaves the critical integral $\int dz/z$;
the coefficient layer sum gives $\ell^{P,1}$. These are limitations
of the argument, not an orthonormal-family counterexample at $\beta=P$.

Concrete lower-bound tests do show that $\beta\le P$ is necessary for
a bound uniform over the whole class of upper-Frostman measures. For
each $\alpha\in(0,1)$ choose the fixed separated two-map Cantor
probability measure of ratio $2^{-1/\alpha}$; for $\alpha=1$ use
$\mathbf1_{[0,1]}dt$. Both have upper growth and lower growth at zero.
Choose a nonnegative smooth annular Fourier profile $\varphi$ and
$N_k=\Lambda^k$ with fixed sufficiently large $\Lambda$, and set
$\widehat f_k(\xi)=N_k^{-d/2}\varphi(\xi/N_k)$, normalized in $L^2$.
Their disjoint Fourier supports make them orthonormal, and
\[
 |U(t)|D|^{-s}f_k(0)|\ge cN_k^{\alpha/P}
 \quad\text{when }0\le t\le c_0N_k^{-2}.
\]
Write $c_*$ for a lower growth constant at zero. Choose
$0<a<b<c_0$ with $Ma^\alpha\le c_*b^\alpha/2$, and then choose
$\Lambda$ sufficiently large. Since these measures are supported in
$[0,1]$,
\[
 \mu([aN_k^{-2},bN_k^{-2}])
 \ge(c_*b^\alpha-Ma^\alpha)N_k^{-2\alpha}
 \ge\tfrac12c_*b^\alpha N_k^{-2\alpha}.
\]
The intervals are disjoint when $\Lambda^2>b/a$. Hence the equal-weight
$K$-element density has strong norm at least $cK^{1/P}$, whereas its
$\ell^\beta$ norm is $K^{1/\beta}$. This disproves $\beta>P$ on one
fixed admissible measure. For completeness the upper envelope is
\[
 \sup_x|U(t)|D|^{-s}f_k(x)|^2
 \lesssim N_k^{2\alpha/P}(1+N_k^2|t|)^{-d}.
\]
On $|t|\asymp N_m^{-2}$, divide the density by $N_m^{2\alpha/P}$.
For arbitrary nonnegative coefficients the result is bounded by
\[
 \sum_{k\le m}n_k\Lambda^{-2\alpha(m-k)/P}
 +\sum_{k>m}n_k\Lambda^{-(2d-2\alpha/P)(k-m)}.
\]
Both scale kernels are summable. The upper Frostman mass of this
time shell is $O(MN_m^{-2\alpha})$, so discrete Young in $\ell^P$
gives the matching upper bound $CM^{1/P}\|n\|_{\ell^P}$.
In particular the equal-weight family does not disprove equality.
The time-shift construction in Section~\ref{sec-necessary} gives
the corresponding necessary condition for all measures satisfying
\eqref{lower-regular}. Neither construction proves failure at equality.
\end{rem}

\subsection{Endpoint limitations}\label{sec-quadratic-limitations}\paragraph{Weak mixed-norm interpolation.}
One might hope that two estimates with output
$L_t^{P_i,\infty}(\rd\mu;L_x^{R_i})$ at the same $s$ could be
real-interpolated with index $P$ into $L_t^P(\rd\mu;L_x^{R,P})$.
There is no such general embedding, even when $P<R$.
We give an explicit counterexample.

Choose $P_0\ne P_1$, $R_0\ne R_1$, all finite and greater than one,
so that
\[
 \frac1{P_i}+\frac d{2\alpha R_i}=\kappa\quad(i=0,1).
\]
These pairs lie on a common line of constant $s$, with
$\kappa=d/(2\alpha)-s/\alpha$.
Normalize a nonzero atomless time measure to total mass one, and choose
disjoint measurable time sets $B_j$ with measure $\tau_j=2^{-j}$.
Let $X_j\subset\R^d$ have volume $\tau_j^{d/(2\alpha)}$ and define
\[
 F_N(t,x)=\sum_{j=1}^N\tau_j^{-\kappa}
                    \ind_{B_j}(t)\ind_{X_j}(x).
\]
For $t\in B_j$,
\[
 \|F_N(t)\|_{L_x^{R_i}}=\tau_j^{-1/P_i}.
\]
The geometric masses show
\[
 \|F_N\|_{L_t^{P_i,\infty}(\rd\mu;L_x^{R_i})}\le C_{P_i}
\]
uniformly in $N$. Let $1/P=(1-\theta)/P_0+\theta/P_1$ and
$1/R=(1-\theta)/R_0+\theta/R_1$.
The intersection bound for the real $K$-method gives
\[
 \|F_N\|_{(L_t^{P_0,\infty}L_x^{R_0},
                 L_t^{P_1,\infty}L_x^{R_1})_{\theta,P}}\le C
\]
uniformly in $N$, since its $K$-functional is bounded by $C\min(1,u)$.
But
\[
 \|F_N\|_{L_t^P(\rd\mu;L_x^R)}^P
 =\sum_{j=1}^N\tau_j\tau_j^{-1}=N,
\]
and the same growth occurs in $L_t^P L_x^{R,P}$, since a characteristic
function has Lorentz norm comparable to its measure to the power $1/R$.
Thus the necessary interpolation-space embedding is false.
This is a function-space counterexample, not a disproof of the
Schr\"odinger estimate itself.

The strong-time Lions--Peetre identity \eqref{real-L} used earlier is not
affected: the counterexample concerns weak outer-time endpoint spaces.
In particular the weak-outer formula asserted in
\cite[(2.2), Section 3.2]{FSW2025} cannot be used as a general interpolation
identity. Our proofs do not use this weak mixed-norm identity.

\paragraph{The spatial endpoint and other dispersions.}
The finite-$r$ Sobolev transfer cannot reach $r=\infty$ at its
critical potential order. It would require
$I_{d/r_0}:L_x^{r_0,q}\to L_x^\infty$, with $q>1$.
For any $u>1$, choose $1/u<b\le1$ and
\[
 h(x)=|x|^{-d/r_0}(\log(e/|x|))^{-b}
       \ind_{\{|x|<1/2\}}.
\]
Then $h\in L^{r_0,u}$, whereas $I_{d/r_0}h(0)=\infty$.
In fact the potential grows without bound as $x\to0$, so its
essential supremum is infinite. This disproves the needed spatial
embedding; it does not disprove an orthonormal Schr\"odinger estimate.

In the high-$q$ proof, $r=\infty$ makes the spatial potential exponent
$c=1$, so two neighboring Banach spatial exponents are unavailable.
At $OA$, the strict stationary-tail integrability in the $q=4$ lemma
is lost. At $q=\infty$, the asymmetric frequency-gap gain vanishes.
No uniform limit in those parameters is proved.
The strong critical estimate on $OD$ and general Frostman endpoint
forms on $OA$ therefore remain unresolved here. The necessity examples
in Section~\ref{sec-necessary} saturate $\beta=q/2$ under lower
regularity; they alone do not settle equality on $OD$.

For $a\ne2$, $s_A=(2-a)d\alpha/(2(d+\alpha))$.
Thus $U_a|D|^{-s_A}$ has no global unweighted energy bound, which is
the local-cell input used above. On a single annulus this costs a
frequency factor; the formal scaling cancels in the endpoint exponent,
but absolute summation of uniformly bounded weak Schatten annuli
loses the number of annuli. A global endpoint theorem for general $a$
or for the wave/Klein--Gordon flow does not follow from this computation.

\section{Proof of Theorem~\ref{thm-necessary}}\label{sec-necessary}

In this section, we prove Theorem~\ref{thm-necessary}. The constructions are modifications of those in \cite{BHLNS2019, BLN2021}. Since $E$ is bounded, the time variable cannot be rescaled. Instead, we work at a large frequency scale $N$, and the dimension $\alpha$ enters through the regularity exponent $s$ and, for the second condition, through \eqref{lower-regular}.

Let $(\phi,\Psi,a,\sigma)$ be one of \eqref{three-cases} and let $s=s_{a\alpha}(q,r)$; thus $s=s_\alpha(q,r)$ for the wave and Klein--Gordon equations. Then
\begin{equation}\label{d-2s}
d-2s=\frac{2d}r+\frac{2a\alpha}q.
\end{equation}
Since $\{f_j\}$ is orthonormal in $\dot H^s$ (resp. $H^s$) if and only if $\{\Psi(D)^sf_j\}$ is orthonormal in $L^2$, we assume that
\begin{align}\label{onf-strichartz-fractal}
	\Bigg\| \sum_{j\in J} n_j |U_\phi\Psi(D)^{-s} g_j|^2 \Bigg\|_{L_t^{q/2}(\mathrm d\mu; L_x^{r/2})} \le C_0 \| \{n_j\}_{j\in J}\|_{\ell^\beta}
\end{align}
holds for all orthonormal families $\{g_j\}$ in $L^2(\mathbb R^d)$ and all $\{n_j\}\in\ell^\beta$.

In all three cases, for $|\xi|\ge1$ we have
\begin{equation}\label{phase-bounds}
|\nabla\phi(\xi)|\lesssim|\xi|^{a-1},\qquad |\nabla^2\phi(\xi)|\lesssim|\xi|^{a-2},\qquad |\partial^\gamma_\xi\Psi(\xi)^{-s}|\lesssim_\gamma|\xi|^{-s-|\gamma|}.
\end{equation}

\subsection{Proof of Theorem~\ref{thm-necessary} (1)}
Let $N\ge1$ be large and $\rho=N^{1-a}$, so that $\rho\le N/8$ for $N$ large. Let $\{v_j\}_{j\in J}$ be a maximal $\rho$-separated subset of the annulus $\{N\le |\xi|\le2N\}$. Then $\#J\sim (N/\rho)^d=N^{ad}$. Let $\chi\in C_c^\infty(\mathbb R^d)$ be nonnegative, supported in $\{|\eta|\le1/4\}$, with $\|\chi\|_{L^2}=1$. Fix $t_0\in\supp\mu$ and define
\[
    \widehat {g_j}(\xi) := \rho^{-d/2}\,\e^{-it_0\phi(\xi)}\chi\big(\rho^{-1}(\xi - v_j)\big).
\]
The functions $\widehat{g_j}$ have disjoint supports, so $\{g_j\}_{j\in J}$ is an orthonormal family in $L^2(\mathbb R^d)$. Changing variables $\xi=v_j+\rho\eta$, we have
\[
\begin{aligned}
 U_\phi(t)\Psi(D)^{-s}g_j(x)
 &=\rho^{d/2}e^{i(x\cdot v_j+(t-t_0)\phi(v_j))}\mathcal I_j(t,x),\\
 \mathcal I_j(t,x)
 &=(2\pi)^{-d/2}\int e^{i\rho(x+(t-t_0)\nabla\phi(v_j))\cdot\eta}
       e^{i(t-t_0)\mathcal R_j(\eta)}\\
 &\qquad\cdot\Psi(v_j+\rho\eta)^{-s}\chi(\eta)\,\rd\eta,
\end{aligned}
\]
where $\mathcal R_j(\eta)=\phi(v_j+\rho\eta)-\phi(v_j)-\rho\nabla\phi(v_j)\cdot\eta$ satisfies $|\mathcal R_j(\eta)|\lesssim\rho^2N^{a-2}$ on $\supp\chi$ by \eqref{phase-bounds}. Moreover, $\Psi(v_j+\rho\eta)^{-s}=\Psi(v_j)^{-s}(1+O(\rho/N))$ and $\Psi(v_j)\sim N$.
Let $c_1>0$ be a small constant. If $|t-t_0|\le c_1$ and $|x|\le c_1\rho^{-1}$, then
\[
\begin{aligned}
 |\rho(x+(t-t_0)\nabla\phi(v_j))|
 &\lesssim c_1+c_1\rho N^{a-1}=2c_1,\\
 |t-t_0||\mathcal R_j|
 &\lesssim c_1\rho^2N^{a-2}=c_1N^{-a}\le c_1.
\end{aligned}
\]
Hence, if $c_1$ is small enough, the real part of the integrand is at least $\frac12\Psi(v_j)^{-s}\chi(\eta)$, and
\[
|U_\phi(t)\Psi(D)^{-s}g_j(x)|\gtrsim\rho^{\frac d2}N^{-s}\quad\text{for } |t-t_0|\le c_1,\ |x|\le c_1\rho^{-1}.
\]
Taking $n_j=1$ for $j\in J$, we obtain
\[
\begin{split}
\sum_{j\in J}|U_\phi(t)\Psi(D)^{-s}g_j(x)|^2
&\gtrsim\#J\,\rho^{d}N^{-2s}\\
&\sim N^{d-2s}\,\mathbf 1_{\{|t-t_0|\le c_1\}}(t)\,\mathbf 1_{\{|x|\le c_1\rho^{-1}\}}(x).
\end{split}
\]
Since $\|\mathbf 1_{\{|x|\le c_1\rho^{-1}\}}\|_{L_x^{r/2}}\simeq\rho^{-2d/r}$ and $\mu((t_0-c_1,t_0+c_1))>0$ does not depend on $N$, \eqref{onf-strichartz-fractal} and \eqref{d-2s} give
\[
N^{\frac{2d}r+\frac{2a\alpha}q}\rho^{-\frac{2d}r}=N^{\frac{2ad}r+\frac{2a\alpha}q}\lesssim(\#J)^{\frac1\beta}\sim N^{\frac{ad}\beta}.
\]
Letting $N\to\infty$, we get $\frac{2\alpha}q+\frac{2d}r\le\frac d\beta$, that is, $\beta\le\beta_{d/2}(q,r)$.

\subsection{Proof of Theorem~\ref{thm-necessary} (2)}
The assertion is automatic when $q=\infty$. If $q<\infty$ and
an estimate held with $\beta=\infty$, it would also hold with every
finite coefficient exponent, since $\|n\|_\infty\le\|n\|_\beta$.
We may therefore assume $2\le q<\infty$ and $q/2<\beta<\infty$
when deriving a contradiction.
We first consider the homogeneous cases $\phi(\xi)=|\xi|^a$ (including $a=1$). Let
\[
\widehat{g_0}(\xi)=\Big(\frac{a|\xi|^{a-d}}{2\pi|\mathbb S^{d-1}|}\Big)^{1/2}\mathbf 1_{[2\pi,4\pi)}(|\xi|^a).
\]
In polar coordinates with $u=|\xi|^a$, we have $|\widehat{g_0}(\xi)|^2\mathrm d\xi=\frac{\mathrm du}{2\pi}\,\mathrm d\omega/|\mathbb S^{d-1}|$, hence
\[
\langle U_\phi(j)g_0,U_\phi(k)g_0\rangle=\frac1{2\pi}\int_{2\pi}^{4\pi}\e^{i(j-k)u}\,\mathrm du=\delta_{jk},\qquad j,k\in\mathbb Z.
\]
For $\lambda\ge1$, let $g(x)=\lambda^{d/2}g_0(\lambda x)$ and $T=\lambda^{-a}$. By scaling, $U_\phi(t)\Psi(D)^{-s}g=\lambda^{\frac d2-s}\big(U_\phi(\lambda^at)\Psi(D)^{-s}g_0\big)(\lambda\,\cdot)$, so that $\{U_\phi(-jT)g\}_{j\in\mathbb Z}$ is orthonormal and
\begin{equation}\label{radial-lower}
 \begin{split}
 \|U_\phi(\tau)\Psi(D)^{-s}g\|_r
 &=\lambda^{d/2-d/r-s}
   \|U_\phi(\lambda^a\tau)\Psi(D)^{-s}g_0\|_r\\
 &\ge c_0\lambda^{a\alpha/q}=c_0T^{-\alpha/q},
 \qquad |\tau|\le T/2,
 \end{split}
\end{equation}
where $c_0=\min_{|\tau|\le1/2}\|U_\phi(\tau)\Psi(D)^{-s}g_0\|_{L^r}>0$. Here we used $s=s_{a\alpha}(q,r)$. The minimum is positive, since $\widehat{g_0}\in L^1\cap L^2$ is supported in an annulus, so $\tau\mapsto U_\phi(\tau)\Psi(D)^{-s}g_0$ is continuous with values in $L^r$, $2\le r\le\infty$, and never vanishes.

Let $J_n=[nT-T/2,nT+T/2)$, $n\in\mathbb Z$, and $m_n=\mu(J_n)$. Only finitely many $m_n$ are nonzero, and $\sum_nm_n=\mu(\mathbb R)=:\|\mu\|>0$. Let $g_n=U_\phi(-nT)g$ and $w_n\ge0$. For $t\in J_n$, discarding all terms except the one with index $n$ and using \eqref{radial-lower}, we get
\[
\Big\|\sum_jw_j|U_\phi(t)\Psi(D)^{-s}g_j|^2\Big\|_{L^{r/2}_x}\ge w_n\|U_\phi(t-nT)\Psi(D)^{-s}g\|_{L^r}^2\ge c_0^2\,w_nT^{-\frac{2\alpha}q}.
\]
Hence, with $n_j=w_j$, the left-hand side of \eqref{onf-strichartz-fractal} is at least
\[
c_0^2\,T^{-\frac{2\alpha}q}\Big(\sum_nw_n^{q/2}m_n\Big)^{2/q}.
\]
Suppose $\beta>q/2$ and let $\varkappa=\frac{2\beta}q>1$. Choosing $w_n=m_n^{2(\varkappa'-1)/q}$, so that $\sum_nw_n^{q/2}m_n=\sum_nm_n^{\varkappa'}$ and $\sum_nw_n^\beta=\sum_nm_n^{\varkappa'}$, the estimate \eqref{onf-strichartz-fractal} becomes
\[
 T^{-2\alpha/q}\left(\sum_nm_n^{\varkappa'}\right)^{2/q}
 \lesssim\left(\sum_nm_n^{\varkappa'}\right)^{1/\beta}.
\]
Since $2/q-1/\beta=2/(q\varkappa')$, raising this inequality to
$q/2$ yields
\begin{equation}\label{radial-contra}
T^{-\alpha}\|\{m_n\}\|_{\ell^{\varkappa'}}\lesssim1 .
\end{equation}
We now use \eqref{lower-regular}. If $m_n>0$, pick $t_n\in J_n\cap\supp\mu$. For indices $n$ of the same parity, the intervals $(t_n-T/2,t_n+T/2)$ are pairwise disjoint and, for $T\le2$, each has $\mu$-measure $\gtrsim T^\alpha$ by \eqref{lower-regular}. Hence $K:=\#\{n:m_n>0\}\lesssim\|\mu\|T^{-\alpha}$. By H\"older's inequality, $\|\mu\|=\sum_nm_n\le K^{1/\varkappa}\|\{m_n\}\|_{\ell^{\varkappa'}}$, so
\[
T^{-\alpha}\|\{m_n\}\|_{\ell^{\varkappa'}}\gtrsim T^{-\alpha}K^{-1/\varkappa}\gtrsim T^{-\alpha(1-1/\varkappa)}\to\infty\quad\text{as }T\to0,
\]
which contradicts \eqref{radial-contra}. Hence $\beta\le q/2$.

For Klein--Gordon evolution, take $\lambda\ge1$, $T=\lambda^{-1}$,
and choose $g$ with $\widehat g\ge0$ so that its Fourier probability measure is
\[
 |\widehat g(\xi)|^2\,\rd\xi
 =\frac{\rd u}{2\pi\lambda}\frac{\rd\omega}{|\mathbb S^{d-1}|},
 \quad u=\langle\xi\rangle\in[2\pi\lambda,4\pi\lambda),
 \quad\omega=\xi/|\xi|.
\]
The same integral over $u$ shows that $(U_\phi(-nT)g)_n$ is
orthonormal. Define
$\widehat g_\lambda(\eta)=\lambda^{d/2}\widehat g(\lambda\eta)$ and
\[
 \begin{aligned}
 F_\lambda(\tau,x)
 &=\lambda^{s-d/2}
   [U_\phi(\lambda^{-1}\tau)\langle D\rangle^{-s}g](\lambda^{-1}x)\\
 &=(2\pi)^{-d/2}\int
 e^{i(x\cdot\eta+\tau\sqrt{\lambda^{-2}+|\eta|^2})}
 (\lambda^{-2}+|\eta|^2)^{-s/2}\widehat g_\lambda(\eta)\,\rd\eta.
 \end{aligned}
\]
The profiles $\widehat g_\lambda$ converge in $L^1\cap L^2$ to
$\widehat g_0$ from the wave case $a=1$, and their supports lie in a
common compact annulus. On that annulus the phase and weight converge
uniformly, so
\[
 F_\lambda(\tau)\longrightarrow e^{i\tau|D|}|D|^{-s}g_0
 \quad\text{in }L^2\cap L^\infty,
 \quad |\tau|\le1/2,
\]
uniformly in $\tau$. Interpolation gives the same convergence in every
$L^r$, $2\le r\le\infty$. Consequently, for large $\lambda$,
\[
 \begin{aligned}
 \|U_\phi(\lambda^{-1}\tau)\langle D\rangle^{-s}g\|_r
 &=\lambda^{d/2-d/r-s}\|F_\lambda(\tau)\|_r\\
 &\ge(c_0/2)\lambda^{\alpha/q}
 =(c_0/2)T^{-\alpha/q},\qquad |\tau|\le1/2.
 \end{aligned}
\]
Here $c_0$ is the constant in \eqref{radial-lower} for $a=1$ with the same $s$. This is \eqref{radial-lower} with the Klein--Gordon regularity
$s=d/2-d/r-\alpha/q$. The time-cell and coefficient argument is unchanged.
\qed

\section{Density matrices}\label{sec-appli}
For fermionic systems the orthonormality assumption arises naturally from the Pauli exclusion principle. The dynamics of infinitely many fermions is conveniently described in terms of density matrices, that is, self-adjoint operators $\gamma$ on $L^2(\mathbb R^d)$, and the free evolution is given by
\begin{equation}\label{free-vN}
i\partial_t\gamma+[\phi(D),\gamma]=0,\qquad \gamma(t)=U_\phi(t)\gamma_0U_\phi(t)^*,\quad t\in\mathbb R,
\end{equation}
which is consistent with $i\partial_tu+\phi(D)u=0$ and
$U_\phi(t)=e^{it\phi(D)}$. For a locally trace-class operator, its
density is defined weakly by
\[
 \int h(x)\rho_\gamma(x)\,\rd x
 =\Tr(M_h\chi\gamma\chi),\qquad h\in C_c^\infty(\R^d),
\]
where $M_h$ is multiplication by $h$, and $\chi\in C_c^\infty$ is real and equals one on $\supp h$.
For kernels with a well-defined diagonal, this agrees with
$\rho_\gamma(x)=\gamma(x,x)$. We refer to
\cite{FLLS2014,FraSa2017,LS2014,LS2015} for the density-matrix formulation
of the orthonormal estimates.

Fix $1\le\beta<\infty$, and let $s\ge0$ in the inhomogeneous case and $s<d/2$ in the homogeneous case. Write
\[
 \Gamma=\Psi(D)^s\gamma_0\Psi(D)^s\in\mathfrak S^\beta,
 \qquad\|\gamma_0\|_{\mathfrak S^{\beta,s}}:=\|\Gamma\|_{\mathfrak S^\beta}.
\]
The compact self-adjoint operator $\Gamma$ has an expansion
$\Gamma=\sum_jn_j|g_j\rangle\langle g_j|$, with $(g_j)$ orthonormal
in $L^2$ and $\|n\|_{\ell^\beta}=\|\Gamma\|_{\mathfrak S^\beta}$.
For $\Psi(\xi)=\langle\xi\rangle$, the inverse weight is bounded,
so $\gamma_0=\Psi(D)^{-s}\Gamma\Psi(D)^{-s}$ is an operator on $L^2$.
In the homogeneous case $\Psi(\xi)=|\xi|$, the inverse weight
$|\xi|^{-s}$ may be unbounded, at the origin if $s>0$ and at infinity if
$s<0$.
We interpret this expression first as the Hermitian
form
\[
 \langle\gamma_0h,k\rangle
 =\langle\Gamma|D|^{-s}h,|D|^{-s}k\rangle
\]
on finite-annulus Fourier test vectors, and define its evolving density
by the completion below. It need not be a bounded $L^2$ operator.

Initially take $\Gamma$ of finite rank with annular Fourier support.
Then $f_j=\Psi(D)^{-s}g_j$ belong to $L^2$, and
\[
 \rho_{\gamma(t)}=\sum_jn_j|U_\phi(t)\Psi(D)^{-s}g_j|^2,
 \qquad
 |\rho_{\gamma(t)}|\le\sum_j|n_j||U_\phi(t)\Psi(D)^{-s}g_j|^2.
\]
Thus an orthonormal estimate
\begin{equation}\label{a-class-orthonormal}
 \left\|\sum_jn_j|U_\phi\Psi(D)^{-s}g_j|^2\right\|_
 {L_t^{q/2}(\rd\mu;L_x^{r/2})}
 \le C\|n\|_{\ell^\beta}
\end{equation}
gives
\begin{equation}\label{density-est}
 \|\rho_{\gamma(t)}\|_{L_t^{q/2}(\rd\mu;L_x^{r/2})}
 \le C\|\gamma_0\|_{\mathfrak S^{\beta,s}}.
\end{equation}
To pass to general $\Gamma$, choose self-adjoint Fourier contractions $Q_N$
supported in finite annuli and converging strongly to the identity. Form
\[
 \Gamma_N=Q_N\left(\sum_{j\le N}n_j|g_j\rangle\langle g_j|\right)Q_N.
\]
Then $\Gamma_N\to\Gamma$ in $\mathfrak S^\beta$. Since the vectors $Q_Ng_j$
need not be orthonormal, we diagonalize the finite-rank operators
$\Gamma_N$ and $\Gamma_N-\Gamma_M$ before applying
\eqref{a-class-orthonormal}. Applying the same
bound to $\Gamma_N-\Gamma_M$ shows that the densities converge in the
target norm and that the limit is independent of the approximants.
This defines the homogeneous weighted density when the initial state
is only a form.

The associated local trace limits can also be identified. Put
$P=q/2$, $R=r/2$, $A(t)=U_\phi(t)\Psi(D)^{-s}$ and
$\gamma_N(t)=A(t)\Gamma_NA(t)^*$. For
$H=\Gamma_N-\Gamma_M$ and a real $\chi\in C_c^\infty$,
since $-|H|\le H\le|H|$, we have
\[
 \|\chi A(t)H A(t)^*\chi\|_{\mathfrak S^1}
 \le\Tr(\chi A(t)|H|A(t)^*\chi).
\]
Diagonalizing $|H|$ and applying \eqref{a-class-orthonormal}, followed
by spatial H\"older and $L_t^P(\rd\mu)\subset L_t^1(\rd\mu)$, yields
\[
 \int\|\chi(\gamma_N(t)-\gamma_M(t))\chi\|_{\mathfrak S^1}\rd\mu(t)
 \le C\|\chi^2\|_{L^{R'}}\mu(\R)^{1-1/P}
       \|\Gamma_N-\Gamma_M\|_{\mathfrak S^\beta}.
\]
Thus the regularized states converge locally in
$L_t^1(\rd\mu;\mathfrak S^1)$. A countable compact cutoff exhaustion
gives compatible local trace-class limits for $\mu$-almost every $t$,
and trace continuity identifies their weak densities with the norm
completion. In the homogeneous case $s<d/2$ ensures
$|D|^{-s}U_\phi(t)^*h\in L^2$ for $h\in C_c^\infty$;
operator-norm convergence of $\Gamma_N$ then identifies these local
limits with the original Hermitian form on such tests. In the
inhomogeneous case the inverse weight is bounded and the same
identification is immediate. This does not require a bounded global
$L^2$ operator in the homogeneous case.

Conversely, testing \eqref{density-est} on
finite-rank annular states recovers the orthonormal estimate there,
and the scalar completion extends it to general Sobolev initial data.

Theorems~\ref{thm-frac-schro}, \ref{thm-wave}, \ref{thm-KG}, and
\ref{thm-quadratic-critical} therefore control densities sampled against
$\mu$. For finite main exponent and second index $u<\infty$, define
$\mathfrak S^{\beta,u,s}$ by the weighted norm
$\|\Psi(D)^s\gamma_0\Psi(D)^s\|_{\mathfrak S^{\beta,u}}$.
The same completion applies to the Lorentz data used here.
For the weak-time conclusion below, $P>1$ and finiteness of $\mu$
give $L_t^{P,\infty}(\rd\mu)\subset L_t^1(\rd\mu)$, so the local
trace-limit argument remains valid.
At $A_{d/2}$, the conclusion is
\[
 \|\rho_{\gamma(t)}\|_{L_t^{P_A}(\rd\mu;L_x^{R_A})}
 \le C \Fmu^{1/P_A}\|\gamma_0\|_{\mathfrak S^{P_A,1}},
 \qquad R_A=r_A/2,
\]
where $s=0$. Theorem~\ref{thm-quadratic-critical}~(2) gives
\eqref{density-est} with $\mathfrak S^{q/2,s}$, and
Theorem~\ref{thm-quadratic-critical}~(3) gives the corresponding weak-time
estimate in $L_t^{q/2,\infty}(\rd\mu;L_x^\infty)$ with $\mathfrak S^{q/2,1,s}$.

\begin{prop}\label{prop-dual-density}
Assume $\mu$ satisfies \eqref{frac_meas}, $2<q<\infty$, $2\le r\le\infty$, and that \eqref{a-class-orthonormal} holds. Then, for every $V\in L_t^{(q/2)'}(\mathrm d\mu; L_x^{(r/2)'})$,
    \begin{equation}\label{dual-density}
        \Bigg\|\int \Psi(D)^{-s}U_\phi(t)^*V(t,\cdot)U_\phi(t)\Psi(D)^{-s}\,\mathrm d\mu(t)\Bigg\|_{\mathfrak S^{\beta'}}\le C\|V\|_{L_t^{(q/2)'}(\mathrm d\mu; L_x^{(r/2)'})}.
    \end{equation}
\end{prop}
\begin{proof}
First take $V\ge0$. The scalar consequence of
\eqref{a-class-orthonormal} defines the positive bounded form
\[
 A=\int\Psi(D)^{-s}U_\phi(t)^*V(t)U_\phi(t)\Psi(D)^{-s}\,\rd\mu(t).
\]
For every positive finite-rank
$\Gamma=\sum_jn_j|g_j\rangle\langle g_j|$, $n_j\ge0$,
\begin{align*}
 \mathrm{Tr}(\Gamma A)
 &=\iint\sum_jn_j|U_\phi(t)\Psi(D)^{-s}g_j(x)|^2
       V(t,x)\,\rd x\rd\mu(t)\\
 &\le C\|\Gamma\|_{\mathfrak S^\beta}
       \|V\|_{L_t^{(q/2)'}(\rd\mu;L_x^{(r/2)'})}.
\end{align*}
Testing a positive operator against positive finite-rank operators
gives its Schatten $\beta'$ norm, also at $\beta'=1,\infty$ with
the usual interpretations. This proves \eqref{dual-density} for
positive $V$. For complex $V$, split the real and imaginary parts
into their positive and negative parts. Each has mixed norm at most
$\|V\|$, and the Schatten triangle inequality gives the stated bound
with at most a factor four.
\end{proof}

\begin{rem}
The Hartree equation for infinitely many particles, $i\partial_t\gamma+[\phi(D)+w*\rho_\gamma,\gamma]=0$, is naturally posed for $t$ in an interval, and our estimates only control the density at times in $E$. Whether fractal-time estimates of the form \eqref{density-est} can be used in a well-posedness theory, for instance to control the solution along a fractal set of times, is a question we do not address here.
\end{rem}

\section*{Acknowledgements}
This work was supported by the National Research Foundation of Korea
(RS-2024-00342160; S. Lee), (RS-2024-00461749, RS-2026-25484566;
J. B. Lee). G. Feng was supported by the China Scholorship Council (CSC).

The authors used an AI assistant (Claude, Anthropic) for language editing,
for checking computations and the consistency of the exposition, and for
suggestions on the organization of the manuscript. All mathematical
content was verified by the authors, who take full responsibility for the
paper.

\bibliographystyle{plain}

\end{document}